\documentclass[10pt, a4paper,reqno]{amsart}
\usepackage{amsthm}
\theoremstyle{}
{\theoremstyle{definition}
\newtheorem{dfn}{Definition}[section]}
\newtheorem{prop}[dfn]{Proposition}
\newtheorem{nota}[dfn]{Notation}
\newtheorem{thm}[dfn]{Theorem}
{\theoremstyle{definition}
\newtheorem{rem}[dfn]{Remark}}
\newtheorem{lem}[dfn]{Lemma}
\newtheorem{cor}[dfn]{Corollary}
\newtheorem{conj}[dfn]{Conjecture}

{\theoremstyle{definition}
\newtheorem{exa}[dfn]{Example}}

\usepackage[top=25truemm,bottom=25truemm,left=30truemm,right=30truemm]{geometry}
\usepackage{amsmath,amssymb,amscd,bm,graphicx,tikz-cd,upgreek,dsfont,amsfonts,tikz,enumitem}
\usetikzlibrary{cd}
\usepackage[all]{xy}
\usepackage[colorlinks=true,
    citecolor=green(pigment),
    linkcolor=alizarin,
    urlcolor=denim]{hyperref}
    
\definecolor{alizarin}{rgb}{0.82, 0.1, 0.26}
\definecolor{azure(colorwheel)}{rgb}{0.0, 0.5, 1.0}
\definecolor{blue(pigment)}{rgb}{0.2, 0.2, 0.6}
\definecolor{denim}{rgb}{0.08, 0.38, 0.74}
\definecolor{mint}{rgb}{0.24, 0.71, 0.54}
\definecolor{parisgreen}{rgb}{0.31, 0.78, 0.47}
\definecolor{persiangreen}{rgb}{0.0, 0.65, 0.58}
\definecolor{seagreen}{rgb}{0.18, 0.55, 0.34}
\definecolor{shamrockgreen}{rgb}{0.0, 0.62, 0.38}
\definecolor{green(pigment)}{rgb}{0.0, 0.65, 0.31}

\usepackage{setspace}
\newcommand{\marginparstretch}{0.6}
\let\oldmarginpar\marginpar
\renewcommand\marginpar[1]{\-\oldmarginpar[\framebox{\setstretch{\marginparstretch}\begin{minipage}{\marginparwidth}{\raggedleft\tiny #1}\end{minipage}}]{\framebox{\setstretch{\marginparstretch}\begin{minipage}{\marginparwidth}{\raggedright\tiny #1}\end{minipage}}}}

\usepackage{scalerel,stackengine}
\stackMath
\newcommand\reallywidecheck[1]{%
\savestack{\tmpbox}{\stretchto{%
   \scaleto{%
    \scalerel*[\widthof{\ensuremath{#1}}]{\kern-.6pt\bigwedge\kern-.6pt}%
    {\rule[-\textheight/2]{1ex}{\textheight}}%WIDTH-LIMITED BIG WEDGE
   }{\textheight}% 
}{0.5ex}}%
\stackon[1pt]{#1}{\scalebox{-1}{\tmpbox}}%
}
\def\ang<#1>{\langle #1 \rangle}
\def\bigang<#1>{\left\langle #1 \right\rangle}

\numberwithin{equation}{section}

\newcommand{\MF}{\operatorname{MF}}
\newcommand{\KQcoh}{\operatorname{KQcoh}}
\newcommand{\DMF}{\operatorname{DMF}}
\newcommand{\KMF}{\operatorname{KMF}}
\newcommand{\HMF}{\operatorname{HMF}}
\newcommand{\AMF}{\operatorname{AMF}}
\newcommand{\Dcoh}{\operatorname{Dcoh}}
\newcommand{\Dmod}{\operatorname{Dmod}}
\newcommand{\Kcoh}{\operatorname{Kcoh}}
\newcommand{\Kmod}{\operatorname{Kmod}}
\newcommand{\Acoh}{\operatorname{Acoh}}
\newcommand{\Amod}{\operatorname{Amod}}
\newcommand{\DcoQcoh}{\operatorname{D^{co}Qcoh}}
\newcommand{\AcoQcoh}{\operatorname{A^{co}Qcoh}}
\newcommand{\Dsg}{\operatorname{D^{sg}}}
\newcommand{\Sing}{\operatorname{Sing}}
\newcommand{\Spec}{\operatorname{Spec}}
\newcommand{\Supp}{\operatorname{Supp}}
\newcommand{\Hom}{\operatorname{Hom}}

\newcommand{\End}{\operatorname{End}}

\newcommand{\coh}{\operatorname{coh}}
\newcommand{\Qcoh}{\operatorname{Qcoh}}

\newcommand{\inj}{\operatorname{inj}}

\newcommand{\fmod}{\operatorname{mod}}

\newcommand{\proj}{\operatorname{proj}}
\newcommand{\Db}{\operatorname{D^b}}
\newcommand{\D}{\operatorname{D}}

\newcommand{\Kb}{\operatorname{K^b}}
\newcommand{\Perf}{\operatorname{Perf}}

\newcommand{\CM}{\operatorname{CM}}
\newcommand{\uCM}{\operatorname{\underline{CM}}}

\newcommand{\id}{\operatorname{id}}

\newcommand{\Image}{\operatorname{Im}}

\newcommand{\Ker}{\operatorname{Ker}}

\newcommand{\Irr}{\operatorname{Irr}}
\newcommand{\irr}{\operatorname{irr}}

\newcommand{\Res}{\operatorname{Res}}
\newcommand{\Ind}{\operatorname{Ind}}

\newcommand{\Fuk}{\operatorname{Fuk}}

\newcommand{\rk}{\operatorname{\sf rk}}

\newcommand{\Cone}{\operatorname{Cone}}

\newcommand{\cA}{\mathcal{A}}
\newcommand{\cB}{\mathcal{B}}
\newcommand{\cC}{\mathcal{C}}
\newcommand{\cD}{\mathcal{D}}
\newcommand{\cE}{\mathcal{E}}
\newcommand{\cF}{\mathcal{F}}

\newcommand{\cH}{\mathcal{H}}
\newcommand{\cI}{\mathcal{I}}

\newcommand{\cL}{\mathcal{L}}

\newcommand{\cO}{\mathcal{O}}

\newcommand{\cS}{\mathcal{S}}
\newcommand{\cT}{\mathcal{T}}

\newcommand{\cW}{\mathcal{W}}

\newcommand{\bA}{\mathbb{A}}
\newcommand{\bC}{\mathbb{C}}
\newcommand{\bG}{\mathbb{G}}

\newcommand{\bP}{\mathbb{P}}
\newcommand{\bQ}{\mathbb{Q}}

\newcommand{\bT}{\mathbb{T}}
\newcommand{\bZ}{\mathbb{Z}}

\newcommand{\K}{\mathfrak{K}}

\newcommand{\G}{{\rm G}}

\newcommand{\simto}{\xrightarrow{\sim}}
\newcommand{\hookto}{\hookrightarrow}

\renewcommand{\l}{\langle}
\renewcommand{\r}{\rangle}

\tikzset{
        DB/.style={circle,draw=black,circle,fill=white,inner sep=0pt, minimum size=4pt},
        DW/.style={circle,draw=black,fill=black,inner sep=0pt, minimum size=4pt},
        cvertex/.style={circle,draw=black,fill=white,inner sep=1pt,outer sep=3pt},
        vertex/.style={circle,fill=black,inner sep=1pt,outer sep=3pt},
        star/.style={circle,fill=yellow,inner sep=0.75pt,outer sep=0.75pt},
        tvertex/.style={inner sep=1pt,font=\scriptsize},
	pvertex/.style={circle,inner sep=1pt,outer sep=2pt,font=\scriptsize},
        gap/.style={inner sep=0.5pt,fill=white}
}

\makeatletter
\@namedef{subjclassname@2020}{%
  \textup{2020} Mathematics Subject Classification}
\makeatother

\begin{document}
\fontsize{9.7pt}{12pt}\selectfont

\title[]{Derived factorization categories of non-Thom--Sebastiani-type  sums of potentials}

\author[Y.~Hirano]{Yuki Hirano}
\address{Y.~Hirano, Department of Mathematics, Kyoto University, Kitashirakawa-Oiwake-cho, Sakyo-ku, Kyoto, 606-8502, Japan}\email{y.hirano@math.kyoto-u.ac.jp}
\author[G.~Ouchi]{Genki Ouchi}
\address{G.~Ouchi, Graduate School of Mathematics, Nagoya University, 
Furocho, Chikusaku, Nagoya,  464-8602, Japan}\email{genki.ouchi@math.nagoya-u.ac.jp}
\date{}

%critical locus

\begin{abstract} 
%For  a smooth quasi-projective variety $X$ with a reductive group action and semi-invariant regular function $W:X\to {\mathbb A}^1$, we describe a  semi-orthogonal decomposition of the  derived factorization categories  associated to the sum  $W+F:X\times {\mathbb A}^m\to{\mathbb A}^1$ of $W$ and a certain semi-invariant regular function $F:X\times{\mathbb A}^m\to{\mathbb A}^1$.  Our situation does not require that $W+F$ is Thom--Sebastiani sum.  
%As an application, we prove that the homotopy category of maximally graded matrix factorizations of an invertible polynomial $f$ of chain type  has a full exceptional collection ${\mathcal E}$, and we also show that the length of  ${\mathcal E}$ is the Milnor number of the Berglund--H\"ubsch transpose $\widetilde{f}$ of $f$. We also give a semi-orthogonal decomposition that is a generalization of a special version of Kuznetsov--Perry's semi-orthogonal decomposition of a certain equivariant category associated to the derived category of a cyclic cover of  a  variety.

We first prove semi-orthogonal decompositions of derived factorization categories arising from sums of potentials of gauged Landau-Ginzburg models, where the sums are not necessarily Thom--Sebastiani type.  We then apply the result to the category $\HMF^{L_f}(f)$  of maximally graded matrix factorizations of  an invertible polynomial $f$ of chain type, and explicitly construct a full strong exceptional collection $E_1,\hdots,E_{\mu}$ in $\HMF^{L_f}(f)$  whose  length $\mu$ is the Milnor number of the Berglund--H\"ubsch transpose $\widetilde{f}$ of $f$. This proves a conjecture, which postulates that for an invertible polynomial $f$ the category $\HMF^{L_f}(f)$ admits a tilting object, in the case when $f$ is a chain polynomial. Moreover, by careful analysis of morphisms between the exceptional objects $E_i$, we  explicitly determine the quiver with relations $(Q,I)$ which represents the endomorphism  ring of the associated tilting object $\oplus_{i=1}^{\mu}E_i$ in $\HMF^{L_f}(f)$, and in particular we obtain an equivalence $\HMF^{L_f}(f)\cong \Db(\fmod kQ/I)$.
\end{abstract}

\maketitle{}

\section{Introduction}
\subsection{Backgrounds}
Let  $f\in S^n:=\bC[x_1, \cdot \cdot \cdot, x_n]$ be a quasi-homogeneous  polynomial.  The {\it Milnor number} of $f$, denoted by $\mu(f)$, is the dimension of the {\it Jacobian ring} $S^n/\l{\partial f}/{\partial x_1},\dots,{\partial f}/{\partial x_n}\r$ as a $\bC$-vector space.
We say that  $f \in S^n$ is an {\it invertible polynomial} if  the following conditions hold:
\begin{itemize}
\item[(i)]There is a $n\times n$-matrix $E=(E_{i,j})$ whose entries $\{E_{i,j}\}$ are non-negative integers such that $E$ is invertible over $\bQ$, and that $f$ is of the form
\[
f(x_1,\hdots,x_n)=\sum_{i=1}^n\Bigl(\prod_{j=1}^n x_{j}^{E_{i,j}}\Bigr).
\] 
\item[(ii)]The {\it Berglund--H\"ubsch transpose} $\widetilde{f}$ of $f$, which is defined by 
\[
\widetilde{f}(x_1,\hdots,x_n)=\sum_{i=1}^n\Bigl(\prod_{j=1}^n x_{j}^{E_{j,i}}\Bigr),
\]
is also quasi-homogeneous.
\item[(iii)]  We have $1\leq \mu(f) <\infty$ and $1\leq \mu(\widetilde{f}) <\infty$.
\end{itemize}
Invertible polynomials are studied in \cite{bh,bh2} to construct pairs of topological mirror Calabi-Yau manifolds as a generalization of the Green--Plesser construction \cite{gp}.
By Kreuzer--Skarke \cite{ks},  an invertible polynomial is the Thom--Sebastiani sums of several invertible polynomials of the following two types:

\vspace{3mm}
({\it chain type}) \hspace{2mm}$x_1^{a_1}+x_1x_2^{a_2}+\cdots+x_{n-2}x_{n-1}^{a_{n-1}}+x_{n-1}x_n^{a_n}$ \hspace{5mm}($a_1\geq 2$ and $a_n\geq 2$)

\vspace{2mm}
({\it loop type}) \hspace{3.9mm}$x_1^{a_1}x_2+x_2^{a_2}x_3+\cdots+x_{n-1}^{a_{n-1}}x_n+x_n^{a_n}x_1$ \hspace{9.1mm}($n\geq 2$)
\vspace{3mm}\\
Here, for invertible polynomials $f \in S^n$ and $g \in S^m$, the {\it Thom--Sebastiani sum} $f\boxplus g$ of $f$ and $g$ is defined by $f\boxplus g:=f\otimes 1+1\otimes g\in S^n\otimes_{\,\bC}S^m\cong S^{n+m}$.
%\in \bC[x_1, \cdot \cdot \cdot, x_m, y_1, \cdot \cdot \cdot, y_n]
We temporarily say that an invertible polynomial $f$ is {\it indecomposable} if it can not be decomposed into the Thom--Sebastiani sum of two invertible polynomials, or equivalently if it is of chain type or loop type.  When  $n=1$, an indecomposable invertible polynomial is of the form $x_1^{a_1}$ for some $a_1\geq 2$, and this polynomial is said to be  of {\it Fermat type}.
For an invertible polynomial $f$, we can associate two triangulated categories in the context of the study of singularities. One is the {\it homotopy category $\HMF^{L_f}_{S^n}(f)$ of maximally graded matrix factorizations} of $f$. Here the {\it maximal grading} $L_f$ of $f$ is the abelian group defined by
\[
L_f:=\left.{\left(\left(\bigoplus_{i=1}^n\bZ \vec{x}_i\right)\oplus \bZ \vec{f}\right)} \middle/ {\left\l\left.{\vec{f}-\sum_{j=1}^nE_{i,j}\vec{x}_j}\,\right| 1\leq i\leq n\right\r}\right..
\]
This category $\HMF^{L_f}_{S^n}(f)$ is equivalent to the maximally graded singularity category $\D^{L_f}_{\mathrm{sg}}(S^n/f)$ of the hypersurface $S^n/f$   and the stable category  $\underline{\CM}^{L_f}(S^n/f)$ of maximally graded Cohen-Macaulay modules over  $S^n/f$. Another one is the derived category $\Db\Fuk^{\to}(f)$ of a directed Fukaya category of $f$. The category $\Db\Fuk^{\to}(f)$ is the categorification of the Milnor lattice of the isolated hypersurface singularity $f$.
Takahashi proposed the following conjecture as a version of homological mirror symmetry for isolated hypersurface singularities (see also \cite{Tak2,et2,et}).

\begin{conj}[{\cite[Conjecture 21]{Tak3}}]\label{HMS}
For an invertible polynomial $f$, there is a finite acyclic quiver $Q$ with admissible relations $I$ such that we have equivalences 
\[\HMF^{L_f}_{S^n}(f) \simeq \Db(\mathrm{mod}\ \bC Q/I) \simeq \Db\Fuk^{\to}(\widetilde{f}).\] 
\end{conj}
For ADE polynomials, the first equivalence in Conjecture \ref{HMS} is proved by Kajiura--Saito--Takahashi \cite{kst1} and the second equivalence is proved by Seidel \cite{seidel}. 
 Conjecture \ref{HMS} for Brieskorn-Pham singularities  and the Thom--Sebastiani sum of ADE polynomials of type A or D  follow  from Futaki--Ueda  \cite{fu1,fu2}.
Recently, Habermann--Smith \cite{hs} and Habermann \cite{hab} proved Conjecture \ref{HMS} for the $n=2$ cases  by constructing tilting objects on the both sides.  Conjecture \ref{HMS} leads to the following conjecture.

%Homological mirror symmetry for isolated hypersurface singularities predicts the equivalence $\HMF^{L_f}(f) \simeq \Db\Fuk^{\to}(\widetilde{f})$. 
%By construction, the triangulated category $\Db\Fuk^{\to}(\widetilde{f})$ has a full exceptional collection $\cE$, and the length of  $\cE$ is equal to the Milnor number $\mu(\widetilde{f})$ of $\widetilde{f}$. The existence of the first equivalence in Conjecture \ref{HMS} is equivalent to the existence of full strong exceptional collection of the category $\HMF^{L_f}(f)$.
%Conjecture \ref{HMS} implies the following.

\begin{conj}\label{fsec}
Let $f\in S^n$ be an invertible polynomial. The category $\HMF^{L_f}_{S^n}(f)$  has a full strong exceptional collection of length $\mu(\widetilde{f})$.
\end{conj}
This conjecture is proved by Kajiura--Saito--Takahashi  \cite{kst1,kst2} for special $n=3$ cases, namely ADE singularities and $14$ exceptional unimodular singularities, and recently Kravets \cite{krav} proved Conjecture \ref{fsec} for the  cases when $n\leq3$.
Conjecture \ref{fsec} implies the following. 

\begin{conj}[{\cite[Conjecture 1.3]{lu}}]\label{tilting conj}
Let $f\in S^n$ be an invertible polynomial. The category $\HMF^{L_f}_{S^n}(f)$  has a tilting object.
\end{conj}

Proving  Conjecture \ref{tilting conj} is  important in Lekili--Ueda's work \cite[Theorem 1.6]{lu} and the study of Cohen-Macaulay representations of graded Gorenstein rings \cite[Problem 3.4]{iyama}.  Conjecture  \ref{fsec} can be weakened to the following.

\begin{conj}\label{exc}
Let $f\in S$ be an invertible polynomial. The category $\HMF^{L_f}_{S}(f)$  has a full exceptional collection of length $\mu(\widetilde{f})$.
\end{conj}
In an unpublished work by Yoko Hirano and Takahashi, they proved that Conjecture \ref{exc} is true when $n\leq3$, and Conjecture \ref{fsec} is true if $n\leq3$ and $f$ is an invertible polynomial  of chain type. 
For an invertible polynomial $f$ of chain type, Aramaki and Takahashi \cite{at} constructed the full exceptional collection of the category $\HMF^{L_f}_{S^n}(f)$ as the algebraic counterpart of the Orlik-Randell conjecture, which is proved by Varolgunes \cite{var}  recently. After the previous version of our paper appeared,  Conjecture \ref{exc} was  proved by Favero--Kaplan--Kelly \cite{fkk}. 

By recent progress of the derived Morita theory for factorizations by \cite{bfk}, Conjecture \ref{fsec} reduces to the case of indecomposable invertible polynomials, i.e. of chain type or  loop type.
In this paper, we prove Conjecture \ref{fsec} for  invertible polynomials of chain type.

%The existence of full strong exceptional collections on the homotopy categories of equivariant matrix factorizations was first proved for indecomposable invertible polynomials of Fermat type by A.~Takahashi \cite{Tak}. In \cite{kst1} and \cite{kst2}, Kajiura--Saito--Takahashi constructed full strong exceptional collections on the homotopy categories of maximally graded matrix factorizations of  ADE polynomials in three variables. 
%In an unpublished work by Yoko Hirano and A.~Takahashi, they proved that the conjecture is true when $n\leq3$, and that there exists a full strong exceptional collection on $\HMF^{L_f}_S(f)$ if $n\leq3$ and $f$ is an invertible polynomial  of chain type. 
%Kravets \cite{krav} proved Conjecture \ref{fsec} for invertible polynomials in three veriables.

\subsection{Main result}

 For each positive integer $i\in\bZ_{\geq1}$, choose a positive integer $a_i\in\bZ_{\geq 1}$. Then for $n\in \bZ_{\geq1}$, we have the associated  polynomial 
 \[f_n:=x_1^{a_1}+x_1x_2^{a_2}+\cdots +x_{n-1}x_n^{a_n},\]
and we set $f_0:=0$. This polynomial $f_n$ is an invertible polynomial of chain type only when $a_n>1$, but we also discuss the case when $a_n=1$. Even when $a_n=1$, we can define the maximal grading $L_n:=L_{f_n}$ of  $f_n$ similarly. More precisely, the  maximal grading $L_n$ is defined by
\begin{equation}\label{max intro}
L_n:=\left.{\left(\left(\bigoplus_{i=1}^n\bZ \vec{x}_i\right)\oplus \bZ \vec{f}_n\right)} \middle/ {\left\l\left.{\vec{f}_n-a_1\vec{x}_1, \,\,\vec{f}_n-\vec{x}_{i-1}-a_i\vec{x}_i}\,\right| 2\leq i\leq n\right\r}\right.,
\end{equation}
and we set  $L_0:=\bZ\vec{f_0}\cong \bZ$. Note that $L_1\cong \bZ$ and $\rk(L_n)=1$ for any $n\geq0$.  
The $n$-dimensional polynomial ring $S^n=\bC[x_1,\hdots,x_n]$ has a natural $L_n$-grading such that $\deg(x_i):=\vec{x}_i$ for any $1\leq i\leq n$ and $\deg(f)=\vec{f}$. For $m\geq-1$, we set 
\[\cC_m:=\HMF_{S^m}^{L_m}(f_{m}),\]
where $\cC_{-1}$ is the zero category. Then we have a natural equivalence $\cC_0=\HMF^{\bZ}_{\bC}(0)\cong \Db(\coh \Spec \bC)$.
%For $l \in L_n$, the category $\cC_n$ has the grading shift functor $(l): \cC_n \simto \cC_n$. 
For each $n$, we define the functors
\begin{eqnarray*}
 \uppsi_i &:& \cC_{n-1} \to \cC_n \hspace{5mm} (0 \leq i \leq a_n-2) \\
\upphi_j &:& \cC_{n-2} \to \cC_n \hspace{5mm} (0 \leq j \leq a_{n-1}-1) 
\end{eqnarray*}
explicitly  (see Section \ref{explicit matrix} for the details). 
Using these functors, we inductively construct  the full strong exceptional collection $\cE_n$ of the category $\cC_n$ as follows: 

For  a sequence  $\cE=(E_1,\hdots,E_r)$ of objects in $\cC_{n-1}$ and for each $0\leq i\leq a_n-2$, we define a sequence $\uppsi_i\cE$ of objects in $\cC_n$ by 
\[
\uppsi_i\cE:=(\uppsi_iE_1,\hdots, \uppsi_iE_r).
\]
Similarly, for a sequence  $\cF=(F_1,\hdots,F_s)$ of objects in $\cC_{n-2}$ and for $0\leq j\leq a_{n-1}-1$, we define a sequence $\upphi_j\cF$ of objects in $\cC_n$ by
\[
\upphi_j\cF:=(\upphi_jF_1,\hdots,\upphi_jF_s).
\]
Moreover, we set 
\[
\cE^{-1}:=\emptyset\mbox{\hspace{5mm}and\hspace{5mm}}
\cE^0:=\{(0 \to \bC \to 0)\}.
\]
Then we inductively define the sequence $\cE^n$ by
\[
\cE^{n}:=\left(\uppsi_0\cE^{n-1},\hdots,\uppsi_{(a_n-2)}\cE^{n-1},\upphi_0\cE^{n-2},\hdots,\upphi_{(a_{n-1}-1)}\cE^{n-2}\right).
\]
The following is our main result in this paper.

\begin{thm}[Theorem \ref{strong thm}]\label{strong thm intro}
For any $n\geq1$, the sequence $\cE^n$ is a full strong exceptional collection in $\cC_n$, and  if $a_n\geq2$ the  length of $\cE^n$ is equal to the Milnor number $\mu(\widetilde{f}_n)$ of $\widetilde{f}_n$. In particular, the category $\cC_n$ has a tilting object.
\end{thm}

To prove Conjecture \ref{HMS}, it is important to determine a quiver with relation $(Q,I)$ such that $\cC_n\cong \Db(\fmod \bC Q/I)$. By a careful analysis of morphisms between exceptional objects in $\cE^n$, we also determine  a quiver with relations $(Q^n,I^n)$ such that the endomorphism ring $\End_{\cC_n}(T_n)$ of the associated tilting object $T_n:=\bigoplus_{E\in \cE^n}E$ is isomorphic to $\bC Q^n/I^n$, and in particular we obtain  an equivalence \[\cC_n\cong \Db(\fmod\bC Q^n/I^n).\]
 We expect that this equivalence can be applied to the study of Conjecture \ref{HMS} for chain polynomials.
 
\subsection{Contribution to Cohen-Macaulay representations}
 
For  an abelian group $L$ of rank one and an $L$-graded Gorenstein ring $R$, the stable category $\underline{\CM}^L(R)$ of $L$-graded Cohen-Macaulay modules over $R$ is one of the  important invariant in the study of Cohen-Macaulay representations of $R$, and  the existence of a tilting object in the category $\underline{\CM}^L(R)$ is widely studied.  For example, if $R=\bigoplus_{i\geq0}R_i$ has a positive grading of $L=\bZ$ with $R_0$ a filed, $\underline{\CM}^{\bZ}(R)$ has a tilting object  either $\dim R=0$  by Yamaura \cite{yamaura} or $\dim R=1$ and $R$ is reduced by Buchweitz--Iyama--Yamaura \cite{biy}. When $\dim R\geq 2$, $\underline{\CM}^{\bZ}(R)$ does not have a tilting object in general, and there is no  general  existence result for a tilting object in $\underline{\CM}^{\bZ}(R)$. 
 Therefore it is important to find  classes of $(R,L)$ containing higher dimensional $R$ such that  $\uCM^L(R)$ has a tilting object, and it is  proved by Herschend--Iyama--Minamoto--Oppermann  that Geigle--Lenzing complete intersection, which are a generalization of weighted projective line \cite{gl}, are one  such class \cite{himo}.
 
 Our result gives a new class   $(R,L)$ such that $\uCM^L(R)$ has a tilting object. More precisely, with the same notation  as in the previous subsection, Theorem \ref{strong thm intro} implies that the category \[\uCM^{L_n}(S^n/f_n)\] has a tilting object. Moreover, the algebras $\bC Q^n/I^n$ are  new examples of finite dimensional algebras  whose derived categories are fractional Calabi-Yau, which might be of independent interest.

 \subsection{Outline of the proof of Theorem \ref{strong thm intro}}
 %We sketch the outline of the proof of Theorem \ref{strong thm} using the same notation as above. 
   First, we  see that $\cC_n$ is equivalent to the derived factorization category, which is a geometric interpretation of $\cC_n$, as follows:
 To the invertible polynomial  $f_n$  of chain type, we associate the {\it maximal symmetry group}
 \begin{equation}\label{max sym intro}
 G_n:=\Bigl\{(\lambda_1,\hdots,\lambda_n)\in(\bG_m)^n \, \Big| \, \lambda_1^{a_1}=\lambda_1\lambda_2^{a_2}=\cdots=\lambda_{n-1}\lambda_n^{a_n}\Bigr\}
\end{equation}
 of $f_n$, which naturally acts on $\bA_n$. Then, the chain polynomial  $f_n$ is a semi-invariant regular function on $\bA^n$ with respect to the $G_n$-action and a character 
 \[
 \chi_{f_n}:G_n\to \bG_m; \,\,(\lambda_1,\hdots,\lambda_n)\mapsto \lambda_1^{a_1}
 \]
 of $G_n$, and  these define  a {\it gauged Landau--Ginzburg model} $(\bA^n,\chi_{f_n},f_n)^{G_n}$. Then we have its {\it derived factorization category}
 \[
\cD^{\coh}_n:=\Dcoh_{G_n}(\bA^n,\chi_{f_n},f_n),
 \]
 and it is  standard  that the category $\cD_n^{\coh}$ is  equivalent to $\cC_n$. 
 
 Next, we prove a semi-orthogonal decomposition of $\cC_n$, which describes  relationships among $\cC_{n-2}$, $\cC_{n-1}$ and $\cC_n$. For this, we show that there are  fully faithful functors 
\begin{eqnarray*}
 \Uppsi_i &:& \cD^{\coh}_{n-1} \to \cD^{\coh}_n \hspace{5mm} (-a_n+2 \leq i \leq 0) \\
\Upphi_j &:& \cD^{\coh}_{n-2} \to \cD^{\coh}_n \hspace{5mm} (-a_{n-1}+1 \leq j \leq 0) 
\end{eqnarray*}
 and a semi-orthogonal decomposition 
 \begin{equation}\label{sod intro}
\cD^{\coh}_n={\Bigg \langle}{\sf Im}(\Uppsi_0),\hdots,{\sf Im}(\Uppsi_{(-a_{n}+2)}),\bigoplus_{j=0}^{{-a_{n-1}+1}}{\sf Im}(\Upphi_{j}){\Bigg \rangle}.
\end{equation}
This semi-orthogonal decomposition is deduced from  a more general semi-orthogonal decomposition associated to the sum of polynomials that are not necessarily Thom--Sebastiani sums. We prove this general semi-orthogonal decompositions by modifying  arguments  appearing in the proofs of  \cite[Corollary 3.4]{degree d} and \cite[Theorem 5.3]{H2}, which are analogous to results of  \cite[Theorem 40]{orlov3}.
%of combining semi-orthogonal decompositions via variations of GIT quotients \cite{vgit,h-l,segal} and Kn\"orrer periodicity \cite{H2, isik,shipman}.
By natural equivalences $\cC_m\cong \cD^{\coh}_m$ for $n-2\leq m\leq n$, we obtain corresponding  fully faithful functors
 $\Uppsi_i:\cC_{n-1}\to \cC_n$ and $\Upphi_j:\cC_{n-2}\to \cC_n$ and the corresponding semi-orthogonal decomposition 
  \[
\cC_n={\Bigg \langle}{\sf Im}(\Uppsi_0),\hdots,{\sf Im}(\Uppsi_{(-a_{n}+2)}),\bigoplus_{j=0}^{{-a_{n-1}+1}}{\sf Im}(\Upphi_{j}){\Bigg \rangle}.
\]
 
 Finally, we describe the fully faithful functors $\Uppsi_i$ and $\Upphi_j$ to obtain the following isomorphisms 
 \[
 \Uppsi_{i}\cong \uppsi_{-i}[i]
\mbox{\hspace{5mm} and \hspace{5mm}}
\Upphi_j\cong \upphi_{-j}[-a_n+1+j].
\]
By these explicit descriptions and the above semi-orthogonal decomposition,  we see that the sequence $\cE^n$ is a full exceptional collection, and  we prove  strongness by  utilizing  Serre duality and  certain exact triangles  (see Lemma \ref{psi phi tri}) in $\cC_n$.
%%%%%%%%%%%%%%%%%%%%%%%%%%%%%%%%%%%%%%%%%%%

\subsection{Organization of the paper} This paper is organized as follows: In Section 2, we provide the basic definitions and properties  of equivariant coherent sheaves, graded modules,  derived factorization categories, and homotopy categories of graded matrix factorizations. In Section 3, we prove  semi-orthogonal decompositions of derived factorization categories associated to the sums of certain functions, which are generalizations of \eqref{sod intro}. In Section 4, we prove Theorem \ref{strong thm intro} and compute the quiver with relations associated to the full strong exceptional collection. In Section 5, we discuss further applications of the general version of the semi-orthogonal decomposition proved in Section 3, and we prove a generalization of Kuznetsov--Perry's semi-orthogonal decompositions of homotopy categories of graded matrix factorizations of the Thom--Sebastiani sum $f\boxplus t^N$ of a quasi-homogeneous polynomial $f$ and a monomial $t^N$ in one variable $t$. In the appendix, we provide a brief summary of equivariant categories, which is necessary in Section 5. 
\subsection{Notation and convention} 

\begin{itemize}
\item Unless stated otherwise, all categories and stacks are over an algebraically closed field $k$ of characteristic zero.

\vspace{2mm}
\item For a functor $F:\cA\to \cB$, we denote by ${\sf Im}(F)$ the essential image of $F$.

\vspace{2mm}
\item For an integer $l\in\bZ$,  $\chi_l:\bG_m\to\bG_m$ denotes the character of $\bG_m$ defined by $\chi_l(a):=a^l$ for $a\in \bG_m$.

\vspace{2mm}
\item For a character $\chi:G\to \bG_m$ of an algebraic group $G$, we denote by $\cO(\chi)$ the $G$-equivariant invertible sheaf on $X$ associated to $\chi$. For a $G$-equivariant quasi-coherent sheaf $F$ on $X$, we set $F(\chi):=F\otimes \cO(\chi)$.

\vspace{2mm}
\item For a dg-category $\cA$, we denote by $[\cA]$ its homotopy category.
\end{itemize}

\subsection{Acknowledgements} 
The authors would like to thank Atsushi Takahashi for giving valuable comments on a draft version of this article, and for telling us his unpublished work and deep insight into the homological LG mirror symmetry. They also thank the referee for carefully reading the paper, and giving valuable comments and suggestions.
Y.H. is supported by JSPS KAKENHI Grant Number 17H06783 and 19K14502. G.O. is supported by 
Interdisciplinary Theoretical and Mathematical Sciences Program (iTHEMS) in RIKEN and JSPS KAKENHI Grant Number 19K14520. A part of the paper was written during visits of both authors to the Max-Planck-Institute for Mathematics in Bonn. The authors would like to thank the Max-Planck-Institute for Mathematics for their hospitality and support.
\section{Preliminaries}

\subsection{Equivariant coherent sheaves}
We recall the basics of equivariant sheaves. Let $G$ be an algebraic group acting on a scheme $X$. 
A  {\it quasi-coherent   $G$-equivariant sheaf} on $X$ is a pair $(F, \theta)$ of a quasi-coherent  sheaf $F$ and an isomorphism $\theta:\pi^*F\xrightarrow{\sim}\sigma^*F$ such that 
$$\iota^*\theta={\rm id}_{F} \hspace{4mm}\rm{and} \hspace{5mm}\big{(}({\rm id}_G\times\sigma)\circ({\it s}\times{\rm id}_X)\big{)}^*\theta\circ({\rm id}_G\times\pi)^*\theta=({\it m}\times{\rm id}_X)^*\theta,$$
where $m:G\times G\rightarrow G$ is the multiplication and $s:G\times G\rightarrow G\times G$ is the switch of two factors. Note that if $(F,\theta)$ is a quasi-coherent $G$-equivariant sheaf, for any $g\in G$ the equivariant structure $\theta$ defines an isomorphism
\[
\theta_g:F\cong \sigma_g^*F,
\]
where $\sigma_g:X\to X$ is  the action by $g\in G$.
We call  $(F, \theta)$  {\it coherent} (resp. {\it locally free}, {\it injective}) if the sheaf $F$ is coherent (resp. locally free, injective).

We denote by $\Qcoh_GX$ (resp. $\coh_GX$) the category of quasi-coherent (resp. coherent) $G$-equivariant sheaves on $X$ whose morphisms are $G$-invariant morphisms. Here, a {\it $G$-invariant morphism} $\varphi:(F_1, \theta_1)\rightarrow(F_2, \theta_2)$ of equivariant sheaves is a morphism of sheaves $\varphi:F_1\rightarrow F_2$ that commutes with $\theta_i$, i.e. $\sigma^*\varphi\circ\theta_1=\theta_2\circ\pi^*\varphi$.

\begin{dfn}
Let $G,H$ be  algebraic groups and $f:H\to G$ a morphism of algebraic groups.
Let $X$ be a $G$-variety, $Y$ an $H$-variety and $\varphi:Y\to X$ an  $f$-equivariant morphism, i.e. $\varphi(h\cdot y)=f(h)\cdot \varphi(y)$ for any $h\in H$ and $y\in Y$. Then we define the pull-back functor
\[
\varphi_f^*:\Qcoh_GX\to \Qcoh_HY
\]
by $\varphi_f^*(F,\theta):=\bigl(\varphi^*F,(f\times \varphi)^*\theta\bigr)$. If $H=G$ and $f={\rm id}_G$,  we denote the pull-back by $\varphi^*$ and define the push-forward
\[ 
\varphi_*:\Qcoh_GY\to \Qcoh_GX
\]
by $\varphi_*(F,\theta):=(\varphi_*F, ({\rm id}_G\times \varphi)_*\theta)$.
If $H$ is a closed subgroup of $G$ and $f$ is the inclusion, the pull-back ${\rm id}_f^*:\Qcoh_GX\to \Qcoh_HX$ is called the {\it restriction functor}, and is denoted by $\Res_H^G$. We write \[\Res ^G:\Qcoh_GX\to \Qcoh X\] for the restriction functor $\Res^G_{\{1\}}$.

\end{dfn}

\begin{dfn}
Let $G$ be an algebraic group acting on a variety $X$. Assume that a restricted action $H\times X\to X$ of a closed normal subgroup $H\subseteq G$ is trivial. 

(1) We  define a functor 
\[
(-)^H:\coh_GX\to \coh_{G/H}X.
\]
as follows: For $(F,\theta)\in \coh_GX$, we define a subsheaf $F^H\subset F$ by the following local sections on any open subspace $U\subseteq X$;
\[
F^H(U):=\{ x\in F(U)\mid \theta_h(x)=x,\mbox{ for }\forall h\in H\}.
\]
Then since $\theta:\pi^*F\simto\sigma^*F$ maps the subsheaf $\pi^*(F^H)\subset\pi^*F$ to the subsheaf $\sigma^*(F^H)\subseteq \sigma^*F$, the pair $(F^H,\theta|_{\sigma^*(F^H)})$ is a $G$-equivariant sheaf that naturally descends to a $G/H$-equivariant sheaf.

(2) For a character $\chi:H\to \bG_m$  and $(F,\theta)\in \coh_GX$, we define a subsheaf $F_{\chi}\subseteq F$ by the following local sections on any open subspace $U\subseteq X$;
\[
F_{\chi}(U):=\{ x\in F(U)\mid \theta_h(x)=\chi(h)x,\mbox{ for }\forall h\in H\}.
\]
Then since $\theta$ preserves $F_{\chi}$, we have a $G$-equivariant sheaf $(F_{\chi},\theta|_{\pi^*F_{\chi}})\in \coh_GX$. We call $(F,\theta)$ is of weight $\chi$ if $F=F_{\chi}$, and  we define a subcategory $(\coh_GX)_{\chi}\subset\coh_GX$ consisting of equivariant sheaves of weights $\chi$.  Then we have a functor
\[
(-)_{\chi}:\coh_GX\to (\coh_GX)_{\chi}.
\]
\end{dfn}

\vspace{3mm}
For later use, we provide a few fundamental lemmas.  For lack of a suitable reference, we give  brief proofs of the lemmas although it is well known to experts. 

\begin{lem}\label{invariant}

Let $G$ be an affine algebraic group acting on varieties $X$ and $Y$, and $H\subseteq G$ a finite normal subgroup. Let $\pi:X\to Y$ be a $G$-equivariant morphism. If $H$-action on $Y$ is trivial and $\pi$ is a principal $H$-bundle, we have an  equivalence
\[
\coh_GX\cong \coh_{G/H}Y.
\] 
\begin{proof}
Since $\pi$ is a finite morphism, the direct image  $\pi_*:\Qcoh_GX\to \Qcoh_GY$ preserves coherent sheaves.  We define a functor 
\[
(\pi_*)^H:\coh_GX\to\coh_{G/H}Y
\]
as the composition of $\pi_*:\coh_GX\to \coh_GY$ and $(-)^H:\coh_GY\to \coh_{G/H}Y$. Then  $(\pi_*)^H$ is right adjoint to $\pi_{p}^*:\coh_{G/H}Y\to \coh_GX$ by \cite[Corollary 2.24]{bfk}, where $p:G\to G/H$ is the natural projection. Hence it is enough to show that the adjunction morphisms
$\eta:{\rm id}\to (\pi_*)^H\circ \pi_p^*$ and $\varepsilon:\pi_p^*\circ(\pi_*)^H\to {\rm id}$ are isomorphisms of functors. Consider the following commutative diagram:
\[\xymatrix{
 \coh_GX\ar[rr]^{(\pi_*)^H}\ar[d]_{{\rm Res}^{G}_H}&&\coh_{G/H}Y\ar[rr]^{\pi_p^*}\ar[d]^{{\rm Res}^{G/H}}&&\coh_GX\ar[d]^{{\rm Res}^{G}_H}\\
\coh_HX\ar[rr]^{(\pi_*)^H}&&\coh Y\ar[rr]^{\pi^*}&&\coh_HX
 }\]
 Then $\eta$ and $\varepsilon$ are isomorphisms if and only if so are ${\rm Res}^{G/H}(\eta)$ and ${\rm Res}^G_H(\varepsilon)$. 
Hence the equivalence $\coh_GX\cong \coh_{G/H}Y$ follows from the equivalence $\coh_HX\cong \coh Y$ that is well known to hold.
\end{proof}
\end{lem}

\begin{lem}\label{trivial action}
Let $G$ be an  affine algebraic group acting on $X$ and $H\subseteq G$ an abelian closed normal subgroup of $G$ such that $H$ acts trivially on $X$. We denote by $p:G\to G/H$  the natural projection, and we write  ${H}^{\vee}$ for the set of characters of $H$. 

\begin{itemize}
\item[$(1)$] 
The functor ${\rm id}_p^*:\coh_{G/H}X\to \coh_G X$
is fully faithful, and it induces an equivalence 
\[
{\rm id}_p^*:\coh_{G/H}X\cong (\coh_GX)_{\chi_{{}_0}},
\]
where $\chi_{{}_0}:H\to \bG_m$ is the trivial character.

\item[$(2)$] The functor
\[\bigoplus_{\chi\in{H}^{\vee}}(-)_{\chi}:\coh_G X\simto \bigoplus_{\chi\in {H}^{\vee}}(\coh_GX)_{\chi}\]
is an equivalence.
\item[$(3)$] For characters $\eta:G\to \bG_m$ and $\chi:H\to \G_m$, the tensor product with $\cO(\eta)$ gives an equivalence;
\[
(-)\otimes \cO(\eta):(\coh_GX)_{\chi}\simto(\coh_GX)_{(\eta|_{{}_{H}})\chi}.
\]
\end{itemize}

\begin{proof}
(1) The functor $(-)^H:\coh_G X\to \coh_{G/H}X$ is right adjoint to ${\rm id}_p^*:\coh_{G/H}X\to \coh_G X$ by \cite[Lemma 2.22]{bfk}.  For any $F\in \coh_{G/H}X$, let   $\eta_F:{F}\to\bigl({\rm id}_p^*(F)\bigr)^H$ be the adjunction morphism.  Denote by $\cO_{[X/G]}\in \coh_GX$ and $\cO_{[X/(G/H)]}\in\coh_{G/H}X$ the structure sheaves with natural equivariant  structures induced by group actions. Then, since there is a natural isomorphism $(\cO_{[X/G]})^H\cong \cO_{[X/(G/H)]}$, we have the following isomorphisms
\[
F\cong F\otimes _{\cO_{[X/(G/H)]}}(\cO_{[X/G]})^H\cong ({\rm id}_p^*(F)\otimes_{\cO_{[X/G] }}\cO_{[X/G]})^H\cong \bigl({\rm id}_p^*(F)\bigr)^H,
\]
where the second isomorphism follows from \cite[Lemma 2.23]{bfk}. Since the composition of the above isomorphisms  is equal to the adjunction morphism $\eta_F$, ${\rm id}_p^*$ is fully faithful. The latter claim is obvious by construction.

(2)  It is enough to show that $F\cong \oplus_{\chi\in H^{\vee}}F_{\chi}$ for any $F\in \coh_GX$. By definition we see that $ \oplus_{\chi\in H^{\vee}}F_{\chi}\subseteq F$, and hence it suffices to show that $F/\left(\oplus_{\chi\in H^{\vee}}F_{\chi}\right)\cong 0$.  Since $H$ is abelian and acts trivially on $X$, we have a decomposition $\coh_HX\cong \bigoplus_{\chi\in {H}^{\vee}}(\coh_HX)_{\chi}$. Hence  $\overline{F}:={\rm Res}^G_H(F)\in \coh_HX$ is decomposed into a direct sum $\overline{F}\cong \oplus_{\chi\in {H}^{\vee}} \overline{F}_{\chi}$. 
Since ${\rm Res}^G_H(F_{\chi})= \overline{F}_{\chi}$, we have ${\rm Res}^G_H(F)\cong {\rm Res}^G_H(\oplus_{\chi\in {H}^{\vee}} F_{\chi})$. This implies  ${\rm Res}^G_H(F/\left(\oplus_{\chi\in H^{\vee}}F_{\chi}\right))\cong 0$, since ${\rm Res}^G_H$ is a right exact functor. 
Hence we have $F/\left(\oplus_{\chi\in H^{\vee}}F_{\chi}\right)\cong 0$.

(3) is obvious.
\end{proof}
\end{lem}

\begin{lem}\label{descent}
Let $G$ be an affine algebraic group acting on a variety $X$. For a character  $\chi:G\to \bG_m$  of $G$,  define a $G\times \bG_m$-action on $X\times \bG_m$ by $(g,h_1)\cdot(x,h_2)=(g\cdot x, \chi(g)^{-1}h_1h_2)$ and denote by $\varphi$ the morphism ${\rm id}_G\times \chi:G\to G\times \bG_m$. Then the morphism $e:X\to X\times \bG_m$ defined by $e(x):=(x,1)$ is $\varphi$-equivariant, and the pull-back
\[
e_{\varphi}^*:\coh_{G\times \bG_m}(X\times \bG_m)\simto\coh_GX
\]
is an equivalence.
\begin{proof}
Let $p:X\times \bG_m\to X$ and $\pi:G\times \bG_m\to G$ be the natural projections. Then $p$ is $\pi$-equivariant, and the composition  $e_{\varphi}^*\circ p_{\pi}^*:\coh_GX\to \coh_GX$ is the identity functor. Hence $e_{\varphi}^*$ is essentially surjective. 
Since the set of morphisms in $\coh_{G\times \bG_m}(X\times \bG_m)$ and $\coh_GX$ are the $G$-invariant subspaces of the set of morphisms in  $\coh_{ \bG_m}(X\times \bG_m)$ and $\coh X$ respectively, the fully-faithfulness of $e_{\varphi}^*$ reduces to the case that $G$ is trivial, which follows from \cite[Lemma 1.3]{tho} (see also \cite[Lemma 2.13]{bfk}).
\end{proof}
\end{lem}

\subsection{Graded modules}\label{graded modules}
To fix notation, we give a quick review of the categories of graded modules over graded rings. Let $L$ be a finitely generated abelian group, and let $S=\bigoplus_{l\in L} S_{l}$ be a noetherian commutative ring with  $L$-grading. An element $m\in M$ of $L$-graded $S$-module $M=\bigoplus_{l\in L} M_{l}$  is called  {\it homogeneous} if $m\in M_{l}$ for some $l\in L$. If $m\in M_{l}$, we say that the {\it degree} of $m$ is $l$, and  we write $\deg(m)=l$. An $L$-graded $S$-module $M$ is called {\it finitely generated} if $M$ has finitely many homogeneous  generators.  Denote by $\fmod^LS$  the category of  finitely generated $L$-graded $S$-modules whose morphisms are degree preserving morphisms. For $M\in \fmod^LS$ and an element $l\in L$, we define the {\it $l$-shift} $M(l)$ of $M$ to be the  $L$-graded $S$-module 
\[
M(l):=\bigoplus_{l'\in L}M(l)_{l'}
\]
defined by $M(l)_{l'}:=M_{l+l'}$. The $l$-shift defines the exact autoequivalence 
\[
(-)(l): \fmod^LS\simto \fmod^LS
\]
 of the abelian category $\fmod^LS$.  
 
  Let $L'$ be another finitely generated abelian groups, and $S '$ an $L'$-graded noetherian ring, and 
 suppose that  $\alpha:L\to L'$ is a group homomorphism. A ring homomorphism $\varphi: S\to S'$ is {\it $\alpha$-equivariant} if  $\varphi(S_{l})\subset S'_{\alpha(l)}$ for any $l\in L$. For $M\in \fmod^LS$, we define the $L'$-graded $S'$-module $\varphi_{\alpha}^*M$ to be the $S'$-module $M\otimes_SS'$ with the $L'$-grading structure given by 
 \[
 (M\otimes_SS')_{l'}:=\left\{ \sum m\otimes s' \, \middle| \, m\in M_l, s'\in S'_{l''} \mbox{ for some }  l\in L, l''\in L' \mbox{ with } \alpha(l)+l''=l'\right\}.
 \]
This defines a right exact functor 
 \[
 \varphi_{\alpha}^*:\fmod^{L}S\to \fmod^{L'}S'.
 \]
 Since $\varphi_{\alpha}^*S(l)\cong S'(\alpha(l))$ for any $l\in L$, this functor restricts to the functor 
 $\proj^LS\to \proj^{L'}S'$.
If $S=S'$, $\varphi=\id$ and $\alpha:L\hookto L'$ is an inclusion, then we write simply  
\begin{equation}\label{extension}
M^{L'}:={\id}_{\alpha}^*M,
\end{equation}
and note that we have
\[
(M^{L'})_{l'}=\begin{cases}
M_{l'}&l'\in L\\
0& l'\notin L.
\end{cases}
\]
If $L=L'$, $\alpha=\id$ and $S'$ is a finitely generated $S$-module, we have a forgetful functor 
\[
\varphi_*:\fmod^{L}S'\to \fmod^LS
\]
associated to the ring homomorphism $\varphi$, and we have  $(\varphi_*M)_l:=M_l$ for any $l\in L$.

 \begin{nota}For a finitely generated abelian group $A$, we set  ${\rm G}(A):=\Spec k[A]$, where $k[A]$  is the group ring associated to $A$. Then ${\rm G}(A)$ is a commutative  affine algebraic group. Conversely, for a commutative affine algebraic group  $G$, we denote by  $G^{\vee}:=\Hom(G,\bG_m)$ the character group  of $G$, and then $G^{\vee}$ is a finitely generated abelian group.  It is standard  that we have natural isomorphisms
 ${\rm G}(A)^{\vee}\cong A$ and ${\rm G}(G^{\vee})\cong G.$
 \end{nota}

 The grading structure on $S$ (resp. $S'$) induces the algebraic group action from $\G(L)$ (resp. $\G(L')$) on $\Spec S$ (resp. $\Spec S'$), and  we have a natural exact equivalence 
 \begin{equation}\label{eqcoh  grmod}
 \Gamma:\coh_{\G(L)}\Spec S\simto \fmod^{L}S
 \end{equation}
 of abelian categories (see e.g. \cite[Section 2.1]{bfk2}) given by taking global sections. This equivalence restricts to the equivalence of $\G(L)$-equivariant locally free sheaves and $L$-graded projective modules. If we write \[\widetilde{\varphi}:\Spec S'\to \Spec S\] the  morphism  associated to $\varphi:S\to S'$ and we set $\G(\alpha):\G(L')\to \G(L)$ the morphism of algebraic groups induced by $\alpha:L\to L'$, then the following diagram is commutative:
\[\xymatrix{
 \coh_{\G(L)}\Spec S\ar[rr]^{\widetilde{\varphi}_{\G(\alpha)}^*}\ar[d]_{\Gamma}&&\coh_{\G(L')}\Spec S'\ar[d]^{\Gamma}\\
\fmod^{L}S\ar[rr]^{\varphi_{\alpha}^*}&&\fmod^{L'}S'.
 }\]
Moreover, if $L=L'$, $\alpha=\id$ and $S'$ is a finitely generated $S$-module,  the following diagram commutes:
\[\xymatrix{
 \coh_{\G(L')}\Spec S'\ar[rr]^{\widetilde{\varphi}_*}\ar[d]_{\Gamma}&&\coh_{\G(L)}\Spec S\ar[d]^{\Gamma}\\
\fmod^{L'}S'\ar[rr]^{\varphi_{*}}&&\fmod^{L}S.
 }\]

\subsection{Derived factorization categories}
In this subsection, we provide a brief summary of  the derived factorization categories.

\begin{dfn}
A {\it gauged Landau-Ginzburg  model}, or simply {\it gauged LG model}, is data $(X,\chi,W)^G$ with
$X$ a scheme, $G$ an algebraic group acting on $X$, $\chi:G\rightarrow\mathbb{G}_m$ a character of $G$  and $W:X\rightarrow\mathbb{A}^1$ a $\chi$-semi-invariant regular function, i.e.  $W(g\cdot x)=\chi(g)W(x)$ for any $g\in G$ and any $x\in X$. If $G$ is trivial, we denote the gauged LG model by $(X,W)$, and call it {\it Landau-Ginzburg model} or {\it LG model}. 

\end{dfn}
 
 For a gauged LG model, we consider its factorizations that can be considered as  ``twisted  complexes".

\begin{dfn}
Let $(X,\chi,W)^G$ be a gauged LG model. A {\it quasi-coherent factorization}  of $(X,\chi,W)^G$ is a sequence
$$F=\Bigl(F_1\xrightarrow{\varphi_1^F} F_0\xrightarrow{\varphi_0^F} F_1(\chi)\Bigr),$$
where, for each $i=0,1$, $F_i$ is a $G$-equivariant quasi-coherent sheaf on $X$ and $\varphi_i^F$ is a $G$-invariant homomorphism  such that $\varphi_0^F\circ\varphi_1^F=W\cdot {\rm id}_{F_1}$ and $\varphi_1^F(\chi)\circ\varphi_0^F=W\cdot {\rm id}_{F_0}$. Equivariant quasi-coherent sheaves $F_0$ and $F_1$ in the above sequence are called the {\it components} of  $F$. If the components $F_i$ of  $F$ are coherent  (resp. locally free coherent, injective) sheaves,  then $F$ is called a {\it coherent factorization} (resp. a {\it matrix factorization}, an {\it injective factorization}). We will often call these sequences just {\it factorizations} of $(X,\chi,W)^G$.
\end{dfn}

\begin{dfn}
For a gauged LG model $(X,\chi,W)^G$, we define the abelian category
\[\Qcoh_G(X,\chi,W)\]
whose objects are quasi-coherent factorizations of $(X,\chi,W)^G$, and whose set of morphisms are defined as follows:
For two objects $E,F\in {\rm Qcoh}_G(X,\chi,W)$, we define   ${\rm Hom}(E,F)$ to be the set of pairs $(f_1,f_0)$ of $f_i\in{\rm Hom}_{\Qcoh_GX}(E_i, F_i)$ such that  the following diagram commutes;
\[\xymatrix{
 E_1\ar[rr]^{\varphi_1^E}\ar[d]_{f_1}&&E_0\ar[rr]^{\varphi_0^E}\ar[d]^{f_0}&&E_1(\chi)\ar[d]^{f_1(\chi)}\\
F_1\ar[rr]^{\varphi_1^F}&&F_0\ar[rr]^{\varphi_0^F}&&F_1(\chi).
 }\]
We define  full subcategories \[{\rm MF}_G(X,\chi,W)\subset\coh_G(X,\chi,W)\subset \Qcoh_G(X,\chi,W)\] of  ${\rm Qcoh}_G(X,\chi,W)$ whose objects are matrix factorizations and coherent factorizations respectively. By construction, theses subcategories are exact subcategories.
\end{dfn}

Since factorizations can be considered as `twisted complexes', we can consider the homotopy category  of factorizations.

\begin{dfn}\label{homotopy def}
Two morphisms $f=(f_1,f_0):E\to F$ and $g=(g_1,g_0):E\to F$ in $\Qcoh_G(X,\chi,W)$ are  {\it homotopy equivalent}, denoted by $f\sim g$, if there exist two homomorphisms in $\Qcoh_GX$
$$h_0:E_0\rightarrow F_1\hspace{3mm}{\rm and}\hspace{3mm}h_1:E_1(\chi)\rightarrow F_0$$
such that $f_0-g_0=\varphi_1^Fh_0+h_1\varphi_0^E$ and $f_1(\chi)-g_1(\chi)=\varphi_0^Fh_1+(h_0(\chi))(\varphi_1^E(\chi))$.
\end{dfn}

\vspace{2mm}
The {\it homotopy category of factorizations} of $(X,\chi,W)^G$, denoted by  \[\KQcoh_G(X,\chi,W),\] is defined by ${\rm Obj}(\KQcoh_G(X,\chi,W)):={\rm Obj}(\Qcoh_G(X,\chi,W))$ and the set of morphisms are defined as the set of homotopy equivalence  classes;
$${\rm Hom}_{{\rm KQcoh}(X,\chi,W)}(E,F):={\rm Hom}_{{\rm Qcoh}(X,\chi,W)}(E,F)/\sim.$$
Similarly, we can define the homotopy category $\Kcoh_G(X,\chi,W)$ (resp. ${\rm KMF}(X,\chi,W)$)  of coherent factorizations (resp. matrix factorizations) of $(X,\chi,W)^G$.

\vspace{3mm}
Next we define the totalization of a bounded complex of factorizations, which is an analogy of the total complex of a double complex.

\begin{dfn}
Let $F^{\text{\tiny{\textbullet}}}=(\cdot\cdot\cdot\rightarrow F^i\xrightarrow{\delta^i}F^{i+1}\rightarrow\cdot\cdot\cdot)$ be a bounded complex of $\Qcoh_G(X,\chi,W)$. For $l=0,1$, set
$$T_l:=\bigoplus_{i+j=-l}F^i_{\overline{j}}(\chi^{{\lceil j/2\rceil}}),$$
and define $$t_l:T_l\rightarrow T_{\overline{l+1}}$$
to be the  homomorphism given by
$$t_l|_{F^i_{\overline{j}}(\chi^{{\lceil j/2\rceil}})}:=\delta^i_{\overline{j}}(\chi^{{\lceil j/2\rceil}})+(-1)^i\varphi^{F^i}_{\overline{j}}(\chi^{{\lceil j/2\rceil}}),$$
where $\overline{n}$ is $n$ modulo $2$, and $\lceil m\rceil$ is the minimum integer which is greater than or equal to a real number $m$.
We define the {\it totalization} Tot$(F^{\text{\tiny{\textbullet}}})\in \Qcoh_G(X,\chi,W))$ of $F^{\text{\tiny{\textbullet}}}$ by 
$${\rm Tot}(F^{\text{\tiny{\textbullet}}}):=\Bigl(T_1\xrightarrow{t_1}T_0\xrightarrow{t_0}T_1(\chi)\Bigr).$$
\end{dfn}

\vspace{3mm}In what follows, we will recall that the homotopy category $\KQcoh_G(X,\chi,W)$ has a  structure of  triangulated category, and  $\Kcoh_G(X,\chi,W)$ and $\KMF_G(X,\chi,W)$  are full triangulated subcategories of $\KQcoh_G(X,\chi,W)$.

\vspace*{2mm}
\begin{dfn}
We define an automorphism $T$ on $\KQcoh_G(X,\chi,W)$, which is called the {\it shift functor}, as follows.
For an object $F\in \Kcoh_G(X,\chi,W)$, we define an object $T(F)$ as
$$T(F):=\Bigl(F_0\xrightarrow{-\varphi^F_0}F_1(\chi)\xrightarrow{-\varphi^F_1(\chi)}F_0(\chi)\Bigr),$$
and for a morphism $f=(f_1,f_0)\in {\rm Hom}(E,F)$ we set $T(f):=(f_0,f_1(\chi))\in{\rm Hom}(T(E),T(F))$. For any integer $n\in\mathbb{Z}$, denote by $(-)[n]$ the functor $T^n(-)$.
\end{dfn}

\vspace*{2mm}
\begin{dfn}
Let $f : E\rightarrow F$ be a morphism in ${\rm Qcoh}_G(X,\chi,W)$. We define its {\it mapping cone} Cone$(f)$ to be the totalization of the complex 
\[(\cdot\cdot\cdot\rightarrow0\rightarrow E\xrightarrow{f} F\rightarrow0\rightarrow\cdot\cdot\cdot)\] with $F$ in degree zero. Then the mapping cone $\Cone(f)$ of $f$ is  of the following form

\[
\left(
F_1\oplus E_0
\xrightarrow{\left(\begin{array}{cc}{\varphi}^F_1&f_0\\0&-{\varphi}^E_0\end{array}\right)}
F_0\oplus E_1(\chi)
\xrightarrow{\left(\begin{array}{cc}{\varphi}^F_0&f_1\\0&-{\varphi}^E_1(\chi)\end{array}\right)}
{F}_1(\chi)\oplus E_0(\chi)
\right).
\]

A {\it distinguished triangle} is a sequence in ${\rm KQcoh}_G(X,\chi,W)$ which is isomorphic to a sequence of the form
$$E\xrightarrow{f}F\xrightarrow{i}{\rm Cone}(f)\xrightarrow{p}E[1],$$
where $i$ and $p$ are  natural injection and  projection respectively.
\end{dfn}

The following  is well known.

\begin{prop}\label{homotopy is tri}
The homotopy category $\KQcoh_G(X,\chi,W)$ is a triangulated category with respect to the above shift functor and the above distinguished triangles. The full subcategories $\Kcoh_G(X,\chi,W)$ and $\KMF_G(X,\chi,W)$ are full triangulated subcategories.
\end{prop}

\vspace*{2mm}
Following Positselski (\cite{posi}, \cite{efi-posi}), we define  derived factorization categories (see also \cite{pv,orlov4}).

\begin{dfn}\label{derived fact}
Denote by  $\Acoh_G(X,\chi,W)$ the smallest thick subcategory of $\Kcoh_G(X,\chi,W)$ containing all totalizations of short exact sequences in $\coh_G(X,\chi,W)$.  We define the {\it (absolute) derived factorization category} of $(X,\chi,W)^G$ as the Verdier quotient 
$$\Dcoh_G(X,\chi,W):=\Kcoh_G(X,\chi,W)/\Acoh_G(X,\chi,W).$$
Similarly, we consider the thick subcategory $\AMF_G(X,\chi,W)$ containing all totalizations of short exact sequences in the exact category $\MF_G(X\chi,W)$, and define the {\it (absolute) derived matrix factorization category} by
$$\DMF_G(X,\chi,W):=\KMF_G(X,\chi,W)/\AMF_G(X,\chi,W).$$
We define a larger category $\DcoQcoh_G(X,\chi,W)$ as follows: Denote by $\AcoQcoh_G(X,\chi,W)$ the smallest thick subcategory  of $\KQcoh_G(X,\chi,W)$ that is closed under taking small direct sums and contain all totalizations of short exact sequences in $\Qcoh_G(X,\chi,W)$. Then we define $\DcoQcoh_G(X,\chi,W)$ by
\[
\DcoQcoh_G(X,\chi,W):=\KQcoh_G(X,\chi,W)/\AcoQcoh_G(X,\chi,W).
\]
Factorizations in $\Acoh_G(X,\chi,W)$ or $\AMF_G(X,\chi,W)$ are said to be {\it acyclic}, and factorizations in $\AcoQcoh_G(X,\chi,W)$ are said to be  {\it coacyclic}. Two quasi-coherent factorizations $E$ and $F$ are said to be {\it quasi-isomorphic} if $E$ and $F$ are isomorphic in $\DcoQcoh_G(X,\chi,W)$.
\end{dfn}

\vspace{2mm}

The following  propositions are standard.

\begin{prop}[{\cite[Propostion 2.25.(1)]{H2}}]\label{ff on coh}
The natural functor 
\[
\Dcoh_G(X,\chi,W)\to \DcoQcoh_G(X,\chi,W)
\]
is fully faithful, and the  the thick closure $\overline{\Dcoh_G(X,\chi,W)}$ of the essential image of the functor is  the subcategory of compact objects.
\end{prop}

\begin{prop}[{\cite[Proposition 3.14]{bfk}}]\label{locally free resolution}
If $X$ is a smooth quasi-projective variety, the natural functor 
\[
\DMF_G(X,\chi,W)\to \Dcoh_G(X,\chi,W)
\]
is an equivalence.
\end{prop}

\begin{prop}[cf.\,{\cite[Lemma 2.24]{bdfik}}]\label{affine}
Assume that $X=\Spec R$ is an affine scheme and $G$ is reductive. For $P\in {\rm KMF}_G(X,\chi,W)$ and $A\in\Acoh_G(X,\chi,W)$, we have 
$${\rm Hom}_{\Kcoh_G(X,\chi,W)}(P,A)=0.$$
In particular, the Verdier  localizing functor 
$$\KMF_G(X,\chi,W)\simto \DMF_G(X,\chi,W)$$
is an equivalence.
\begin{proof}

Since $G$ is reductive, the restriction functor $\Res^G:\Kcoh_G(X,\chi,W)\to \Kcoh(X,W)$ is faithful. Hence the problem reduces to the case when $G$ is trivial, and it follows from \cite[Lemma 2.24]{bdfik}. 
\end{proof}
\end{prop}

The following result follows from Lemma \ref{trivial action}.
\begin{lem}\label{trivial action 2}
Notation is same as in Lemma \ref{trivial action}. Let $\overline{\chi}:G/H\to \bG_m$ be a character, and define a character $\chi:G\to \bG_m$ by $\overline{\chi}\circ p$. Let $W:X\to \bA^1$ be a  $\chi$-semi-invariant regular function. Then $W$ is also $\overline{\chi}$-semi-invariant with respect to the induced $G/H$-action on $X$.
For $\eta\in {H}^{\vee}$, we denote by  $\coh_G(X,\chi,W)_{\eta}$  the subcategory of $\coh_G(X,\chi,W)$ consisting of factorizations whose components  lie in $(\coh_GX)_{\eta}$.

\begin{itemize}
\item[$(1)$] 
The functor ${\rm id}_p^*:\coh_{G/H}(X,\overline{\chi},W)\to \coh_G(X,\chi,W)$
is fully faithful, and it induces an equivalence 
\[
{\rm id}_p^*:\coh_{G/H}(X,\overline{\chi},W)\cong \coh_G(X,\chi,W)_{\eta_{{}_0}},
\]
where $\eta_{{}_0}:H\to \bG_m$ is the trivial character.

\item[$(2)$] There is a decomposition of $\coh_G(X,\chi,W)$ into a direct sum 
\[\coh_G(X,\chi,W)\simto \bigoplus_{\eta\in {H}^{\vee}}\coh_G(X,\chi,W)_{\eta}.\]
\item[$(3)$] For   a character $\phi:G\to \bG_m$ of $G$, the tensor product with $\cO(\phi)$ gives an equivalence
\[
(-)\otimes \cO(\phi):\coh_G(X,\chi,W)_{\eta}\simto\coh_G(X,\chi,W)_{(\phi|_{{}_{H}})\eta}.
\]
\end{itemize}
\end{lem}

The following says that the derived factorization categories are  generalizations of derived categories of coherent sheaves.

\begin{cor}\label{derived category} Assume that $\bG_m$ acts trivially on $X$. Let $n>0$ be a positive integer, and denote by $\chi_n:\bG_m\to\bG_m$ the character defined by $\chi_n(a):=a^n$. We write ${\pmb \upmu_n}:=\l\zeta\r\subset \bG_m$ for the subgroup  of $\bG_m$ generated by a primitive $n$-th root of unity $\zeta$. Then we have an orthogonal decomposition; 
\[
\Dcoh_{\bG_m}(X,\chi_{n},0)\cong \bigoplus_{i=1}^n\Db(\coh X)
\]
\begin{proof}
Since the kernel of the surjection $\chi_n:\bG_m\to\bG_m$ is equal to the subgroup ${\pmb \upmu_n}\subset \bG_m$, we have the following short exact sequence:
\[\
1\to {\pmb \upmu_n}\hookto \bG_m\xrightarrow{\chi_n}\bG_m\to 1
\]
Then, by Lemma \ref{trivial action 2}, we obtain a decomposition $\coh_{\bG_m}(X,\chi_n,0)\cong \bigoplus_{i=1}^n\coh_{\bG_m/{\pmb \upmu_n}}(X,\chi_1,0)$, and this  induces the following direct sum decomposition of the derived factorization category
\[
\Dcoh_{\bG_m}(X,\chi_{n},0)\cong \bigoplus_{i=1}^n\Dcoh_{\bG_m/{\pmb \upmu_n}}(X,\chi_1,0).
\]
Hence the result follows from \cite[Proposition 2.14]{H2}.
\end{proof}
\end{cor}

For a noetherian scheme $X$ with an action from an algebraic group $G$, the {\it $G$-equivariant singularity category} $\D_G^{\rm sg}(X)$  is defined by the Verdier quotient
\[
\D_G^{\rm sg}(X):=\Db(\coh_G X)/\Perf_GX.
\]
Similarly, for an abelian group $L$ and a commutative $L$-graded noetherian ring $R$, we define the {\it $L$-graded singularity category} $\D_{\rm sg}^{L}(R)$ by 
\[
\D_{\rm sg}^{L}(R):=\Db(\fmod^LR)/\Kb(\proj^LR).
\]
The natural equivalence in \eqref{eqcoh  grmod} induces  the following  equivalence:
\[
\D_{G}^{\rm sg}(\Spec R)\simto \D_{\rm sg}^{G^{\vee}}(R).
\]
\begin{thm}[{\cite[Theorem 2]{orlov4}}, {\cite[Theorem 3.6]{H2}}]\label{fact sing} Let $(X,\chi,W)^G$ be a gauged LG model. Assume that $X$ is a smooth variety, $G$ is a reductive affine algebraic group and $W$ is flat. Denote by $X_0$ the zero scheme of $W$. Then we have an equivalence
\[
\Dcoh_G(X,\chi,W)\simto \D_G^{\rm sg}(X_0).
\]

\end{thm}

\subsection{Homotopy category of graded matrix factorizations}
We use the same notation  as in Section \ref{graded modules}, and let $f\in S$ be a homogeneous element with $\deg(f)=l\in L$. 
\begin{dfn} 
An {\it $L$-graded factorization} of $f$ is a sequence 
\begin{equation}\label{fact def}
F=\left(F_1\xrightarrow{\varphi_1}F_0\xrightarrow{\varphi_0} F_1(l)\right),
\end{equation}
where each $F_i$ is an $L$-graded finitely generated $S$-module and $\varphi_i$ are degree preserving homomorphisms such that $\varphi_0\circ \varphi_1=f\cdot\id_{F_1}$ and $\varphi_1(l)\circ \varphi_0=f\cdot \id_{F_0}$. Each graded module $F_i$ in \eqref{fact def} is called a {\it component} of $F$,  and we say that an $L$-graded factorization $F$ is an {\it $L$-graded  matrix factorization} if both components are projective $S$-modules. 

For $L$-graded factorizations $E$ and $F$, a {\it morphism} from  $E$ to $F$ is a pair $(\alpha_1,\alpha_0)$ of morphisms $\alpha_i: E_i\to F_i$ in $\fmod^LS$ such that it is compatible with $\varphi^E_i$ and $\varphi_i^F$ ($i=0,1$).
\end{dfn} 

We define the {\it homotopy category of $L$-graded  factorizations} of $f$
\[
\Kmod^{L}(S,f)
\]
whose objects are $L$-graded factorizations of $f$ and the morphisms are homotopy equivalence classes of morphisms, where the homotopy equivalence is defined similarly to Definition \ref{homotopy def}. We define the {\it homotopy category of $L$-graded  matrix factorizations} of $f$  to be the full subcategory
\[
\HMF^L_S(f)\subset \Kmod^{L}(S,f)
\]
of  $L$-graded matrix factorizations of $f$.  Similarly to Proposition \ref{homotopy is tri}, the categories $\Kmod^{L}(S,f)$ and  $\HMF^L_S(f)$ have natural structures of  triangulated categories.   We define the category
\[
\Dmod^{L}(S,f):=\Kmod^{L}(S,f)/\Amod^{L}(S,f),
\]
to be  the Verdier quotient of $\Kmod^{L}(S,f)$ by the subcategory $\Amod^{L}(S,f)\subset \Kmod^{L}(S,f)$, where $\Amod^{L}(S,f)$ is defined similarly to $\Acoh_G(X,\chi,W)$ in Definition \ref{derived fact}. By the identical argument in Propositions \ref{locally free resolution} and \ref{affine}, if $S$ is regular, we have a natural equivalence
\[
\HMF^{L}(f)\simto \Dmod^L(S,f).
\]

The element $f\in S$ induces a $\chi_l$-semi-invariant regular function $f:\Spec S\to \bA^1$, where $\chi_l\in \G(L)^{\vee}$ the character of $\G(L)$ corresponding to the isomorphism $\G(L)^{\vee}\cong L$. The equivalence \eqref{eqcoh grmod} implies the following:
 
 \begin{prop}\label{coh mod}
 We have  natural equivalences
 \begin{eqnarray*}
\Dcoh_{\G(L)}(\Spec S,\chi_l,f)&\simto&\Dmod^{L}(S,f)  \\
 \KMF_{\G(L)}(\Spec S,\chi_l,f)&\simto& \HMF^{L}_{S}(f)
 \end{eqnarray*}
of triangulated categories, and the following diagram of natural equivalences commutes
\[
 \begin{tikzcd}
 \KMF_{\G(L)}(\Spec S,\chi_l,f)\arrow[r]\arrow[d]&\HMF^{L}_{S}(f)\arrow[d]\\
\Dcoh_{\G(L)}(\Spec S,\chi_l,f)\arrow[r]&\Dmod^{L}(S,f)  .
 \end{tikzcd}
 \]
 \end{prop}
 
 If $S$ is regular, we write
 $
 \underline{\CM}^{L}(S/f)
 $
 for the stable category of the Frobenius category  of maximal Cohen-Macaulay $L$-graded modules over the hypersurface $S/f$.
 The following is well known.

 \begin{thm}[\cite{buc,eis,orlov3}]\label{mf=cm}
Let $L$ be an abelian group, and  suppose that a polynomial ring $S^n:=k[x_1,\hdots,x_n]$ has an $L$-grading. Let $f\in S^n$ be a non-zero homogeneous element.  Then there are  exact equivalences
 \[
 \HMF^{L}_{S^n}(f)\simto \underline{\CM}^{L}(S^n/f)\simto\D_{\rm sg}^{L}(S^n/f)
 \]
 \end{thm} 
 
 By Auslander--Reiten duality (see \cite[Theorem 3.8]{kst2}) or more general result by Favero--Kelly \cite[Theorem 1.2]{fk}, we have the following.
 
 \begin{thm}
 Set $\vec{\bf x}^n:=\sum_{i=1}^n\vec{x}_i\in L_f$. Then the autoequivalence 
 \[
 S:=(-)(-\vec{\bf x}^n)[n]:\HMF^{L_f}_{S^n}(f)\simto \HMF^{L_f}_{S^n}(f)
 \]
 is the Serre functor on  $\HMF^{L_f}_{S^n}(f)$. 
  \end{thm}

\subsection{Functors of derived factorization categories}
 In this subsection, we recall fundamental functors between factorization categories. Let $X$ and $Y$ be noetherian schemes, and  $G$ an  algebraic group acting on $X$ and $Y$. Let $\chi:G\to \bG_m$ be a character of $G$.
 
 \subsubsection{Direct image and  inverse image}
 
 Let $f:X\to Y$ be a $G$-equivariant morphism. If $W:Y\to \bA^1$ is a $\chi$-semi-invariant regular function, there are gauged LG models $(X,\chi,f^*W)$ and $(Y,\chi,W)$. Then the direct image $f_*:\Qcoh_GX\to \Qcoh_GY$ and the inverse image $f^*:\Qcoh_GY\to \Qcoh_GX$ of $G$-equivariant sheaves naturally define the direct image and the inverse image between factorizations;
 \[
 f_*:\Qcoh_G(X,\chi,f^*W)\to \Qcoh_G(Y,\chi,W)
 \]
\[
f^*:\Qcoh_G(Y,\chi,W)\to \Qcoh_G(X,\chi,f^*W).
\]
The inverse image $f^*$ preserves coherent factorizations, and, if $f$ is proper, the direct image $f_*$ also preserves coherent factorizations. Since these functors preserves homotopy equivalences of morphisms, we have the induced functors of the homotopy categories of factorizations;
  \[
 f_*:\KQcoh_G(X,\chi,f^*W)\to \KQcoh_G(Y,\chi,W)
 \]
\[
f^*:\KQcoh_G(Y,\chi,W)\to \KQcoh_G(X,\chi,f^*W).
\]
Furthermore, if $f$ is affine,  $f_*$ maps $\AcoQcoh_G(X,\chi,f^*W)$ to $\AcoQcoh_G(Y,\chi,W)$, and so $f_*$ induces a functor
\[
f_*:\DcoQcoh_G(Y,\chi,W)\to \DcoQcoh_G(X,\chi,f^*W).
\]
For general $f$, by \cite[Corollary 2.25]{bdfik} we can  define  the right derived functor of $f_*$;
 \[
{\bf R}f_*:\DcoQcoh_G(Y,\chi,W)\to \DcoQcoh_G(X,\chi,f^*W).
\]
If $f$ is proper, the functor ${\bf R}f_*$ preserves coherent factorizations.

 On the other hand, if $f$ is flat, $f^*$ preserves coacyclic factorizations, and hence $f^*$ naturally induces the functor
\[
f^*:\DcoQcoh_G(Y,\chi,W)\to \DcoQcoh_G(X,\chi,f^*W).
\]
If $Y$ is a smooth variety,  by \cite[Proposition 3.14]{bfk} we can define the  left derived functor 
\[
{\bf L}f^*:\DcoQcoh_G(Y,\chi,W)\to \DcoQcoh_G(X,\chi,f^*W)
\]
of $f^*$ without the assumption of flatness of $f$, and  this functor ${\bf L}f^*$ preserves coherent factorizations.
See e.g. \cite{ls,bfk,H1} for more details.

\subsubsection{Tensor product}

Let $W_1:X\to \bA^1$ and $W_2:X\to \bA^1$ be  $\chi$-semi-invariant functions. We define the tensor products 
\[
(-)\otimes (-):\coh_G(X,\chi,W_1)\times \coh_G(X,\chi,W_2)\to \coh_G(X,\chi,W_1+W_2)
\]
of factorizations by
$$E\otimes F:=\Bigl(\bigoplus_{i=0,1}(F_{i}\otimes E_{\overline{i+1}})\xrightarrow{\varphi_1^{E\otimes F}}\bigoplus_{i=0,1}(F_{i}\otimes E_{i})(\chi^i)\xrightarrow{\varphi_0^{E\otimes F}}      \bigoplus_{i=0,1}(F_{i}\otimes E_{\overline{i+1}})(\chi)\Bigr),$$
where $\overline{n}$ is $n$ modulo 2,  
$${\varphi_1^{E\otimes F}}=\begin{pmatrix}
\varphi^E_1\otimes1&1\otimes\varphi^F_1\\ -1\otimes\varphi^F_0&\varphi^E_0\otimes1
\end{pmatrix}
$$
and 
$$\varphi_0^{E\otimes F}=\begin{pmatrix}
\varphi^E_0\otimes1 & -(1\otimes\varphi^F_1)(\chi)\\
1\otimes\varphi^F_0 & (\varphi^E_1\otimes1)(\chi)
\end{pmatrix}.$$
Since this bi-functor preserves the homotopy equivalences of morphisms, we have the induced functor 
\[
(-)\otimes (-):\Kcoh_G(X,\chi,W_1)\times \Kcoh_G(X,\chi,W_2)\to \Kcoh_G(X,\chi,W_1+W_2).
\]
If $F\in \MF_G(X,\chi,W_1)$,  the functor $F\otimes (-):\Kcoh_G(X,\chi,W_2)\to \Kcoh_G(X,\chi,W_1+W_2)$ preserves acyclic factorizations, and thus it induces the functor
\[
(-)\otimes (-):\DMF_G(X,\chi,W_1)\times \Dcoh_G(X,\chi,W_2)\to \Dcoh_G(X,\chi,W_1+W_2).
\]
If $X$ is a smooth variety, by Proposition \ref{locally free resolution} we can define the {\it left derived functor}
\[
(-)\otimes^{\bf L} (-):\Dcoh_G(X,\chi,W_1)\times \Dcoh_G(X,\chi,W_2)\to \Dcoh_G(X,\chi,W_1+W_2)
\]
of  tensor products.

\subsection{Supports of factorizations}

We recall the supports of factorizations, and we provide a lemma that we need for a slight generalization of global Kn\"orrer periodicity in the next subsection. 
In this subsection, we only consider the case that $G$ is trivial. Let $X$ be a smooth variety, and $W:X\to \bA^1$  a  regular function on $X$. 

For a point $p\in X$ in the scheme $X$, we set $X_p:=\Spec\cO_{X,p}$. Since the functor $(-)_p:\coh X\to \coh X_p$ defined by taking the stalk at $p$ is  exact,  it induces a functor between derived factorization categories;
\[
(-)_p:\Dcoh(X,W)\to \Dcoh(X_p,W_p),
\]
where $W_p:X_p\to \bA^1$ is  the stalk of $W$ at $p$. Obviously this functor preserves matrix factorizations;
\[
(-)_p:\DMF(X,W)\to \DMF(X_p,W_p).
\]
\begin{dfn}
For $F\in \Dcoh(X,W)$, we define the subset $\Supp(F)\subset X$ by
\[
\Supp(F):=\{\,p\in X\mid F_p\neq 0 \mbox{ in } \Dcoh(X_p,W_p)\}.
\]
\end{dfn}

By Proposition \ref{locally free resolution}, we have an equivalence $\Phi:\DMF(X,W)\simto\Dcoh(X,W)$, and $\Phi$ commutes with the functor $(-)_p$. Hence we have 
\[
\Supp(F)=\Supp(\Phi^{-1}(F)),
\]
and thus $\Supp(F)$ is a closed subset of $X$ by \cite[Proposition 2.20.(2)]{H3}. We say that a thick subcategory $\cT\subset \Dcoh(X,W)$ is {\it closed under tensor action from $\DMF(X,0)$} if for any  $E\in \DMF(X,0)$ and any $F\in \cT$ we have $E\otimes F\in \cT$.  By \cite[Lemma 2.25]{H3}, for any subset $S\subseteq X$ of $X$, the full subcategory 
\[
\{F\in \Dcoh(X,W)\mid \Supp(F)\subseteq S\}
\]
is a triangulated subcategory that is closed under direct summands and tensor action from $\Dcoh(X,0)$. 
The following lemma will be necessary in the proof of Theorem \ref{global knorrer}.(3). 

\begin{lem}[{\cite[Theorem 1.1, Proposition 5.3.(2)]{H3}}]\label{tensor result}
Let $F\in \Dcoh(X,W)$. Assume that $W$ is a non-zero-divisor, and denote by $X_0\subset X$ the zero scheme of $W$. 

\begin{itemize}
\item[$(1)$] The support $\Supp(F)$ of $F$ is contained in the singular locus $\Sing(X_0)$ of $X_0$.

\item[$(2)$] There is an object $E\in  \Dcoh(X,W)$ such that $\Supp(E)=\Sing(X_0)$.

\item[$(3)$] If $\cT\subseteq \Dcoh(X,W)$ is a thick subcategory  closed under tensor action from $\DMF(X,0)$, then  there exists a unique specialization-closed subset $Z\subseteq X$ such that \[\cT=\{F\in \Dcoh(X,W)\mid \Supp(F)\subseteq Z\}.\] 
In particular, if $\cT$ contains an object $F$ with $\Supp(F)=\Sing(X_0)$, then $\cT=\Dcoh(X,W)$.
\end{itemize}
\end{lem}

\subsection{Kn\"orrer periodicity}
In this subsection, we prove a slight generalization of the Kn\"orrer periodicity by \cite{H2,isik,orlov1,shipman}, that  will be necessary in the proof of Theorem \ref{main result}. 

 Let $X$ be a smooth quasi-projective variety over $k$, and let $G$ be a reductive affine algebraic group acting on $X$. Let $\mathcal{E}$ be a $G$-equivariant locally free sheaf on $X$ of finite rank, and choose a $G$-invariant regular section $s\in\Gamma(X,\mathcal{E}^{\vee})^G$ of the dual $\cE^{\vee}:=\cH om_X(\cE,\cO_X)$ of $\cE$. Denote by $Z_s\subset X$ the zero scheme of $s$. Let $\chi:G\rightarrow \mathbb{G}_m$ be a character of $G$.
 Then $\mathcal{E}(\chi)$ induces a vector bundle  $\rm{V}(\mathcal{E}(\chi)):=\underline{\rm Spec}(Sym(\mathcal{E}(\chi)^{\vee}))$ over $X$ with a $G$-action induced by the equivariant structure of $\mathcal{E}(\chi)$. Let $q:{\rm V}(\mathcal{E}(\chi))\rightarrow X$ and $p:{\rm V}(\mathcal{E}(\chi))|_Z\rightarrow Z_s$ be natural projections, and let $i:{\rm V}(\mathcal{E}(\chi))|_{Z_s}\hookto{\rm V}(\mathcal{E}(\chi))$ and $j:Z_s\hookto X$ be the closed immersions.
Now we have the following cartesian square:
\[
 \begin{tikzcd}
 {\rm V}(\mathcal{E}(\chi))|_{Z_s}\arrow[r,hookrightarrow, near start, pos=0.4,"i"]\arrow[d,"p"]&\rm{V}(\mathcal{E}(\chi))\arrow[d,"q"]\\
Z_s\arrow[r,hookrightarrow,"j"]&X.
 \end{tikzcd}
 \]
  The regular section $s$ induces a $\chi$-semi-invariant regular function $Q_s:{\rm V}(\mathcal{E}(\chi))\rightarrow \mathbb{A}^1$. Let $W:X\rightarrow \mathbb{A}^1$ be a $\chi$-semi invariant regular function. 
    
  \begin{thm}[cf. \cite{H2,isik,orlov1,shipman}]\label{global knorrer}
 Assume that one of the following conditions holds:
 \begin{itemize}
 \item[$(1)$] The restricted function $W|_{Z_s}:Z_s\rightarrow \mathbb{A}^1$ is flat.
 
 \item[$(2)$] There is a reductive algebraic group $H$ such that $G=H\times\bG_m$ and $1\times\bG_m\subset G$ acts trivially on $X$. Moreover, $W=0$ and $\chi:H\times\bG_m\to \bG_m$ is the projection.
 
 \item[$(3)$]  $W|_{Z_s}=0$, $Z_s$ is smooth, and $\Dcoh_G(Z_s,\chi,W|_{Z_s})$ is idempotent complete.
 \end{itemize}
 Then we have an equivalence
 \[i_*p^*:{\rm Dcoh}_G(Z_s,\chi,W|_{Z_s})\xrightarrow{\sim}{\rm Dcoh}_G({\rm V}(\mathcal{E}(\chi)),\chi, q^*W+Q_s).\]

 \begin{proof}
 (1) This is \cite[Theorem 1.2]{H2}. 
 
 (2) This is \cite{isik,shipman} for  trivial $H$, and the result  for general $H$ is  \cite[Proposition 4.8]{H2}.
 
 (3)  To ease notation we denote $Z:=Z_s$, ${\rm V}:={\rm V}(\mathcal{E}(\chi))$ and $Q:=Q_s:{\rm V}\to\bA^1$. We prove the result by the following three steps.
 
  \medskip
\noindent
 {\bf Step 1:} In the first step, we consider the case when $G=0$, and we prove the functor 
  \[
 i_*p^*:{\rm D^{co}Qcoh}(Z,W|_{Z})\to{\rm D^{co}Qcoh}({\rm V}, q^*W+Q)
 \]
 between larger categories is an equivalence.
By \cite[Lemma 4.5]{H2}, this  functor is fully faithful, and has a right adjoint \[p_*i^!:{\rm D^{co}Qcoh}({\rm V}, q^*W+Q)\to{\rm D^{co}Qcoh}(Z,W|_Z).\] by \cite[Theorem 3.8]{efi-posi}. Hence it suffices to show  that $p_*i^!$ is also fully faithful.
 For this, we show that, for any $F\in {\rm D^{co}Qcoh}({\rm V}, q^*W+Q_s)$, the adjunction morphism 
 \begin{equation}
 \sigma_F:i_*p^*p_*i^!(F)\to F\nonumber
 \end{equation}
  is an isomorphism, or equivalently, the cone $C(\sigma_F):=\Cone(\sigma_F)$  is the zero object. 
 
 By  \cite[Proposition 1.10]{efi-posi}, the set of objects in  the category $\Dcoh({\rm V}, q^*W+Q)$  is a set of compact generators in ${\rm D^{co}Qcoh}({\rm V}, q^*W+Q)$. Therefore, by \cite[Lemma 2.1.1]{ss},  the smallest triangulated subcategory of ${\rm D^{co}Qcoh}({\rm V}, q^*W+Q)$  that contains $\Dcoh({\rm V}, q^*W+Q)$ and  is closed under (infinite) coproducts is equal to the whole category ${\rm D^{co}Qcoh}({\rm V}, q^*W+Q)$. Moreover, since $i_*p^*$ and $i^!$ admit right adjoint functors, these functors commute with  coproducts, and so does the direct image $p_*$  by \cite[Lemma 1.4]{nee2}. Thus the composition $i_*p^*\circ p_*i^!$ also commutes with  coproducts. Therefore, to show that the adjunction $\sigma_F$ is an isomorphism, we may assume that $F$ lies in $\Dcoh({\rm V}, q^*W+Q)$. Denote by $\cT:=\overline{{\sf Im}(i_*p^*)}$ the thick closure of the triangulated subcategory ${\sf Im}(i_*p^*)\subset \Dcoh({\rm V}, q^*W+Q)$. We claim that $\sigma_F$ is an isomorphism if  $F\in \cT$. Indeed, if $F\in {\sf Im}(i_*p^*)$,  $\sigma_F$ is an isomorphism since $(p_*i^!)\circ(i_*p^*)\cong {\rm id}$. For arbitrary $F\in \cT$, there is an object $E\in \cT$ such that $F\oplus E\in {\sf Im}(i_*p^*)$. Then since $C(\sigma_{F})\oplus C(\sigma_{E})\cong C(\sigma_{(F\oplus E)})\cong 0$, we have $C(\sigma_{F})\cong 0$. Thus, it is enough to show the following claim:
 
\medskip
\noindent
 {\it Claim:} We have $\cT=\Dcoh({\rm V}, q^*W+Q)$.
 
 \medskip

First, we claim that $\cT$ is closed under tensor action from $\DMF({\rm V},0)$. Since $X$ is quasi-projective, $X$ has an ample line bundle $\cL_X$. If we denote $\cL:=q^*\cL_X$, $\cL$ is an ample line bundle on ${\rm V}$ since $q$ is quasi-affine. Then the category $\DMF({\rm V},0)$ is generated by objects of the form $\widetilde{\cL^n}:=(0\to \cL^n\to 0)$ for $n\in \bZ$. Hence it suffices to show that $\cT$ is closed under tensor product with  $\widetilde{\cL^n}$. For  any $E\in \cT$, there is an object $E'\in \cT$ and $F\in \Dcoh(Z,W|_{Z})$ such that $i_*p^*(F)\cong E\oplus E'$. Then we have $(E\otimes \widetilde{\cL^n})\oplus(E'\otimes \widetilde{\cL^n})\cong i_*p^*(F)\otimes\widetilde{\cL^n}\cong i_*(p^*(F)\otimes i^*\widetilde{\cL^n})\cong i_*p^*\bigl(F\otimes (0\to j^*\cL^n_X\to 0)\bigr)\in {\sf Im}(i_*p^*)$. Hence $E\otimes \widetilde{\cL^n}\in \cT$.  
 
 Since $\cT$ is a thick subcategory that is closed under tensor action from $\DMF({\rm V},0)$, by Lemma \ref{tensor result}.(3), there is a specialization-closed subset $Y\subset \Sing (V_0)$ of the singular locus of the zero scheme $V_0\subset {\rm V}$ of $q^*W+Q:{\rm V}\to \bA^1$ such that 
 \[
 \cT=\{E\in \DMF({\rm V}, q^*W+Q)\mid \Supp(E)\subset Y\}.
 \]   
In order to prove the claim,   we only need   to show that $Y=\Sing(V_0)$ by Lemma \ref{tensor result}. Since $Z$ is smooth, $\Sing(V_0)=Z$, where $Z$ is considered as a closed subscheme of ${\rm V}$ via the zero section $Z\hookto {\rm V}|_Z$. For a point $x\in Z$, we denote by the same notation $x$  the points in ${\rm V}|_Z$ and ${\rm V}$ that correspond to $x$ via the zero section. Then we can show that the localized functor $(i_x)_*\circ (p_x)^*:\Dcoh(Z_x,0)\to \Dcoh({\rm V}_x, (q^*W+Q)_x)$ is also fully faithful by using \cite[Lemma 4.4]{H2} and an  argument that is identical to  \cite[Lemma 4.5]{H2}, where $p_x:\bigr({\rm V}|_Z\bigl)_x\to Z_x$ and $i_x:\bigr({\rm V}|_Z\bigl)_x\to{\rm V}_x$ are the localized morphisms. 
 Consider a matrix factorization $\widetilde{\cO_Z}$ defined by   $\widetilde{\cO_Z}:=(0\to \cO_Z\to0)\in \Dcoh(Z,0)$. Then it is easy to see that   $\Supp(\widetilde{\cO_Z})=Z$ by \cite[Proposition 2.30]{ls}.  Since we have $i_*p^*(\widetilde{\cO_Z})_x\cong (i_x)_*\circ (p_x)^*((\widetilde{\cO_Z})_x)$ by flat base changes and the localized functor $(i_x)_*\circ (p_x)^*$ is fully faithful, we have $\Supp(i_*p^*(\widetilde{\cO_Z}))=\Supp(\widetilde{\cO_Z})=Z$. Hence  $Y=Z$, since $i_*p^*(\widetilde{\cO_Z})\in \cT$.

  \medskip
\noindent
 {\bf Step 2:}  In the second step, we show that the functor 
 \[
 i_*p^*:{\rm D^{co}Qcoh}_G(Z,\chi,W|_Z)\to{\rm D^{co}Qcoh}_G({\rm V},\chi, q^*W+Q)
 \]
 is an equivalence. The restriction functors $\Res^G$ between equivariant quasi-coherent sheaves induces the following functors 
 \[\Res^G:{\rm D^{co}Qcoh}_G(Z,\chi,W|_Z)\to {\rm D^{co}Qcoh}(Z,W|_Z)\] 
 \[\Res^G:{\rm D^{co}Qcoh}_G({\rm V},\chi, q^*W+Q)\to{\rm D^{co}Qcoh}({\rm V}, q^*W+Q),\] 
 and these functors have right adjoint functors, denoted by $\Ind^G$ (see  \cite[Definition 2.14]{bfk} for the definition of $\Ind^G$). Then by Example \ref{adj} we have the induced comonads, denoted by $\bT$, on  ${\rm D^{co}Qcoh}(Z,W|_Z)$ and ${\rm D^{co}Qcoh}({\rm V}, q^*W+Q)$. Since $G$ is reductive, by \cite[Lemma 4.56]{H1}, the comparison functors
 \[
\Gamma: {\rm D^{co}Qcoh}_G(Z,\chi,W|_Z)\to {\rm D^{co}Qcoh}(Z,W|_Z)_{\bT}
 \]
 \[
\Gamma: {\rm D^{co}Qcoh}_G({\rm V},\chi, q^*W+Q)\to {\rm D^{co}Qcoh}({\rm V}, q^*W+Q)_{\bT}
 \]
 are equivalences (note that in our setting, by \cite[Remark 4.4]{H1}, ${\rm DQcoh}_G(-)$ in \cite{H1} is equivalent to ${\rm D^{co}Qcoh}_G(-)$ since ${\rm V}$ and $Z$ are smooth). By  \cite[Lemma 2.11]{H1} we see that the functor $i_*p^*:{\rm D^{co}Qcoh}(Z,W|_Z)\to{\rm D^{co}Qcoh}({\rm V}, q^*W+Q)$ is a linearizable functor, and we have an induced functor 
 \[
 (i_*p^*)_{\bT}:{\rm D^{co}Qcoh}(Z,W|_Z)_{\bT}\to{\rm D^{co}Qcoh}({\rm V}, q^*W+Q)_{\bT}
 \]
 that commutes with $i_*p^*:{\rm D^{co}Qcoh}_G(Z,\chi,W|_Z)\to{\rm D^{co}Qcoh}_G({\rm V},\chi, q^*W+Q)$ and the comparison functors.  Proposition \ref{main prop 2} implies that the induced functor $(i_*p^*)_{\bT}$ is an equivalence by Step 1, and hence 
 the functor $i_*p^*:{\rm D^{co}Qcoh}_G(Z,\chi,W|_Z)\to{\rm D^{co}Qcoh}_G({\rm V},\chi, q^*W+Q)$ is also an equivalence.

   \medskip
\noindent
 {\bf Step 3:} In the final step, we finish the proof. By Step 2, we have already shown that the functor $i_*p^*:{\rm Dcoh}_G(Z,\chi,W|_Z)\to{\rm Dcoh}_G({\rm V},\chi, q^*W+Q)$ is fully faithful. By Step 2 and  Proposition \ref{ff on coh}, the equivalence $ i_*p^*:{\rm D^{co}Qcoh}_G(Z,\chi,W|_Z)\to{\rm D^{co}Qcoh}_G({\rm V},\chi, q^*W+Q)$ gives the  equivalence
 \[
 i_*p^*:\overline{{\rm Dcoh}_G(Z,\chi,W|_Z)}\xrightarrow{\sim}\overline{{\rm Dcoh}_G({\rm V},\chi, q^*W+Q)}
 \]
of compact objects, where $\overline{(-)}$ denotes the thick closure of $(-)$. But since ${\rm Dcoh}_G(Z,\chi,W|_Z)$ is idempotent complete, we have $\overline{{\rm Dcoh}_G(Z,\chi,W|_Z)}={\rm Dcoh}_G(Z,\chi,W|_Z)$. Hence the functor $ i_*p^*:{\rm Dcoh}_G(Z,\chi,W|_Z)\to{\rm Dcoh}_G({\rm V},\chi, q^*W+Q)$ is essentially surjective. This finishes the proof.
 \end{proof}
 \end{thm}

\subsection{DG-enhancements and  derived Morita theory}

We recall the definition of the Thom--Sebastiani sum of gauged LG models:
\begin{dfn}\label{def of TS sum} The {\it Thom--Sebastiani sum} of  two gauged LG models $(X_1,\chi_1,W_1)^{G_1}$ and $(X_2,\chi_2,W_2)^{G_2}$  is defined to be the gauged LG model \[(X_1\times X_2,\chi_1\times_{\bG_m}\chi_2,W_1\boxplus W_2)^{G_1\times_{\bG_m}G_2 },\]
where $G_1\times_{\bG_m}G_2:=\{(g_1,g_2)\in G_1\times G_2\mid \chi_1(g_1)=\chi_2(g_2)\}$,  $(\chi_1\times_{\bG_m}\chi_2)(g_1,g_2):=\chi_1(g_1)=\chi_2(g_2)$, and $W_1\boxplus W_2:=p_1^*W_1+p_2^*W_2$ (here  $p_i$ are natural projections). We also call the potential $W_1\boxplus W_2$ the {\it Thom--Sebastiani sum} of  $W_1$ and $W_2$.
\end{dfn}

We recall the derived Morita theory for factorizations by \cite{bfk} which is a general version of  \cite{dyc,pv2}. For this, we recall dg-enhancements of derived factorization categories.

\begin{dfn}
Let $(X,\chi,W)^G$ be a gauged LG model. We define the dg-category
\[{\sf Qcoh}_G(X,\chi,W)\]
to be the category of quasi-coherent factorizations of $(X,\chi,W)^G$ whose Hom-complexes are defined  as follows: For any $E,F\in{\sf Qcoh}_G(X,\chi,W)$, we 
 define the complex ${\rm Hom}(E,F)^{\text{\tiny{\textbullet}}}$ of morphisms from $E$ to $F$ as the following  graded vector space
$${\rm Hom}(E,F)^{\text{\tiny{\textbullet}}}:=\bigoplus_{n\in\mathbb{Z}}{\rm Hom}(E,F)^n$$
with a differential $d^i:{\rm Hom}(E,F)^i\rightarrow{\rm Hom}(E,F)^{i+1}$ given by
$$d^i(f):=\varphi^F\circ f-(-1)^i f\circ\varphi^E,$$
where 
$${\rm Hom}(E,F)^{2m}:={\rm Hom}(E_1,F_1(\chi^m))\oplus{\rm Hom}(E_0,F_0(\chi^m))$$
$${\rm Hom}(E,F)^{2m+1}:={\rm Hom}(E_1,F_0(\chi^m))\oplus{\rm Hom}(E_0,F_1(\chi^{m+1})).$$
\end{dfn}

\vspace{3mm}
We denote by 
${\sf Inj}_G(X,\chi,W)$ (resp. 
${\sf MF}_G(X,\chi,W)$)
 the  full dg-subcategory of  ${\sf Qcoh}_G(X,\chi,W)$ consisting of injective factorizations (resp. matrix factorizations). 
We also consider the  full dg-subcategory ${\sf inj}_G(X,\chi,W)\subset{\sf Inj}_G(X,\chi,W)$ consisting of injective factorizations that are quasi-isomorphic to coherent factorizations. Note that  the homotopy category $[{\sf MF}_G(X,\chi,W)]$ of the dg-category ${\sf MF}_G(X,\chi,W)$ is equal to  $\KMF_G(X,\chi,W)$.
Hence, if $X$ is an affine scheme and $G$ is reductive,  then we have the following equivalence  by Lemma \ref{affine};
\[
\bigl[{\sf MF}_G(X,\chi,W)\bigr]\simto\DMF_G(X,\chi,W).
\] 
The following standard result provides  dg-enhancements of  general gauged LG models.
 
 \begin{prop}[{\cite[Lemma 2.12]{H2}}, {\cite[Corollary 2.25]{bdfik}}]\label{dg enh}
 Let $(X,\chi,W)^G$ be a gauged LG model, and assume that $X$ is noetherian. Then  the natural functors 
   \[\bigl[{\sf Inj}_G(X,\chi,W)\bigr]\to \DcoQcoh_G(X,\chi,W)\]
 \[
 \bigl[{\sf inj}_G(X,\chi,W)\bigr]\to \Dcoh_G(X,\chi,W)
 \]
 are equivalences.
 
 \end{prop}

In what follows, we recall the derived Morita theory for derived factorization categories.  We  freely use the terminology and notation from \cite{toen}. 
Recall that the natural inclusion functor $Ho(dg\mathchar`-cat^{tr})\to Ho(dg\mathchar`-cat)$ has a left adjoint functor $\widehat{(-)}_{pe}:Ho(dg\mathchar`-cat)\to Ho(dg\mathchar`-cat^{tr})$ that sends a dg-category $\cT$ to its triangulated hull ${\widehat{\cT}}_{pe}$. For two objects $\cT,\cS\in Ho(dg\mbox{-}cat)$, the {\it Morita product} $\cT\circledast\cS\in Ho(dg\mathchar`-cat^{tr})$ of $\cT$ and $\cS$ is defined by 
 \[
 \cT\circledast\cS:=\widehat{(\cT\otimes_k\cS)}_{pe}.
 \]
 The following is a special case of the derived Morita theory in \cite{bfk}.

 \begin{thm}[{\cite[Corollary 5.18]{bfk}}]\label{morita general}
 Let $X_1$ and $X_2$ be smooth varieties, and $G_1$ and $G_2$ affine algebraic groups. For each $i=1,2$, assume that $G_i$ acts on $X_i$, and  let $W_i:X_i\to \bA^1$ be a $\varphi_i$-semi-invariant regular functions, where $\varphi_i:G_i\to \bG_m$ is a character. Assume the following conditions:
  \begin{itemize}
  \item[$(i)$] The singular locus $\Sing(Z_{W_1}\times Z_{W_2})$ of the product of the zero schemes $Z_{W_i}$  of $W_i$ is contained in $\Sing(Z_{W_1\boxplus W_2})$.
  
   \item[$(ii)$] The  category $\Dcoh_{G_1\times_{\bG_m} G_2}(X_1\times X_2,\varphi_1\times_{\bG_m}\varphi_2, W_1\boxplus W_2)$ is idempotent complete. 
   \end{itemize}
 Then we have the following quasi-equivalence
 \[
 {\sf inj}_{G_1}(X_1,\varphi_1,W_1)\circledast {\sf inj}_{G_2}(X_2,\varphi_2,W_2)\simto  {\sf inj}_{G_1\times_{\bG_m} G_2}(X_1\times X_2,\varphi_1\times_{\bG_m}\varphi_2,W_1\boxplus W_2).
 \]
Here the affine algebraic group $G_1\times _{\bG_m}G_2$ and its character $\varphi_1\times_{\bG_m}\varphi_2$ are defined as in Definition \ref{def of TS sum}.
 \end{thm}

\section{Semi-orthogonal decompositions from sums of potentials}

In this section, we first prove semi-orthogonal decompositions of derived factorization categories of sums of potentials. The key ingredients of this semi-orthogonal decompositions are Theorem \ref{global knorrer} and semi-orthogonal decompositions arising from variations of GIT quotients \cite{segal, h-l, vgit}.
We  then discuss its application to invertible polynomials of chain type.   

\subsection{Semi-orthogonal decompositions of the sums of potentials: General version}\label{general sod section}

Denote by $\zeta\in k$ a primitive $N$-th root of unity, and define a cyclic group ${\pmb \upmu_N}=\l\zeta\r$ to be the subgroup of $\bG_m$ generated by $\zeta$. 
 Let   $G$ be a reductive affine algebraic group together with an injection $\iota:{\pmb \upmu_N}\hookrightarrow G$  of algebraic groups such that its image $\iota({\pmb \upmu_N})\subset G$ is a normal subgroup. Denote by $\pi: G\to G/{\pmb \upmu_N}$ the natural projection.
Let $X$ be a smooth quasi-projective variety with an action from $G/{\pmb \upmu_N}$ such that $X$ has a $G/{\pmb \upmu_N}$-invariant affine open covering. Let $\psi:G\to \bG_m$ be a character of $G$ such that the composition $\psi\circ\iota:{\pmb \upmu_N}\to \bG_m$ is the natural inclusion.  Since $\psi^N\circ\iota:{\pmb \upmu_N}\to \bG_m$ is trivial, there exists a character $\overline{\psi^N}:G/{\pmb \upmu_N}\to \bG_m$ such that $\overline{\psi^N}\circ \pi=\psi^N$.   For a character $\phi: G/{\pmb \upmu_N}\to \bG_m$, we define  characters $\overline{\chi}:G/{\pmb \upmu_N}\to \bG_m$ and $\chi:G\to \bG_m$ by $\overline{\chi}:=\phi\overline{\psi^N}$ and  $\chi:=\overline{\chi}\circ \pi$ respectively. 
We fix a positive integer $m>0$ and a vector ${\bm d}:=(d_1,\hdots,d_m)\in\bZ^m_{>0}$, and  we  define a set 
 \[\cI({\bm d},N):=\{ (i_1,\hdots,i_m)\in\bZ^m_{\geq0} \mid  d_1i_1+\cdots+d_mi_m=N\}.\]
 Let $\bA^m_t$ be the $m$-dimensional affine space with coordinate $t=(t_1,\hdots,t_m)$, and define a $G$-action on $X\times \bA^m_t$ by 
\[
g\cdot(x,t_1,\hdots,t_m):=\Bigl(\pi(g)\cdot x,\psi^{d_1}(g) t_1,\hdots,\psi^{d_m}(g)t_m\Bigr).
\]

Let $W:X\to \bA^1$ be  a $\overline{\chi}$-semi-invariant regular function and $F:X\times\bA^m_t\to\bA^1$ a   non-constant regular function of the form
\[
 F=f_{I_1}t^{I_1}+\hdots +f_{I_r}t^{I_r}\in\Gamma(X,\cO_X)[t_1,\hdots,t_m],
\]
where $I_k=(i_{k,1},\hdots,i_{k,m})\in \cI({\bm d},N)$, $t^{I_k}:=t_1^{i_{k,1}}\cdots t_m^{i_{k,m}}$, and each $f_{I_k}\in\Gamma(X,\cO(\phi))^{G/{\pmb \upmu_N}}$  is a $\phi$-semi-invariant regular function on $X$. Then $F$ is a $\chi$-semi-invariant regular function on $X\times \bA^m_t$ with respect to the above $G$-action.
Denote by $\bP({\bm d}):=[\bA^m_t\!\setminus\!\{\bm 0\}/\bG_m]$ the quotient stack of the $\bG_m$-action on $\bA^m_t\!\setminus\!\{\bm 0\}$ with the weight of $t_i$ being $d_i$, and denote by $Z_F\subset X\times \bP({\bm d})$  the hypersurface  defined by $F$. Then $Z_F$ is isomorphic to the quotient stack $[Y_F/\bG_m]$, where $Y_F\subset Y:=X\times (\bA^m_t\!\setminus\!\{\bm 0\})$ is the zero scheme in $Y$ defined by $F$. Denote by $\cO(1)\in \coh\bP({\bm d})$ the line bundle corresponding to $\cO({\rm id}_{\bG_m})$ via the natural equivalence $\coh\bP({\bm d})\cong \coh_{\bG_m}\bA^m_t\!\setminus\!\{\bm 0\}$, and let $\cL\in \coh_{G/{\pmb \upmu_N}}Z_F$ be the pull-back of $\cO(1)\in \coh_{G/{\pmb \upmu_N}}\bP({\bm d})$ by the morphism $Z_F\to \bP({\bm d})$ defined as the composition of the inclusion $Z_F\hookto X\times \bP({\bm d})$ and the projection $X\times \bP({\bm d})\to \bP({\bm d})$. By abuse of notation, we write $W$ for the pull-backs of $W:X\to \bA^1$ by  the natural projections $X\times\bA^m_t\to X$ and $X\times \bP({\bm d})\to X$.

\begin{thm}\label{main result} Fix $\ell\in\bZ$, and set $\mu:=\sum_{i=1}^md_i$. Assume 
 either of the following conditions:
 \begin{itemize}
 \item[$(i)$] The restricted function  $W|_{Y_F}:Y_{F}\to \bA^1$  is flat.
 \item[$(ii)$] The restricted function  $W|_{Y_F}:Y_{F}\to \bA^1$ is the zero map, $Y_F$ is smooth and the category $\Dcoh_{G/{\pmb \upmu_N}}(Z_{F},\overline{\chi},W|_{Z_F})$ is idempotent complete.
 \end{itemize}
 Then we have the following:
 \begin{itemize}
\item[$(1)$] If $N<\mu$, there are fully faithful functors
\[
\Phi_{\ell}:\Dcoh_{G}(X\times\bA^m_t,\chi,W+F)\hookrightarrow\Dcoh_{G/{\pmb \upmu_N}}(Z_{F},\overline{\chi},W|_{Z_F}),
\]
\[
\Psi_{\ell}:\Dcoh_{G/{\pmb \upmu_N}}(X,\overline{\chi},W)\hookrightarrow\Dcoh_{G/{\pmb \upmu_N}}(Z_{F},\overline{\chi},W|_{Z_F})
\]
and there is a semi-orthogonal decomposition 
\[
\Dcoh_{G/{\pmb \upmu_N}}(Z_{F},\overline{\chi},W|_{Z_F})=\l\, {\sf Im}(\Psi_{\ell}({N-\mu+1})),\hdots,{\sf Im}(\Psi_{\ell}), {\sf Im}(\Phi_{\ell}) \, \r,
\]
where $\Psi_{\ell}(k):=\bigl((-)\otimes\cL^k\bigr)\circ\Psi_{\ell}$.

\vspace{2mm}
\item[$(2)$] If $N=\mu$, we have an equivalence
\[
\Phi_{\ell}:\Dcoh_{G}(X\times\bA^m_t,\chi,W+F)\simto\Dcoh_{G/{\pmb \upmu_N}}(Z_{F},\overline{\chi},W|_{Z_F}).
\]

\vspace{2mm}
\item[$(3)$] If $N>\mu$, 
there are fully faithful functors
\[
\Phi_{\ell}:\Dcoh_{G/{\pmb \upmu_N}}(Z_{F},\overline{\chi},W|_{Z_F})\hookrightarrow\Dcoh_{G}(X\times\bA^m_t,\chi,W+F),
\]
\[
\Psi_{\ell}:\Dcoh_{G/{\pmb \upmu_N}}(X,\overline{\chi},W)\hookrightarrow\Dcoh_{G}(X\times\bA^m_t,\chi,W+F)
\]
and there is a semi-orthogonal decomposition 
\[
\Dcoh_{G}(X\times\bA^m_t,\chi,W+F)=\l\, {\sf Im}(\Psi_{\ell}),\hdots,{\sf Im}(\Psi_{\ell}(\mu-N+1)), {\sf Im}(\Phi_{\ell}) \, \r,
\]
where $\Psi_{\ell}(k):=\bigl((-)\otimes \cO(\psi^k)\bigr)\circ\Psi_{\ell}$.
\end{itemize}
\end{thm}

\subsection{Proof of Theorem \ref{main result}}\label{proof of main result}

We will prove Theorem \ref{main result} by modifying  the arguments in  \cite[Section 3]{degree d} and \cite[Section 5]{H2} rather than by generalizing them. The key new technique is to consider an unramified cyclic cover in \eqref{cover}. Set \[Q:=X\times \bA^m_t\times \bA^1_u.\] We define a $\bG_m\times (G/{{\pmb \upmu_N}})$-action on $Q$  by
\begin{equation}\label{action def}
(a,g)\cdot (x,t_1,\hdots,t_m,u):= (g\cdot x, a^{d_1}t_1,\hdots,a^{d_m}t_m, a^{-N}{\overline{\psi^N}}(g)u).
\end{equation}
Then the regular function \[W_Q:=W+Fu:Q\to \bA^1\] on $Q$ is  $(1\times \overline{\chi})$-semi-invariant. 
Denote by $\lambda:\bG_m\to \bG_m\times (G/{{\pmb \upmu_N}})$ the one-parameter subgroup defined by
$\lambda(a):=(a,1)$. We denote by $Z_{\lambda}\subset Q$ the fixed locus  with respect to the $\lambda$-action.   Then we have  $Z_{\lambda}=\{ (x,t_1,\hdots,t_m,u)\in Q \mid t_1=\cdots=t_m=u=0\}\cong X$. Set $S_+:=\{q\in Q\mid  \displaystyle{\lim_{a\to 0}}\lambda(a)q\in Z_{\lambda}\}$ and 
$S_-:=\{q\in Q\mid  \displaystyle{\lim_{a\to 0}}\lambda(a)^{-1}q\in Z_{\lambda}\}.$ If we define $Q_{\pm}$  to be the complement of $S_{\pm}$  in $Q$,  then  we have 
\[
Q_+=X\times \bA^m_t\times (\bA^1_u\setminus{\{0\}})\hspace{3mm}\mbox{and}\hspace{3mm} Q_{-}=X\times (\bA^m_t\setminus{\{\bm0\}})\times \bA^1_u.
\]
Then the stratifications 
\[
(\K^+:Q=Q_+\sqcup S_+)\hspace{3mm}\mbox{and}\hspace{3mm}  (\K^-: Q=Q_{-}\sqcup S_{-})
\]
are elementary HKKN stratifications, and the pair of these stratifications defines an elementary wall crossing in the sense of \cite{vgit}. Denote by $t(\K^{\pm})$   the $\lambda^{\pm1}$-weight  of the $\bG_m$-action on the fiber 
$\underline{\Spec}({\rm Sym}(\omega_{S_{\pm}/Q}))_x$
of the geometric vector bundle associated to the $\lambda^{\pm1}$-equivariant  relative canonical bundle $\omega_{S_{\pm}/Q}$  at  a point $x$ in the fixed locus $Z_{\lambda}=Z_{\lambda^{-1}}$. Since $Z_{\lambda}$ is connected, the numbers $t(\K^{\pm})$ do not depend on the choice of $x$ (see \cite[Lemma 2.1.18]{vgit}). Then we have 
\[
t(\K^+)=-N \hspace{3mm}\mbox{and}\hspace{3mm}  t(\K^-)=-\mu.
\]
Let $\l\lambda\r\subset \bG_m\times (G/{{\pmb \upmu_N}})$ be the image of $\lambda:\bG_m\to\bG_m\times (G/{{\pmb \upmu_N}})$, and denote by $C(\lambda)$ the centralizer of $\l\lambda\r\subset \bG_m\times (G/{{\pmb \upmu_N}})$. Then we have $C(\lambda)=\bG_m\times (G/{{\pmb \upmu_N}})$. Furthermore, we define $G_{\lambda}:=C(\lambda)/\l\lambda\r$, and let \[\theta:\bG_m\times (G/{{\pmb \upmu_N}})\to \bG_m\] be the character of $\bG_m\times (G/{{\pmb \upmu_N}})=C(\lambda)$ defined by
$\theta(a,g):=a$. Then we have a natural equivalence 
\begin{equation}\label{fixed locus}
\Dcoh_{G_{\lambda}}(Z_{\lambda},1\times\overline{\chi},W_Q|_{Z_{\lambda}})\cong \Dcoh_{G/{{\pmb \upmu_N}}}(X,\overline{\chi},W).
\end{equation}

Let \[i_{\pm}:Q_{\pm}\hookrightarrow Q \hspace{3mm}\mbox{ and }\hspace{3mm} j^{\pm}:S_{\pm}\hookrightarrow Q\] be the open immersions and the closed immersions respectively, and denote by \[\pi_{\pm}:S_{\pm}\to X\]
the natural projections. Then $\pi_+$ is a trivial vector bundle of rank $m$, and $\pi_-$ is a trivial vector bundle of rank $1$. Identifying $X$ with the fixed locus  $Z_{\lambda}=X\times \{\bm 0\}\times\{0\}\subset Q$, we have an induced $\bG_m\times (G/{{\pmb \upmu_N}})$-action on $X$ given by $(a,g)\cdot x=g\cdot x$. Then by Lemma \ref{trivial action} the projection  $p:\bG_m\times (G/{{\pmb \upmu_N}})\to G/{{\pmb \upmu_N}}$  induces a fully faithful functor 
\[
{\rm id}_p^*:\coh_{G/{{\pmb \upmu_N}}}X\hookrightarrow \coh_{\bG_m\times (G/{{\pmb \upmu_N}})}X
\]
since $\bG_m$-action  on $X$ is trivial, and  we have the following decomposition 
\begin{equation}\label{2}
\coh_{\bG_m\times (G/{{\pmb \upmu_N}})}X\cong\bigoplus_{k\in \bZ}(\coh_{\bG_m\times (G/{{\pmb \upmu_N}})}X)_{(\theta|_{\bG_m})^k}.
\end{equation}
We define  $\Dcoh_{\bG_m\times (G/{{\pmb \upmu_N}})}(X,1\times G,W)_k\subset \Dcoh_{\bG_m\times G}(X,1\times G,W)$ to be the subcategory consisting of factorizations $F$ whose components $F_i$ lie in $(\coh_{\bG_m\times G}X)_{(\theta|_{\bG_m})^k}$. Note that the subcategory $\Dcoh_{\bG_m\times (G/{{\pmb \upmu_N}})}(X,1\times G,W)_k$ is the subcategory of factorizations with $\lambda$-weights $k$ (or equivalently, with $(-\lambda)$-weights $-k$) in the sense of \cite[Definition 3.4.3]{vgit}. Then the above   decomposition (\ref{2}) gives rise to a decomposition of the derived factorization category
\[
\Dcoh_{\bG_m\times (G/{{\pmb \upmu_N}})}(X,1\times \overline{\chi},W)\cong\bigoplus_{k\in\bZ}\Dcoh_{\bG_m\times (G/{{\pmb \upmu_N}})}(X,1\times \overline{\chi},W)_k.
\]
Denote by
\[\tau_k:\Dcoh_{G/{{\pmb \upmu_N}}}(X,\overline{\chi},W)\simto \Dcoh_{\bG_m\times (G/{{\pmb \upmu_N}})}(X,1\times \overline{\chi},W)_k\]
 the equivalence defined by $\tau_k(F):= {\rm id}_p^*(F)\otimes \cO(\theta^k)$. We also denote by \[\iota_k:\Dcoh_{\bG_m\times (G/{{\pmb \upmu_N}})}(X,1\times \overline{\chi},W)_k\hookto \Dcoh_{\bG_m\times (G/{{\pmb \upmu_N}})}(X,1\times \overline{\chi},W)\] the natural inclusion.
 
  \begin{dfn}
We define  functors \[\Upsilon^{\pm}_{k}:\Dcoh_{\bG_m\times (G/{{\pmb \upmu_N}})}(X,1\times \overline{\chi},W)_{ k}\to \Dcoh_{\bG_m\times (G/{{\pmb \upmu_N}})}(Q,1\times \overline{\chi}, W_Q)\]
as follows: First consider  the (underived) inverse image
\[
\pi_{\pm}^*:\Dcoh_{\bG_m\times (G/{{\pmb \upmu_N}})}(X,1\times \overline{\chi},W)\to \Dcoh_{\bG_m\times (G/{{\pmb \upmu_N}})}(S_{\pm},1\times \overline{\chi},\pi_{\pm}^*W)
\]
of the  projection map $\pi_{\pm}:S_{\pm}\to X$.
Next, since the pull-back $(j^{\pm})^*W_Q:S_{\pm}\to \bA^1$  equals to the pull-back $\pi_{\pm}^*W$, we have the (underived) direct image 
\[
j^{\pm}_*:\Dcoh_{\bG_m\times (G/{{\pmb \upmu_N}})}(S_{\pm},1\times \overline{\chi},\pi_{\pm}^*W)\to \Dcoh_{\bG_m\times (G/{{\pmb \upmu_N}})}(Q,1\times \overline{\chi},W_Q).
\]
Then we define $\Upsilon_k^{\pm}$ to be the composition $j^{\mp}_*\circ \pi_{\mp}^*\circ \iota_k$.
\end{dfn}

\begin{lem}[{\cite[Lemma 3.4.5]{vgit}}]
The functor \[\Upsilon^{\pm}_{k}:\Dcoh_{\bG_m\times (G/{{\pmb \upmu_N}})}(X,1\times \overline{\chi},W)_k\to \Dcoh_{\bG_m\times (G/{{\pmb \upmu_N}})}(Q,1\times \overline{\chi}, W_Q)\] is fully faithful.
\end{lem}

Next, we consider windows of $\Dcoh_{\bG_m\times (G/{{\pmb \upmu_N}})}(Q,1\times \overline{\chi}, W_Q)$. For an interval $I\subset\bZ$, let 
\begin{equation}\label{I window}
\cW_{\lambda^{\pm1},I}\subset\Dcoh_{\bG_m\times (G/{{\pmb \upmu_N}})}(Q,1\times \overline{\chi}, W_Q)
\end{equation}
be the $I$-grade-windows with respect to $\lambda^{\pm}$ in the sense of  \cite[Definition 3.1.2]{vgit}. 
For each integer $\ell\in \bZ$, we set $I^{\pm}_{\ell}:=[\ell+t(\K^{\pm})+1,\ell]\subset \bZ$, and we define  subcategories
\[
\cW_{\ell}^{\pm}:=\cW_{\lambda^{\pm1},I_{\ell}^{\pm}}\subset \Dcoh_{\bG_m\times (G/{{\pmb \upmu_N}})}(Q,1\times \overline{\chi}, W_Q).
\]
Then, by \cite[Corollary 3.2.2, Proposition 3.3.2]{vgit},  the pull-backs 
\[
i^*_{\pm}(\ell):=i_{\pm}^*|_{\cW_{\ell}^{\pm}}:\cW_{\ell}^{\pm}\to \Dcoh_{\bG_m\times (G/{{\pmb \upmu_N}})}(Q_{\pm},1\times \overline{\chi}, W_Q)
\]
 of the open immersions $i_{\pm}$ are equivalences. Since $\cW^{\pm}_{\ell}=\cW_{{\lambda^{\mp}},[-\ell,-\ell-t(\K^{\pm})-1]}$, we have the following:
 
 \begin{itemize}
\item If $t(\K^+)\leq t(\K^-)$, we have $\cW_{-\ell}^-\subseteq \cW^+_{-t(\K^+)+\ell-1}$. In this case, we define a fully faithful functor 
\[
\Phi_{\ell}:\Dcoh_{\bG_m\times (G/{{\pmb \upmu_N}})}(Q_{-},1\times \overline{\chi}, W_Q)\hookrightarrow\Dcoh_{\bG_m\times (G/{{\pmb \upmu_N}})}(Q_{+},1\times \overline{\chi}, W_Q)
\]
to be the composition $i^*_+({-t(\K^+)+\ell-1})\circ(i^*_{-}({-\ell}))^{-1}$.\vspace{2mm}
\item If $t(\K^+)\geq t(\K^-)$, we have $\cW^+_{-\ell}\subseteq\cW^-_{-t(\K^-)+\ell-1}$. In this case, we define a fully faithful functor
\[
\Phi_{\ell}:\Dcoh_{\bG_m\times (G/{{\pmb \upmu_N}})}(Q_{+},1\times \overline{\chi}, W_Q)\hookrightarrow\Dcoh_{\bG_m\times (G/{{\pmb \upmu_N}})}(Q_{-},1\times \overline{\chi}, W_Q)
\]
to be the composition $i^*_-(-t(\K^-)+\ell-1)\circ(i^*_+(-\ell))^{-1}$
 \end{itemize}
Furthermore, by \cite[Lemma 3.4.6]{vgit}, we have the following.
 \begin{itemize}
 \item If $t(\K^+)<t(\K^-)$ and $-t(\K^-)+\ell\leq j\leq -t(\K^+)+\ell-1$, the essential image of $\Upsilon_j^{+}$ lies in $\cW_{-t(\K^+)+\ell-1}^+$. In this case, we define a fully faithful functor
 \[
 \Psi_j:\Dcoh_{G/{{\pmb \upmu_N}}}(X,\overline{\chi},W)\hookrightarrow\Dcoh_{\bG_m\times (G/{{\pmb \upmu_N}})}(Q_{+},1\times \overline{\chi}, W_Q)
 \]
 to be  the composition $i^*_+({-t(\K^+)+\ell-1})\circ \Upsilon^+_j\circ \tau_j$. Since the essential image of $\Upsilon_{j+1}^{+}$ is equal to  that of the composition $\bigl((-)\otimes \cO(\theta^{})\bigr)\circ\Upsilon_j^{+}$, we have
 \[{\sf Im}(\Psi_{j+1})={\sf Im}\Bigl(\bigl((-)\otimes \cO(\theta^{})\bigr)\circ\Psi_j\Bigr).\]
 \item If $t(\K^+)>t(\K^-)$ and $-t(\K^+)+\ell\leq j\leq -t(\K^-)+\ell-1$, the essential image of $\Upsilon_j^{-}$ lies in $\cW_{-t(\K^-)+\ell-1}^-$. In this case, we define a fully faithful functor
 \[
  \Psi_j:\Dcoh_{G/{{\pmb \upmu_N}}}(X,\overline{\chi},W)\hookrightarrow\Dcoh_{\bG_m\times (G/{{\pmb \upmu_N}})}(Q_{-},1\times \overline{\chi}, W_Q)
 \]
 to be the composition $i^*_-({-t(\K^-)+\ell-1})\circ \Upsilon^-_j\circ\tau_j$. Since the essential image of $\Upsilon_{j+1}^{-}$ is equal to that of the composition $\bigl((-)\otimes \cO(\theta^{-1})\bigr)\circ\Upsilon_j^{-}$, we have
 \[{\sf Im}(\Psi_{j+1})={\sf Im}\Bigl(\bigl((-)\otimes \cO(\theta^{-1})\bigr)\circ\Psi_j\Bigr).\]
\end{itemize}

\vspace{2mm}
Then,   by iterated applications of  \cite[Proposition 3.4.7]{vgit}, the above fully faithful functors give  rise to the following semi-orthogonal decompositions.

\begin{prop}\label{main prop} Fix $\ell\in\bZ$. \begin{itemize}
\item[$(1)$] If $N>\mu$,  there is a semi-orthogonal decomposition 
\[
\Dcoh_{\bG_m\times (G/{{\pmb \upmu_N}})}(Q_{+},1\times \overline{\chi}, W_Q)=\l\, {\sf Im}(\Psi_{N+\ell-1}),\hdots,{\sf Im}(\Psi_{\mu+\ell}), {\sf Im}(\Phi_{\ell}) \, \r.
\]

\vspace{2mm}
\item[$(2)$] If $N=\mu$, we have an equivalence
\[
\Phi_{\ell}:\Dcoh_{\bG_m\times (G/{{\pmb \upmu_N}})}(Q_{+},1\times \overline{\chi}, W_Q)\simto\Dcoh_{\bG_m\times (G/{{\pmb \upmu_N}})}(Q_{-},1\times \overline{\chi}, W_Q).
\]

\vspace{2mm}
\item[$(3)$] If $N<\mu$, 
 there is a semi-orthogonal decomposition 
\[
\Dcoh_{\bG_m\times (G/{{\pmb \upmu_N}})}(Q_{-},1\times \overline{\chi}, W_Q)=\l\, {\sf Im}(\Psi_{\mu+\ell-1}),\hdots,{\sf Im}(\Psi_{N+\ell}), {\sf Im}(\Phi_{\ell}) \, \r.
\]
\end{itemize}
\end{prop}

\vspace{5mm}
Now it is enough to show  the following two lemmas:

\begin{lem}\label{main lemma} With the same notation as above, we have the following.
\begin{itemize}
\item[$(1)$] We have an equivalence
\[
\Phi_+:\Dcoh_{\bG_m\times (G/{{\pmb \upmu_N}})}(Q_+,1\times \overline{\chi}, W_Q)\cong\Dcoh_{G}(X\times\bA^m_t,\chi,W+F).
\]

\item[$(2)$]
 Under the assumption of Theorem \ref{main result}, we have an equivalence
\[
\Phi_-:\Dcoh_{\bG_m\times (G/{{\pmb \upmu_N}})}(Q_{-},1\times \overline{\chi}, W_Q)\cong\Dcoh_{G/{{\pmb \upmu_N}}}(Z_{F},\overline{\chi},W|_{Z_F}).
\]

\end{itemize}
\begin{proof}
(1) We consider  $\bG_m\times G$-action on $Q_+=X\times\bA^m_t\times(\bA^1_u\!\setminus\!\{0\})$ defined by 
\[
(a,g)\cdot (x,t_1,\hdots,t_m,u)\mapsto  (\pi(g)\cdot x, {\psi}(g)^{d_1}t_1,\hdots, {\psi}(g)^{d_m}t_m,{\psi}(g)^{-1}au)
\]
Denote by $\widetilde{Q}_+$ this new $\bG_m\times G$-variety.
Furthermore, we consider an unramified cyclic  cover 
\begin{equation}\label{cover}
p_+:\widetilde{Q}_+\to Q_+
\end{equation}
defined by
$p_+(x,t_1,\hdots,t_m,u):=(x,t_1u^{d_1}\hdots,t_mu^{d_m},u^{-N})$.
Note that $\bG_m\times (G/{{\pmb \upmu_N}})$-action on $Q_+$ lifts to an action from $\bG_m\times G$ on $Q_+$, and the subgroup $1\times {{\pmb \upmu_N}}\subset \bG_m\times G$ acts trivially on $Q_+$. Then $Q_+$ is a $\bG_m\times G$-variety, and $p_+$ is a $\bG_m\times G$-equivariant morphism that is a principal ${\pmb \upmu_N}$-bundle, where ${\pmb \upmu_N}$-action on $\widetilde{Q}_+$ is given by \[\zeta\cdot(x,t_1,\hdots,t_m,u):=(x,\zeta^{d_1}t_1,\hdots,\zeta^{d_m}t_m,\zeta^{-1} u)\]
since $\psi\circ \iota:{\pmb \upmu_N}\to \bG_m$ is the natural inclusion.
Therefore we have the following equivalences
\begin{eqnarray}
\Dcoh_{\bG_m\times (G/{{\pmb \upmu_N}})}(Q_+,1\times \overline{\chi}, W_Q)&\cong&\Dcoh_{\bG_m\times G}(\widetilde{Q}_+,1\times\chi,\pi^*W_Q)\nonumber\\
&\cong&\Dcoh_{ G}(X\times \bA^m_t,\chi, W+F),\nonumber
\end{eqnarray}
where the first equivalence follows from Lemma \ref{invariant} and the second one follows from Lemma \ref{descent}.

(2)  
We define $Y:=X\times (\bA^m_{t}\!\setminus\!\{\bm 0\})$. Then $Q_{-}=Y\times\bA^1_u$, and so $Y$ has the induced $\bG_m\times (G/{\pmb \upmu_N})$-action. Recall that $\chi_{-N}:\bG_m\to\bG_m$ denotes the character defined by $\chi_{-N}(a):=a^{-N}$, and set $\cE:=\cO(\chi_{-N}\times(\phi^{-1}))$. Then $F\in \Gamma(Y,\cE^{\vee})^{\bG_m\times (G/{\pmb \upmu_N})}$, and $Q_{-}$ is the $\bG_m\times (G/{\pmb \upmu_N})$-vector bundle on $Y$  associated to $\cE(1\times\overline{\chi})\cong\cO(\chi_{-N}\times\overline{\psi^N})$. More precisely, we have an isomorphism
\[
Q_{-}\cong \underline{\rm Spec}({\rm Sym}_Y(\cE(1\times\overline{\chi})^{\vee})).
\]  
Then the quotient stack $Z_F=[Y_F/\bG_m]\subset [Y/\bG_m]\cong X\times\bP({\bm d})$ is the closed substack of $X\times\bP({\bm d})$. Therefore, we have the following equivalences 
\begin{eqnarray}
\Dcoh_{\bG_m\times (G/{\pmb \upmu_N})}(Q_{-},1\times \overline{\chi}, W_Q)&\cong&\Dcoh_{\bG_m\times (G/{\pmb \upmu_N})}(Y_{F},1\times\overline{\chi}, W|_{Y_{F}})\nonumber\\
&\cong&\Dcoh_{G/{\pmb \upmu_N}}(Z_F,\overline{\chi},W|_{Z_F}),\nonumber
\end{eqnarray}
where the first equivalence follows from Theorem \ref{global knorrer}.
\end{proof}
\end{lem}

\begin{lem}
Let $\Phi_{\pm}$ be the  equivalences in the above lemma. We have the following  functor isomorphisms:
\begin{eqnarray}
\Phi_+\Bigl((-)\otimes \cO(\theta)\Bigr)&\cong&  \Phi_+(-)\otimes\cO(\psi)\nonumber\\
\Phi_-\Bigl((-)\otimes \cO(\theta)\Bigr)&\cong &\Phi_-(-)\otimes\cL\nonumber
\end{eqnarray}

\begin{proof}
Denote by $\gamma:\bG_m\times G\to \bG_m\times (G/{{\pmb \upmu_N}})$ and $\eta:G\to \bG_m\times G$ the morphisms defined by $\gamma(a,g):=(a,\pi(g))$ and $\eta({g}):=(\psi({g}),{g})$, respectively. Let $e:X\times \bA^m_t\to\widetilde{Q}_+=X\times \bA^m_t\times (\bA^1_u\setminus \{0\}) $ be the inclusion defined by $e(x,t):=(x,t,1)$.  Then $\Phi_+=e_{\eta}^*\circ{(p_+)}_{\gamma}^*$ by definition, and hence the former functor isomorphism  follows from the following sequence of isomorphisms;
\begin{eqnarray}
\Phi_+\Bigl((-)\otimes \cO(\theta)\Bigr)&\cong& e_{\eta}^*\Bigl({(p_+)}_{\gamma}^*\bigl((-)\otimes \cO(\theta)\bigr)\Bigr)\nonumber\\
&\cong& e_{\eta}^*\Bigl({(p_+)}_{\gamma}^*(-)\otimes \cO(\theta\circ \gamma)\Bigr)\nonumber\\
&\cong& e_{\eta}^*\Bigl({(p_+)}_{\gamma}^*(-)\Bigr)\otimes \cO(\theta\circ \gamma\circ \eta)\nonumber\\
&=& \Phi_+(-)\otimes \cO(\psi).\nonumber
\end{eqnarray}
The latter functor isomorphism is also checked by an easy direct calculation.
\end{proof}
\end{lem}

\subsection{One-dimensional reduction of  factorizations}
In this subsection, we refine the semi-orthogonal decompositions in  Theorem \ref{main result}  in the case when $m=d_1=\mu=1$. We keep the notation above, and let $s\in \Gamma(X,\cO(\phi))^{G/{\pmb \upmu_N}}$ be a non-zero $\phi$-semi-invariant regular function on $X$. Denote by $Z_s\subset X$  the zero scheme of $s$. Then if $s$ is non-constant and $m=d_1=1$, $X\times \bP({\bm d})$ is isomorphic to $X$ and the algebraic stack $Z_F$ is  isomorphic to  $Z_s$. Hence by Proposition \ref{main prop}, Lemma \ref{main lemma}.(2) and  Lemma \ref{zero cat} below, we obtain a semi-orthogonal decomposition that compares derived factorization categories associated to the potentials $W:X\to \bA^1$ and $W+st^N:X\times \bA^1_t\to \bA^1$.

\begin{cor}\label{1d reduction}  There is a fully faithful functor
\[
\Phi:\Dcoh_{G/{\pmb \upmu_N}}(X,\overline{\chi},W)\hookrightarrow\Dcoh_{G}(X\times\bA^1_t,\chi,W+st^N).\]
Furthermore, we have the following semi-orthogonal decomposition:
\begin{itemize}
\item[$(1)$] 
If $s=1$ and $N>1$, we have a semi-orthogonal decomposition
\[
\Dcoh_{G}(X\times\bA^1_t,{\chi},W+t^N)=\langle\,{\sf Im}(\Phi_0),\hdots,{\sf Im}(\Phi_{-N+2})  \,\rangle.
\]
\item[$(2)$] If $s$ is non-constant and  $W|_{Z_s}:Z_s\to \bA^1$ is flat,  there is a fully faithful functor 
\[
\Psi:\Dcoh_{G/{\pmb \upmu_N}}(Z_s,\overline{\chi},W|_{Z_s})\hookrightarrow\Dcoh_{G}(X\times\bA^1_t,\chi,W+st^N),
\]
and we have a semi-orthogonal decomposition
\[
\Dcoh_{G}(X\times\bA^1_t,\chi,W+st^N)=\langle\,{\sf Im}(\Phi_0),\hdots,{\sf Im}(\Phi_{-N+2}), {\sf Im}(\Psi)\,\rangle.
\]
\end{itemize}
Here $\Phi_i$ denotes  the composition $((-)\otimes\cO(\psi^{i}))\circ\Phi$.
\end{cor}

\begin{lem}\label{zero cat} 
 If $s=1$, the category $\Dcoh_{\bG_m\times (G/{{\pmb \upmu_N}})}(Q_{-},1\times \overline{\chi}, W_Q)$ in the proof of Theorem \ref{main result} is the zero category. 

\begin{proof}
 Since $\bG_m\times (G/{{\pmb \upmu_N}})$ is reductive, by \cite[Lemma 2.33]{H2} the restriction functor 
\[
\Dcoh_{\bG_m\times (G/{{\pmb \upmu_N}})}(Q_{-},1\times \overline{\chi}, W_Q)\to \Dcoh(Q_{-}, W_Q)
\]
is faithful. Hence it is enough to show that $\Dcoh(Q_{-}, W_Q)=0$. Since $W_Q$ is flat, $\Dcoh(Q_{-}, W_Q)$ is equivalent to the singularity category
$\Dsg(Q_0)$ of the zero scheme $Q_0$ of $W_Q:Q_{-}\to \bA^1$ by Theorem \ref{fact sing}. Since $Q_0$ is smooth, we have  $\Dcoh(Q_{-}, W_Q)\cong\Dsg(Q_0)=0$.
\end{proof}
\end{lem}

\begin{rem}
By \cite{Tak} or \cite{orlov3}, there are quasi-fully faithful dg-functors $F_i: {\sf inj}_{\bG_m}(\Spec k,\chi_1,0)\to {\sf inj}_{\bG_m}(\bA^1_t,\chi_N,t^N)$ ($-N+2\leq i\leq 0$) and  the semi-orthogonal decomposition
\[
{\sf inj}_{\bG_m}(\bA^1_t,\chi_N,t^N)=\langle\, {\sf Im}(F_0),\hdots,{\sf Im}(F_{-N+2})\,\rangle.
\]
Moreover   \cite[Corollary 5.18]{bfk} implies  the following semi-orthogonal decomposition 
\[
\overline{{\sf inj}}_{G/{\pmb \upmu_N}}(X,\overline{\chi},W)\cong{\sf inj}_{G/{\pmb \upmu_N}}(X,\overline{\chi},W)\circledast{\sf inj}_{\bG_m}(\Spec k,\chi_1,0),
\]
where $\overline{{\sf inj}}$ denotes the dg-subcategory of ${\sf Inj}$ consisting of compact objects. Hence by \cite[Corollary 5.18]{bfk} and \cite[Lemma 2.49]{bfk2} there are quasi-fully faithful functors 
\[
\widetilde{\Phi}_i:\overline{{\sf inj}}_{G/{\pmb \upmu_N}}(X,\overline{\chi},W)\to\overline{{\sf inj}}_G(X\times \bA^1_t,\chi,W+t^N)
\] for $-N+2\leq i\leq 0$ and the semi-orthogonal decomposition
\[
\overline{{\sf inj}}_G(X\times \bA^1_t,\chi,W+t^N)=\langle\, {\sf Im}(\widetilde{\Phi}_0),\hdots,{\sf Im}(\widetilde{\Phi}_{N-2})\,\rangle.
\]
Hence Corollary \ref{1d reduction}.(1) follows from the derived Morita theory if $\Dcoh_{G}(X\times\bA^1_t,{\chi},W+t^N)$ and $\Dcoh_{G/{\pmb \upmu_N}}(X,\overline{\chi},W)$ are idempotent complete.
\end{rem}

%%%%%%%%%%%%%%%%%%%%%%%%%%%%%%%%%%%%%%%%%%%

\subsection{Semi-orthogonal decomposition for chain polynomials}\label{chain sod section}
For each positive integer $i\in\bZ_{\geq1}$, choose a positive integer $a_i\in\bZ_{\geq 1}$ with $a_1\geq2$. Then for each $n\in \bZ_{\geq0}$, we have the associated  polynomial \[f_n:=x_1^{a_1}+x_1x_2^{a_2}+\cdots x_{n-1}x_n^{a_n},\]
where we set $f_0:=0$.
We denote by $L_n$  {\it the maximal grading  of $f_n$}, which is defined  in \eqref{max intro}. 
Then we have  $L_0=\bZ\vec{0}\cong \bZ$, $L_1\cong \bZ$ and $\rk(L_n)=1$ for any $n\geq0$.  
Then the $n$-dimensional polynomial ring $S^n=k[x_1,\hdots,x_n]$ has a natural $L_n$-grading such that $\deg(x_i):=\vec{x}_i$ and $\deg(f)=\vec{f}$. For $n\geq-1$, we set 
\begin{eqnarray*}
\cC_n&:=&\HMF_{S^n}^{L_n}(f_{n})\\
\cD_n&:=&\Dmod^{L_n}(S^n,f_n)
\end{eqnarray*}
where $\cC_{-1}$ and $\cD_{-1}$ are defined to be the zero categories. It is well known that the category $\cC_n$ is idempotent complete for any $n\geq0$ \cite[Proposition 2.7]{at}.

Any element $l\in L_n$ can be represented by $\sum_{i=1}^nl_i\vec{x}_i$ for some $l_i\in \bZ$, and it induces a character  $\chi_l\colon G_n\to \bG_m$ defined by 
\[
 \chi_l(\lambda_1,\hdots,\lambda_n) :=\prod_{i=1}^n\lambda_i^{l_i}.
\]
Note that this character does not depend on the choice of $l_i$, and this defines an isomorphism
\begin{equation}\label{L=G}
L_n\simto (G_n)^{\vee}=\Hom(G_n,\bG_m),
\end{equation}
where $G_n$ is the algebraic group defined in \eqref{max sym intro}.
Thus, if we write
\[
\cD_n^{\coh}:=\Dcoh_{G_n}(\bA^n,\chi_{f_n},f_n)
\]
as in the introduction, by Proposition \ref{coh mod}, we have  natural equivalences
\begin{equation}\label{C=D}
\cC_n\simto\cD_n\simto\cD_n^{\coh}.
\end{equation}

Now we apply  Theorem \ref{main result} in the case when $\ell=a_n-1$ for chain polynomials as follows:
We use the same notation as in Section \ref{proof of main result}, and we consider the case when  $m=\mu=d_1=1$, $X=\bA^{n-1}_x$, $N=a_n$ and $G=G_n$. 
The natural projection 
$\pi_n:(\bG_m)^n\to (\bG_m)^{n-1}$
induces a surjective  morphism 
\begin{equation}\label{pi def}
\pi:=\pi_n:G_n\to G_{n-1}; (\lambda_1,\hdots,\lambda_n)\mapsto (\lambda_1,\hdots,\lambda_{n-1})
\end{equation}
of algebraic groups, and we have $\Ker(\pi)=\{(1,\hdots,1,\lambda_n) \in G_n\mid \lambda_d\in \pmb \upmu_{a_n}\}\cong \pmb \upmu_{a_n}$. Hence there is an exact sequence
\[
1\to {\pmb \upmu_{a_n}}\xrightarrow{\iota} G_n\xrightarrow{\pi} G_{n-1}\to 1,
\]
where $\iota:{\pmb \upmu_{a_n}}\to  G_n$ is defined by $\iota(\lambda):=(1,\hdots,1,\lambda)$, and so $G_{n-1}\cong G/\upmu_{a_n}$.
Let $\psi:G=G_n\to \bG_m$, $\overline{\psi^{a_n}}:G/\upmu_{a_n}\cong G_{n-1}\to \bG_m$ and $\phi:G/{\pmb\upmu_{a_n}}\cong G_{n-1}\to \bG_m$ be the characters  in Section \ref{proof of main result} defined by
\begin{eqnarray}
\psi(\lambda_1,\hdots,\lambda_n)&:=&\lambda_n\nonumber\\
\overline{\psi^{a_n}}(\lambda_1,\hdots,\lambda_{n-1})&:=&\lambda_1^{a_1}\lambda_{{n-1}}^{-1}\nonumber\\
\phi(\lambda_1,\hdots,\lambda_{n-1})&:=&\lambda_{n-1}.\nonumber
\end{eqnarray}
respectively. Then $\psi\circ\iota:{\pmb \upmu_{a_n}}\to \bG_m$ is the natural inclusion, and $\overline{\psi^{a_n}}\circ\pi=\psi^{a_n}$. Recall that  $\overline{\chi}:G_{n-1}\to \bG_m$ is the character defined by $\overline{\chi}:=\phi\overline{\psi^{a_n}}$ and  $\chi=\overline{\chi}\circ \pi:G_n\to \bG_m$. Then $\overline{\chi}=\chi_{f_{n-1}}$ and $\chi=\chi_{f_n}$.
We write $x_n$ instead of $t$, and so 
\[
 Q=\bA^{n-1}_x\times \bA^1_{x_n}\times \bA^1_u\cong \Spec S^{n-1}[x_n,u],
\]
and $Q$ has an action from  $\bG_m\times G/{\pmb\upmu_{a_n}}\cong \bG_m\times G_{n-1}$ defined by \eqref{action def}.
We set 
\begin{eqnarray}
R\,\,\,\,&:=&S^{n-1}[x_n,u]\\
R_+&:=&S^{n-1}[x_n,u^{\pm1}]\\
R_-&:=&S^{n-1}[x_n^{\pm1},u].
\end{eqnarray}
Then $Q\cong \Spec R$, and  we have
\begin{eqnarray}
Q_+&= &\bA^{n-1}_x\times \bA^1_{x_n}\times (\bA^1_u\setminus\{0\})\cong \Spec R_+\\
Q_-&= &\bA^{n-1}_x\times (\bA^1_{x_n}\setminus\{0\})\times \bA^1_u\cong \Spec R_-.
\end{eqnarray}
Consider the case when
  $W=f_{n-1}$ and  $F=x_{n-1}x_n^{a_n}$, and thus $W_Q=f_{n-1}+x_{n-1}x_n^{a_n}u$.  
 By Proposition \ref{main  prop} in the case when $\ell:=-a_n+1$, the fully faithful functors
  \begin{equation}\label{Psi i}
  \Psi_i:\cD_{n-1}^{\coh}\hookto\Dcoh_{\bG_m\times G_{n-1}}(Q_+,1\times\chi_{f_{n-1}},W_Q)
  \end{equation}
  for $-a_n+2\leq i\leq 0$ and
\begin{equation*}\label{Phi l}
\Phi_{(-a_n+1)}:\Dcoh_{\bG_m\times G_{n-1}}(Q_-,1\times\chi_{f_{n-1}},W_Q)\hookto\Dcoh_{\bG_m\times G_{n-1}}(Q_+,1\times\chi_{f_{n-1}},W_Q),
\end{equation*}
gives rise to a semi-orthogonal decomposition
\begin{equation*}\label{sod Q}
\Dcoh_{\bG_m\times G_{n-1}}(Q_{+},1\times \chi_{f_{n-1}}, W_Q)=\l\, {\sf Im}(\Psi_{0}),\hdots,{\sf Im}(\Psi_{-a_n+2}), {\sf Im}(\Phi_{-a_n+1}) \, \r.
\end{equation*}
Moreover by Lemma \ref{main lemma} we have equivalences
\begin{eqnarray*}
\Phi_+:\Dcoh_{\bG_m\times G_{n-1}}(Q_+,1\times \chi_{f_{n-1}},W_Q)&\simto& \cD_{n}^{\coh}\\
\Phi_-:\Dcoh_{\bG_m\times G_{n-1}}(Q_-,1\times \chi_{f_{n-1}},W_Q)&\simto&
\Dcoh_{G_{n-1}}(\bA^{n-2}_x,\chi_{f_{n-1}},f_{n-2}).
\end{eqnarray*}
Since the subgroup ${\pmb\upmu_{a_n}}\hookto G_{n-1};\lambda\mapsto(1,\hdots,1,\lambda)$ of $G_{n-1}$ trivially acts on $\bA^{n-2}_x$, by Lemma \ref{trivial action} we have orthogonal decomposition
\[
\coh_{G_{n-1}}\bA^{n-2}_x\cong\bigoplus_{j=0}^{a_{n-1}-1}{\sf Im}\Bigl(\bigl((-)\otimes \cO(\phi^j)\bigr)\circ \id_{\pi_{n-1}}^*\Bigr), 
\]
    where $\id_{\pi_{n-1}}^*:\coh_{G_{n-2}}\bA^{n-2}_x\hookto\coh_{G_{n-1}}\bA^{n-2}_x$ is the fully faithful functor associated to the surjection $\pi_{n-1}:G_{n-1}\to G_{n-2}$. Since  $\phi^{j+a_{n-1}}|_{\pmb\upmu_{a_{n-1}}}=\phi^j|_{\pmb\upmu_{a_{n-1}}}$ for any $j\in\bZ$, we have an equality ${\sf Im}\Bigl(\bigl((-)\otimes \cO(\phi^{j+a_{n-1}})\bigr)\circ \id_{\pi_{n-1}}^*\Bigr)={\sf Im}\Bigl(\bigl((-)\otimes \cO(\phi^j)\bigr)\circ \id_{\pi_{n-1}}^*\Bigr)$  by Lemma \ref{trivial action}. Thus we have 
\[
\coh_{G_{n-1}}\bA^{n-2}_x\cong\bigoplus_{j=0}^{-a_{n-1}+1}{\sf Im}\Bigl(\bigl((-)\otimes \cO(\phi^j)\bigr)\circ \id_{\pi_{n-1}}^*\Bigr).
\]

We define the functors 
\begin{eqnarray}
\Uppsi^{\coh}:\cD_{n-1}^{\coh}&\to& \cD_n^{\coh}\\
\Upphi^{\coh}:\cD_{n-2}^{\coh}&\to& \cD_n^{\coh}
\end{eqnarray}
to be the following compositions of  fully faithful functors
\begin{eqnarray}
\Uppsi^{\coh}&:=&\Phi_+\circ \Psi_0 \\
\Upphi^{\coh}&:=& \Phi_+\circ\Phi_{(-a_n+1)}\circ (\Phi_-)^{-1}\circ \id_{\pi_{n-1}}^* .
\end{eqnarray}
Denote by 
 \begin{alignat}{2}
 \Uppsi&:\cC_{n-1}\hookto\cC_n,&\quad\quad
 \Upphi&:\cC_{n-2}\hookto\cC_n,\\
  \Uppsi^{\fmod}&:\cD_{n-1}\hookto\cD_n,&\quad\quad
 \Upphi^{\fmod}&:\cD_{n-2}\hookto\cD_n,
 \end{alignat}
the corresponding functors via the natural equivalences \eqref{C=D}.

\begin{thm}\label{sod for chain} Let $n\geq1$ be a positive integer.
The fully faithful functors
 \[
 \Uppsi:\cC_{n-1}\hookto\cC_n\mbox{ \hspace{3mm} and \hspace{3mm}}
 \Upphi:\cC_{n-2}\hookto\cC_n,
 \]
gives a semi-orthogonal decomposition
\[
\cC_n={\Bigg \langle}{\sf Im}(\Uppsi_0),\hdots,{\sf Im}(\Uppsi_{(-a_{n}+2)}),\bigoplus_{j=0}^{{-a_{n-1}+1}}{\sf Im}(\Upphi_{j}){\Bigg \rangle},
\]
 where $\Uppsi_i$ is the composition  $\bigl((-)\otimes (i\vec{x}_n)\bigr)\circ\Uppsi$ and $\Upphi_j$ is the composition $\left((-)\otimes(j\vec{x}_{n-1})\right)\circ\Upphi$.
\begin{proof}
 By the above argument, we have a semi-orthogonal decomposition
\[
\cD_n^{\coh}={\Bigg \langle}{\sf Im}(\Uppsi^{\coh}_0),\hdots,{\sf Im}(\Uppsi^{\coh}_{(-a_{n}+2)}),\bigoplus_{j=0}^{{-a_{n-1}+1}}{\sf Im}(\Upphi^{\coh}_{j}){\Bigg \rangle},
\]
where $\Uppsi^{\coh}_i$ is the composition  $\bigl((-)\otimes \cO(\psi^i)\bigr)\circ\Uppsi^{\coh}$ and $\Upphi^{\coh}_j$ is the composition $\left((-)\otimes\cO(\phi^j)\right)\circ\Upphi^{\coh}$. Since the character $\psi:G_n\to \bG_m$ (resp. $\phi:G_{n-1}\to \bG_m$) corresponds to the element $\vec{x}_n\in L_n$ (resp. $\vec{x}_{n-1}\in L_{n-1}$) via the isomorphism \eqref{L=G},  the result follows from the equivalence 
$\cD_n^{\coh}\cong \cC_n$ in \eqref{C=D}.
\end{proof}
\end{thm}

\section{Full strong exceptional collections for chain polynomials}

In this section, we prove Conjecture \ref{fsec} for chain polynomials. We use the same notation as in Section \ref{chain sod section}.

\subsection{Inductive construction of full strong exceptional collections}\label{explicit matrix}
In this section, we construct a sequence $\cE^n$ of objects in $\cC_n$, which turn out to be full strongly exceptional. Set 
\begin{equation}\label{E def}
E^0:=(0\to k\to 0)\in\cC_0.
\end{equation}
We define two objects $\uppsi E^0\in \cC_1$ and $\upphi E^0\in \cC_2$
by 
\begin{eqnarray*}
\uppsi E^0&:=&\left( S^1(-\vec{x}_1)\xrightarrow{x_1} S^1\xrightarrow{x_1^{a_1-1}} S^1(\vec{f}-\vec{x}_1)\right)\\
\upphi E^0 &:=&\left( S^2(-\vec{x}_1)\xrightarrow{x_1} S^2\xrightarrow{x_1^{a_1-1}+x_2^{a_2}} S^1(\vec{f}-\vec{x}_1)\right),
\end{eqnarray*}
and for each $0\leq i\leq a_1-2$ and $0\leq j\leq a_1-1$, we set 
\begin{eqnarray*}
\hspace{37mm}\uppsi_i E^0&:=&\uppsi E^0(-i\vec{x}_1)[i]\in \cC_1\\
\hspace{37mm}\upphi_j E^0 &:=& \upphi E^0(-j\vec{x}_1+(-a_2+1)\vec{x}_2)[a_2+j-1]\in \cC_2.
\end{eqnarray*}

Similarly, for $n\geq1$ and an object $F=\left(F_1\xrightarrow{\varphi_1}F_0\xrightarrow{\varphi_0}F_1(\vec{f}\,)\right)\in\cC_n$, we define  objects 
\[
\uppsi F\in \cC_{n+1} \mbox{\hspace{5mm}and\hspace{5mm}} \upphi F\in\cC_{n+2}
\] 
as follows: To ease notation, we set  
\begin{equation}\label{easy notation}
x:=x_{n}, \hspace{2mm}y:=x_{n+1}, \hspace{2mm}z:=x_{n+2}, \hspace{2mm}b:=a_{n+1}, \mbox{ and } \hspace{2mm}c:=a_{n+2},
\end{equation}
 and so  $f_{n+1}=f_n+xy^b$ and $f_{n+2}=f_n+xy^b+yz^c$. The objects $\uppsi F$  and $\upphi F$ are defined to be the matrix factorizations defined by

{\scriptsize
\begin{eqnarray*}\uppsi F&:=&\left(
\widetilde{F}_1\oplus\left(\widetilde{F}_0(-\vec{y})\right)
\xrightarrow{\left(\begin{array}{cc}\widetilde{\varphi}_1&y\\-xy^{b-1}&\widetilde{\varphi}_0\end{array}\right)}
\widetilde{F}_0\oplus\left(\widetilde{F}_1(\vec{f}-\vec{y})\right)
\xrightarrow{\left(\begin{array}{cc}\widetilde{\varphi}_0&-y\\xy^{b-1}&\widetilde{\varphi}_1\end{array}\right)}
\widetilde{F}_1(\vec{f})\oplus\left(\widetilde{F}_0(\vec{f}-\vec{y})\right)
\right)\\
\upphi F&:=&\left(
\widehat{F}_1\oplus\left(\widehat{F}_0(-\vec{y})\right)
\xrightarrow{\left(\begin{array}{cc}\widehat{\varphi}_1&y\\-g&\widehat{\varphi}_0\end{array}\right)}
\widehat{F}_0\oplus\left(\widehat{F}_1(\vec{f}-\vec{y})\right)
\xrightarrow{\left(\begin{array}{cc}\widehat{\varphi}_0&-y\\g&\widehat{\varphi}_1\end{array}\right)}
\widehat{F}_1(\vec{f})\oplus\left(\widehat{F}_0(\vec{f}-\vec{y})\right)
\right)
\end{eqnarray*}}
\begin{flushleft}where $(\widetilde{-}):=(-)^{L_{n+1}}\otimes_{S^n}S^{n+1}$, $(\widehat{-}):=(-)^{L_{n+2}}\otimes_{S^n}S^{n+2}$ and $g:=xy^{b-1}+z^c$. Sometimes we simply write $\varphi_i$ for $\widetilde{\varphi}_i$ or $\widehat{\varphi}_i$. For each $0\leq i\leq b-2$ and each $0\leq j\leq b-1$, we set 
\end{flushleft}
\begin{eqnarray}\label{psi phi obj}
\hspace{35mm}\uppsi_iF\hspace{-1mm}&:=&\hspace{-1mm}\uppsi F(-i\vec{y})[i]\in\cC_{n+1}\\
\hspace{35mm}\upphi_jF\hspace{-1mm}&:=&\hspace{-1mm}\upphi F\left(-j\vec{y}+(-c+1)\vec{z}\right)[c+j-1]\in\cC_{n+2}.
\end{eqnarray}
Then $\uppsi_0F=\uppsi F$ and $\upphi_0F=\upphi F\left((-c+1)\vec{y}\right)[c-1]$.

For  a sequence  $\cE=(E_1,\hdots,E_r)$ of objects in $\cC_n$ and for each $0\leq i\leq b-2$, we define a sequence $\uppsi_i\cE$ of objects in $\cC_{n+1}$ by 
\[
\uppsi_i\cE:=(\uppsi_iE_1,\hdots, \uppsi_iE_r).
\]
Similarly, for a sequence  $\cF=(F_1,\hdots,F_s)$ of objects in $\cC_n$ and for $0\leq j\leq b-1$, we define a sequence $\upphi_j\cF$ of objects in $\cC_{n+2}$ by
\[
\upphi_j\cF:=(\upphi_jF_1,\hdots,\upphi_jF_s).
\]
We inductively define sequences $\cE^n$ of objects in $\cC_n$ for any $n\geq-1$ as follows.
First set 
\[
\cE^{-1}:=\emptyset\mbox{\hspace{5mm}and\hspace{5mm}}
\cE^0:=\{E^0\}.
\]
Then we define the sequence $\cE^n$ by
\[
\cE^{n}:=\left(\uppsi_0\cE^{n-1},\hdots,\uppsi_{(a_n-2)}\cE^{n-1},\upphi_0\cE^{n-2},\hdots,\upphi_{(a_{n-1}-1)}\cE^{n-2}\right).
\]
The following is our main result in this paper.

\begin{thm}\label{strong thm}
For any $n\geq1$, the sequence $\cE^n$ is a full strong exceptional collection in $\cC_n$, and, if $a_n\geq2$, the length of $\cE^n$ is equal to the Milnor number  of $\widetilde{f}_n$. In particular, $\cC_n$ has a tilting object.
\end{thm}

In the study of Cohen-Macaulay representation, the existence of a tilting object in the stable category of graded Cohen-Macaulay modules over a graded Gorenstein ring is a fundamental   problem \cite[Problem 3.4]{iyama}. By Theorem \ref{strong thm} and Theorem \ref{mf=cm}, we have the following. 

\begin{cor}
The  category  $\underline{\CM}^{L_n}(S^n/f_n)$ has a tilting object. 
\end{cor}

\subsection{Explicit descriptions of  $\Uppsi$ and $\Upphi$}
In this section, we  describe the functors $\Uppsi:\cC_{n-1}\to\cC_n$ and $\Upphi:\cC_{n-2}\to \cC_n$ in Theorem \ref{sod for chain} explicitly in terms of matrix factorizations. This descriptions enable us to apply  certain distinguished triangles in $\cC_n$ (see Lemma \ref{psi phi tri} below) to our proof of Theorem \ref{strong thm}. We use the notation in the above sections, e.g. Section \ref{proof of main result}, \ref{chain sod section}.
\begin{lem}\label{psi coh des}
The fully faithful functor 
$
\Uppsi^{\coh}:\cD_{n-1}^{\coh}\to\cD_n^{\coh}
$
is isomorphic to the composition 
\[
\cD_{n-1}^{\coh}\xrightarrow{\id_{\pi_{n}}^*}\Dcoh_{G_n}(\bA^{n-1}_x,\chi_{f_{n}},f_{n-1})\xrightarrow{j_*}\cD_n^{\coh}
\]
where $\pi_n:G_n\to G_{n-1}$ is the surjection in \eqref{pi def} and $j:\bA_x^{n-1}\hookto \bA^n_x;(x_1,\hdots,x_{n-1})\mapsto (x_1,\hdots,x_{n-1},0)$ is the natural closed immersion.
\begin{proof}

By construction the functor $\Psi_0$ in \eqref{Psi i} is given by  
\[
\Psi_0\left(F_1\xrightarrow{\varphi_1} F_0\xrightarrow{\varphi_0}  F_1(\chi_{f_{n-1}})\right)=\left(\Uptheta(F_1)\xrightarrow{\Uptheta(\varphi_1)} \Uptheta(F_0)\xrightarrow{\Uptheta(\varphi_0)}  \Uptheta(F_1)(1\times\chi_{f_{n-1}})\right),
\]
where $\Uptheta: \coh_{G_{n-1}}\bA^{n-1}_x \to \coh_{\bG_m\times G_{n-1}}Q$ is the exact functor defined to be the compositions $i_+^*\circ j^-_*\circ \pi_-^*\circ \id_{p_2}^*$ of the following exact functors:
\begin{eqnarray*}
\coh_{G_{n-1}}\bA^{n-1}_x&\xrightarrow{\quad\id_{p_2}^*\quad}& \coh_{\bG_m\times G_{n-1}}\bA^{n-1}_x\\
&\xrightarrow{\,\quad\pi_-^*\quad\,}&\coh_{\bG_m\times G_{n-1}}S_-\\
&\xrightarrow{\,\quad j^-_*\quad\,}&\coh_{\bG_m\times G_{n-1}}Q\\
&\xrightarrow{\,\quad i_+^*\quad\,}&\coh_{\bG_m\times G_{n-1}}Q_+.
\end{eqnarray*}
By the flat base change, the functor $\Uptheta$ is isomorphic to the composition $q^*\circ j_*\circ \id_{p_2}^*$, where 
\[
j_*: \coh_{\bG_m\times G_{n-1}}\bA^{n-1}_x\to \coh_{\bG_m\times G_{n-1}}\bA^{n-1}_x\times\bA^1_{x_n}
\]
 is the direct image by the natural closed immersion $j:\bA^{n-1}_x\hookto \bA^{n-1}_x\times\bA^1_{x_n}$ and 
\[
q^*:\coh_{\bG_m\times G_{n-1}}\bA^{n-1}_x\times\bA^1_{x_n}\to \coh_{\bG_m\times G_{n-1}}Q_+
\]
is the pull-back by the natural projection $q:Q_+=\bA^{n-1}_x\times\bA^1_{x_n}\times(\bA^1_u\setminus\{0\})\to \bA^{n-1}_x\times\bA^1_{x_n}$.

On the other hand, recall that the equivalence 
\[\Phi_+:\Dcoh_{\bG_m\times G_{n-1}}(Q_{+},1\times \chi_{f_{n-1}}, W_Q)\simto \Dcoh_{G_n}(\bA^n_x,\chi_{f_n},f_n)\]
 in Lemma \ref{main lemma} is defined by
\[
\Phi_+\left(F_1\xrightarrow{\varphi_1} F_0\xrightarrow{\varphi_0}  F_1(1\times \chi_{f_{n-1}})\right)=\left(\Upsigma(F_1)\xrightarrow{\Upsigma(\varphi_1)} \Upsigma(F_0)\xrightarrow{\Upsigma(\varphi_0)}  \Upsigma(F_1)(\chi_{f_{n}})\right),
\]
where the functor $\Upsigma:\coh_{\bG_m\times G_{n-1}}Q_{+}\simto  \coh_{G_n}\bA^n_x$ is the exact equivalence defined by the composition
\begin{eqnarray*}
\left((p_{+*})^{\pmb\upmu_{a_n}}\right)^{-1}:\coh_{\bG_m\times G_{n-1}}Q_{+}&\simto& \coh_{\bG_m\times G_n}\widetilde{Q}_{+}\\
\widetilde{e}_{\varphi}^*:\coh_{\bG_m\times G_n}\widetilde{Q}_{+}&\simto&\coh_{G_n}\bA^n_x,
\end{eqnarray*}
where $\widetilde{Q}_+$ and $p_+:\widetilde{Q}_+\to Q$ are as in the proof of Lemma \ref{main lemma}, $\widetilde{e}:\bA^n_x\hookto \widetilde{Q}_+=\bA^n_x\times(\bA^1_u\setminus\{0\})$ is the embedding given by $x \mapsto (x,1)$ and $\varphi:G_n\to \bG_m\times G_n$ is the character defined by $\varphi(g):=(\psi(g),g)$. By \cite[Corollary 2.24]{bfk}, the equivalence $\left((p_{+*})^{\pmb\upmu_{a_n}}\right)^{-1}$ is isomorphic to $(p_+)_{1\times\pi}^*$, and thus $\Upsigma$ is isomorphic to $(p_+\circ \widetilde{e})^*_{(1\times\pi)\circ \varphi}\cong e_{\psi\times \pi}^*$, where $e:\bA^n_x\hookto Q_+=\bA^n_x\times(\bA^1_u\setminus\{0\})$ is the map $x\mapsto (x,1)$ and $\psi\times \pi:G_n\to \bG_m\times G_{n-1}$ is given by $(\psi\times\pi)(\lambda_1,\hdots,\lambda_n)=\left(\lambda_n,(\lambda_1,\hdots,\lambda_{n-1})\right)$. Therefore, using \cite[Lemma 2.19]{bfk}, we see that $\Uppsi^{\coh}$ is induced by the composition
\[
\coh_{G_{n-1}}\bA^{n-1}_x\xrightarrow{\id_{\pi}^*} \coh_{G_n}\bA^{n-1}_x\xrightarrow{j_*}\coh_{G_n}\bA^{n}_x.
\]
This finishes the proof.
\end{proof}
\end{lem}

\begin{lem}\label{uppsi descri}
For $F=E^0\in \cC_0$ in \eqref{E def} or any $F\in \cC_{n-1}$ with $n\geq 2$, we have an isomorphism
\[\Uppsi(F)\cong \uppsi F.\]
\begin{proof}
We only prove the case when $F\in \cC_{n-1}$ with $n\geq2$. By Lemma \ref{psi coh des}, the functor $\Uppsi^{\fmod}:\cD_{n-1}\to \cD_n$ is the composition
\[
\cD_{n-1}\xrightarrow{(-)^{L_n}}\Dmod^{L_n}(S^{n-1},f_{n-1})\xrightarrow{q_*}\cD_n,
\]
where $q:S^n\to S^n/\l x_n\r\simto S^{n-1}$ is the quotient map. Thus it suffices to show that the object $q_*(F^{L_n})$ is isomorphic to $\uppsi F$ in $\cD_n$. We define the functor $(\widetilde{-}):\proj^{L_{n-1}}S^{n-1}\to \proj^{L_{n}}S^{n}$  by
 $\widetilde{P}:=P^{L_n}\otimes_{S^{n-1}}S^n$. Then we have a short exact sequence 
\[
0\to \widetilde{P}\xrightarrow{x_n}\widetilde{P}\xrightarrow{q_P}q_*(P^{L_n})\to0
\]
in $\fmod^{L_n}S^n$, where $q_P:=q\otimes P^{L_n}$. This  induces  a short exact sequence
\begin{equation}\label{mf resol psi}
0\longrightarrow E^1\xrightarrow{\,\,\,\alpha\,\,\,} E^0\xrightarrow{\,\,\,\beta\,\,\,} q_*(F^{L_n})\longrightarrow0
\end{equation}
in the abelian category $\fmod^{L_{n}}(S^n,f_n)$ of graded factorizations, where $E^i$ and morphisms $\alpha,\beta$ are defined as follows: $E^0$ is defined by
\[
E^0:=\left( E^0_1
\xrightarrow{\left(\begin{array}{cc}\widetilde{\varphi}_1&x_{n}\\-x_{n-1}x_n^{a_n-1}&\widetilde{\varphi}_0\end{array}\right)}
E^0_1
\xrightarrow{\left(\begin{array}{cc}\widetilde{\varphi}_0&-x_{n}\\x_{n-1}x_n^{a_n-1}&\widetilde{\varphi}_1\end{array}\right)}
E^0_1(\vec{f}_{n})
\right),
\]
where  $E^0_1:=\widetilde{F}_1\oplus (\widetilde{F}_0(-\vec{x}_n))$ and  $E^0_0:=\widetilde{F}_0\oplus (\widetilde{F}_1\bigl(\vec{f}_{n}-\vec{x}_n)\bigr)$, and $E^1$ is defined by 
\[
E^1:=\left( E^1_1
\xrightarrow{\left(\begin{array}{cc}\widetilde{\varphi}_1&1\\-x_{n-1}x_n^{a_n}&\widetilde{\varphi}_0\end{array}\right)}
E^1_1
\xrightarrow{\left(\begin{array}{cc}\widetilde{\varphi}_0&-1\\x_{n-1}x_n^{a_n}&\widetilde{\varphi}_1\end{array}\right)}
E^1_1(\vec{f}_{n})
\right),
\]
where $E^1_1:=\widetilde{F}_1(-\vec{x}_n)\oplus (\widetilde{F}_0(-\vec{x}_n))$ and $E^1_0:=\widetilde{F}_0(-\vec{x}_n)\oplus (\widetilde{F}_1\bigl(\vec{f}_{n}-\vec{x}_n\bigr))$. The morphism $\alpha=(\alpha_1,\alpha_0):E^1\to E^0$ is given by 
\[
\alpha_1:E^1_1\xrightarrow{\left(\begin{array}{cc}x_{n}&0\\0&1\end{array}\right)}E^0_1\mbox{\hspace{8mm}and\hspace{8mm}}\alpha_0:E^1_0\xrightarrow{\left(\begin{array}{cc}x_{n}&0\\0&1\end{array}\right)}E^0_0,
\]
and the morphism $\beta=(\beta_1,\beta_0):E^0\to \Upgamma(F)$ is given by
\[
\beta_1:E^0_1\xrightarrow{\left(\begin{array}{cc}q_{F_1}&0\end{array}\right)}q_*(F_1^{L_n})\mbox{\hspace{8mm}and\hspace{8mm}}
\beta_0:E^0_0\xrightarrow{\left(\begin{array}{cc}q_{F_0}&0\end{array}\right)}q_*(F_0^{L_n}).
\]
By \cite[Lemma 2.7(a)]{ls} the short exact sequence \eqref{mf resol psi} induces a  triangle 
\[
 E^1\xrightarrow{\,\,\,\alpha\,\,\,} E^0\xrightarrow{\,\,\,\beta\,\,\,} q_*(F^{L_n})\longrightarrow E^1[1]
\]
in $\cD_n$.   If we set 
\[h_0:E^1_0\xrightarrow{\left(\begin{array}{cc}0&0\\1&0\end{array}\right)} E^1_1\mbox{\hspace{8mm}and\hspace{8mm}}h_1:E^1_1(\vec{f}_n)\xrightarrow{\left(\begin{array}{cc}0&0\\-1&0\end{array}\right)}E^1_0,\] 
then we have $\id_{E^1_0}=\varphi_1^{E^1}h_0+h_1\varphi_0^{E^1}$ and $\id_{E^1_1}=h_0\varphi_1^{E^1}+\varphi_0^{E^1}(-\vec{f}_n)h_1(-\vec{f}_n)$, which means that  the identity map $\id: E^1\to E^1$  is homotopy equivalent to the zero map. Hence  $E^1$ is isomorphic to the zero object in $\cD_n$, and thus we have 
\[\Uppsi(F)\simto q_*(F^{L_n})\simto E^0=\uppsi F.\]
This finishes the proof.
\end{proof}
\end{lem}

By Proposition \ref{coh mod} we have  equivalences
\begin{eqnarray}
\Dcoh_{\bG_m\times G_{n-1}}(Q,1\times\chi_{f_{n-1}},W_Q)&\cong& \Dmod^{\bZ\oplus L_{n-1}}(R,W_Q)\label{coh Q}\\
\Dcoh_{\bG_m\times G_{n-1}}(Q_+,1\times\chi_{f_{n-1}},W_Q)&\cong& \Dmod^{\bZ\oplus L_{n-1}}(R_+,W_Q)\\
\Dcoh_{\bG_m\times G_{n-1}}(Q_-,1\times\chi_{f_{n-1}},W_Q)&\cong& \Dmod^{\bZ\oplus L_{n-1}}(R_-,W_Q).\label{Q- coh mod}
 \end{eqnarray}
 For an interval $I\subset \bZ$, denote by $\cF_I$ the set of  finite direct sums $\oplus_i F^i$ of $(\bZ\oplus L_{n-1})$-graded rank one free $R$-modules $F^i$ with $F^i\cong R(a,l)$ for some $a\in I$ and $l\in L_{n-1}$. Then the $I$-grade-window $\cW_{\lambda,I}$ with respect to $\lambda$ in \eqref{I window} corresponds to  the following subcategory
 \begin{equation}\label{mf window}
 \left\{F=\left(F_1\xrightarrow{\varphi_1} F_0\xrightarrow{\varphi_0}  F_1\left(0,\vec{f}_{n-1}\right)\right)\middle| F_i\in \cF_I\right\}\subset \HMF^{\bZ\oplus L_{n-1}}_{R}(W_Q),
 \end{equation}
which is denoted by the same notation $\cW_{\lambda,I}$, of  $\HMF^{\bZ\oplus L_{n-1}}_{R}(W_Q)$ via the natural equivalences.  Thus, since $I_{0}^+=[-a_n+1,0]$ and $I_{a_n-1}^{-}=[a_n-1]:=[a_n-1,a_n-1]$, we have
\begin{eqnarray}
\cW_{0}^{+}&=&\cW_{\lambda,[-a_n+1,0]}=\left\{F\in\HMF^{\bZ\oplus L_{n-1}}_{R}(W_Q)\middle|  F_i\in \cF_{[-a_n+1,0]}\right\}\\
\cW_{a_n-1}^{-}&=&\cW_{\lambda,[-a_n+1]}=\left\{F\in\HMF^{\bZ\oplus L_{n-1}}_{R}(W_Q)\middle|  F_i\in \cF_{[-a_n+1]}\right\}.\label{window minus}
\end{eqnarray}

\begin{lem}\label{upphi descri} For $F=E^0\in \cC_0$ in \eqref{E def} or any $F\in\cC_{n-2}$ with $n\geq 3$, we have an isomorphism 
\[\Upphi(F)\cong \upphi F\left((-a_n+1)\vec{x}_n\right).\]
\begin{proof}
The functor $\Upphi:\cC_{n-2}\to\cC_n$ is the composition
\begin{eqnarray*}
\cC_{n-2}&\xrightarrow{(-)^{L_{n-1}}\,}&\HMF^{L_{n-1}}(S^{n-2},f_{n-2})\\
&\xrightarrow{(\Phi_{-}^{\fmod})^{-1}}&\Dmod^{\bZ\oplus L_{n-1}}(R_-,W_Q)\\
&\xrightarrow{\,\,\,\quad\sim\quad\,\,\,}&\HMF^{\bZ\oplus L_{n-1}}_{R_-}(W_Q)\\
&\xrightarrow{\Phi_{(-a_n+1)}^{\proj}}&\HMF^{\bZ\oplus L_{n-1}}_{R_+}(W_Q)\\
&\xrightarrow{\,\,\,\,\,\Phi_+^{\proj}\,\,\,\,\,\,}&\cC_n,
\end{eqnarray*}
where the first functor is induced by the functor $(-)^{L_{n-1}}:\fmod^{L_{n-2}}S^{n-2}\to\fmod^{L_{n-1}}S^{n-2}$ (recall notation \eqref{extension}) associated to the inclusion $L_{n-2}\hookto L_{n-1}$ corresponding to the surjection $\pi_{n-1}:G_{n-1}\to G_{n-2}$, and $\Phi_{-}^{\fmod}$, $\Phi_{(-a_n+1)}^{\proj}$ and $\Phi_+^{\proj}$ are the functors corresponding to $\Phi_-$, $\Phi_{(-a_n+1)}$ and $\Phi_+$ in Section \ref{chain sod section} respectively via natural equivalences.

First, we describe the composition 
\begin{equation}\label{first comp}
(\Phi_-^{\fmod})^{-1}\circ (-)^{L_{n-1}}:\cC_{n-2}\to \Dmod^{\bZ\oplus L_{n-1}}(R_-,W_Q).\end{equation}
The quasi-inverse 
\[
(\Phi_-)^{-1}:\Dcoh_{G_{n-1}}(\bA^{n-2}_x,\chi_{f_{n-1}},f_{n-2})\simto \Dcoh_{\bG_m\times G_{n-1}}(Q_-,1\times\chi_{f_{n-1}},W_Q)
\]
of $\Phi_-$ in Lemma \ref{main lemma} is given by 
\[
(\Phi_-)^{-1}\left(F_1\xrightarrow{\varphi_1} F_0\xrightarrow{\varphi_0}  F_1( \chi_{f_{n-1}})\right)=\left(\Updelta(F_1)\xrightarrow{\Updelta(\varphi_1)} \Updelta(F_0)\xrightarrow{\Updelta(\varphi_0)}  \Updelta(F_1)(1\times \chi_{f_{n-1}})\right),
\]
where $\Updelta:\coh_{G_{n-1}}\bA^{n-2}_x\to \coh_{\bG_m\times G_{n-1}}Q_-$ is the following composition:
\begin{eqnarray*}
\coh_{G_{n-1}}\bA^{n-2}_x&\xrightarrow{{(p_1)}^*_{p_2}}& \coh_{\bG_m\times G_{n-1}}\bA^{n-2}_x\times(\bA^1_{x_n}\setminus\{0\})\\
&\xrightarrow{\,\,\,p^*_{1,2}\,\,\,}& \coh_{\bG_m\times G_{n-1}}\bA^{n-2}_x\times (\bA^1_{x_n}\setminus\{0\})\times\bA^1_u\\
&\xrightarrow{\,\,\,\,\,k_*\,\,\,\,\,}& \coh_{\bG_m\times G_{n-1}}Q_-,
\end{eqnarray*}
where $p_1:\bA^{n-2}_x\times(\bA^1_{x_n}\setminus\{0\})\to \bA^{n-2}_x$, $p_2: \bG_m\times G_{n-1}\to G_{n-1}$ and $p_{1,2}:\bA^{n-2}_x\times (\bA^1_{x_n}\setminus\{0\})\times\bA^1_u\to\bA^{n-2}_x\times (\bA^1_{x_n}\setminus\{0\})$ are natural projections and $k:\bA^{n-2}_x\times (\bA^1_{x_n}\setminus\{0\})\times\bA^1_u\hookto Q_-$ is the closed immersion. Thus the composition  \eqref{first comp} sends an object $F=(F_1\xrightarrow{\varphi_1}F_0\xrightarrow{\varphi_1}F_1(\vec{f}_{n-2}))\in \cC_{n-2}$ to the object
\begin{equation}\label{gamma F}
\Upgamma(F)=\left(\Upgamma(F_1)\xrightarrow{\Upgamma(\varphi_1)}\Upgamma(F_0)\xrightarrow{\Upgamma(\varphi_0)}\Upgamma(F_1)(\vec{W}_Q)\right)
\end{equation}
where $\Upgamma(-):=q_*\circ e^* \circ (-)^{{\bZ\oplus L_{n-1}}}$ is the composition 
\begin{eqnarray*}
\proj^{L_{n-2}}S^{n-2}&\xrightarrow{ (-)^{{\bZ\oplus L_{n-1}}}}&\proj^{\bZ\oplus L_{n-1}}S^{n-2}\\
&\xrightarrow{\,\,\,\,\quad e^*\quad\,\,\,\,}&\proj^{\bZ\oplus L_{n-1}}S^{n-2}[x_n^{\pm1},u]\\
&\xrightarrow{\,\,\,\,\quad q_*\quad\,\,\,\,}&\fmod^{\bZ\oplus L_{n-1}}R_-
\end{eqnarray*}
and $e:S^{n-2}\hookto S^{n-2}[x_n^{\pm1},u]$ is the natural inclusion and $q:R_-\to (S^{n-1}/\l x_{n-1}\r)[x_n^{\pm1},u]\simto S^{n-2}[x_n^{\pm1},u]$ is the quotient morphism. 

Next we replace the object $\Upgamma(F)$ in \eqref{gamma F} with a matrix factorization.
For  a graded projective module $P\in \proj^{L_{n-2}}S^{n-2}$, we write 
\[(\,\overline{-}\,):\proj^{L_{n-2}}S^{n-2}\to \proj^{\bZ\oplus L_{n-1}}R_-\]
 for the functor defined by 
\[\overline{P}:=(P^{\bZ\oplus L_{n-2}})\otimes _{S^{n-2}}R_-\in \proj^{\bZ\oplus L_{n-1}}R_-.\] 
 Then we have a short exact sequence
\[
0\longrightarrow \overline{P}(-d)\xrightarrow{x_{n-1}}\overline{P}\xrightarrow{q_P} \Upgamma(P)\longrightarrow 0
\]
in $\fmod^{\bZ\oplus L_{n-1}}R_-$, where $d:=\deg(x_{n-1})=(0,\vec{x}_{n-1})\in \bZ\oplus L_{n-1}$ and $q_P:=q\otimes P:\overline{P}\to\Upgamma(P)$ is the surjection  induced by the  quotient map $q:R_-\to S^{n-2}[x_n^{\pm1},u]$. This  induces  a short exact sequence
\[
0\longrightarrow E^1\xrightarrow{\,\,\,\alpha\,\,\,} E^0\xrightarrow{\,\,\,\beta\,\,\,} \Upgamma(F)\longrightarrow0
\]
in the abelian category $\fmod^{\bZ\oplus L_{n-1}}(R_-,W_Q)$ of graded factorizations, where $E^i$ and morphisms $\alpha,\beta$ are defined as follows: $E^0$ is defined by
\[
E^0:=\left( E^0_1
\xrightarrow{\left(\begin{array}{cc}\overline{\varphi}_1&x_{n-1}\\-g&\overline{\varphi}_0\end{array}\right)}
E^0_0
\xrightarrow{\left(\begin{array}{cc}\overline{\varphi}_0&-x_{n-1}\\g&\overline{\varphi}_1\end{array}\right)}
E^0_1(0,\vec{f}_{n-1})
\right),
\]
where $g:=x_{n-2}x_{n-1}^{a_{n-1}-1}+x_n^{a_n}u$, $E^0_1:=\overline{F}_1\oplus (\overline{F}_0(-d))$ and  $E^0_0:=\overline{F}_0\oplus (\overline{F}_1\bigl((0,\vec{f}_{n-1})-d)\bigr)$, and $E^1$ is defined by 
\[
E^1:=\left( E^1_1
\xrightarrow{\left(\begin{array}{cc}\overline{\varphi}_1&1\\-gx_{n-1}&\overline{\varphi}_0\end{array}\right)}
E^1_0
\xrightarrow{\left(\begin{array}{cc}\overline{\varphi}_0&-1\\gx_{n-1}&\overline{\varphi}_1\end{array}\right)}
E^1_1(0,\vec{f}_{n-1})
\right),
\]
where $E^1_1:=\overline{F}_1(-d)\oplus (\overline{F}_0(-d))$ and $E^1_0:=\overline{F}_0(-d)\oplus (\overline{F}_1\bigl((0,\vec{f}_{n-1})-d\bigr))$. The morphism $\alpha=(\alpha_1,\alpha_0):E^1\to E^0$ is given by 
\[
\alpha_1:E^1_1\xrightarrow{\left(\begin{array}{cc}x_{n-1}&0\\0&1\end{array}\right)}E^0_1\mbox{\hspace{8mm}and\hspace{8mm}}\alpha_0:E^1_0\xrightarrow{\left(\begin{array}{cc}x_{n-1}&0\\0&1\end{array}\right)}E^0_0,
\]
and the morphism $\beta=(\beta_1,\beta_0):E^0\to \Upgamma(F)$ is given by
\[
\beta_1:E^0_1\xrightarrow{\left(\begin{array}{cc}q_{F_1}&0\end{array}\right)}\Upgamma(F_1)\mbox{\hspace{8mm}and\hspace{8mm}}
\beta_0:E^0_0\xrightarrow{\left(\begin{array}{cc}q_{F_0}&0\end{array}\right)}\Upgamma(F_0).
\]
By the same argument as in the proof of Lemma \ref{uppsi descri}, we see that there is an isomorphism
\[
\Upgamma(F)\simto E^0\in\HMF^{\bZ\oplus L_{n-1}}_{R_-}(W_Q).
\]

Finally we compute the image of $E^0\in\HMF^{\bZ\oplus L_{n-1}}_{R_-}(W_Q)$ by the composition
\[
\HMF^{\bZ\oplus L_{n-1}}_{R_-}(W_Q)\xrightarrow{\Phi_{(-a_n+1)}^{\proj}}\HMF^{\bZ\oplus L_{n-1}}_{R_+}(W_Q)\xrightarrow{\Phi_+^{\proj}}\cC_n.
\]
The first functor $\Phi_{(-a_n+1)}^{\proj}$ is the composition
\[
\HMF^{\bZ\oplus L_{n-1}}_{R_-}(W_Q)\xrightarrow{(i_{-}^*)^{-1}}\cW_{a_n-1}^-\subset \cW_0^+\xrightarrow{i_+^*}\HMF^{\bZ\oplus L_{n-1}}_{R_+}(W_Q),
\]
where $i_-:R\hookto R_-$ and $i_+:R\hookto R_+$ are natural inclusions. Write 
\[
(\reallywidecheck{-}):\proj^{L_{n-2}}S^{n-2}\to \proj^{\bZ\oplus L_{n-1}}R
\]
for the functor given by $\reallywidecheck{P}:=P^{\bZ\oplus L_{n-1}}\otimes_{S^{n-2}}R$, where $R$ is the $S^{n-2}$-algebra via the natural inclusion $S^{n-2}\hookto R$.
We define an object $\upgamma F\in \HMF^{\bZ\oplus L_{n-1}}_R(W_Q)$ by

\[
\upgamma F:=\left( \upgamma F_1
\xrightarrow{\left(\begin{array}{cc}\reallywidecheck{\varphi}_1&x_{n-1}\\-g&\reallywidecheck{\varphi}_0\end{array}\right)}
\upgamma F_0
\xrightarrow{\left(\begin{array}{cc}\reallywidecheck{\varphi}_0&-x_{n-1}\\g&\reallywidecheck{\varphi}_1\end{array}\right)}
\upgamma F_1(0,\vec{f}_{n-1})
\right),
\]
where $g=x_{n-2}x_{n-1}^{a_{n-1}-1}+x_n^{a_n}u$, and  $\upgamma F_1:=\reallywidecheck{F}_1\oplus (\reallywidecheck{F}_0(-d))$ and  $\upgamma F_0:=\reallywidecheck{F}_0\oplus (\reallywidecheck{F}_1\bigl((0,\vec{f}_{n-1})-d)\bigr)$.
Then $\upgamma F\in \cW_{\lambda,[0]}$, where $[0]:=[0,0]$ and $\cW_{\lambda,[0]}$ is the  window subcategory defined in \eqref{mf window}. Since $\deg(x_n)=(1,0)$,   by the equality \eqref{window minus}
\[
\upgamma F((-a_n+1)\vec{x}_n)\in \cW_{a_n-1}^-,
\]
and $i_-^*(\upgamma F((-a_n+1)\vec{x}_n))\cong E^0((-a_n+1)\vec{x}_n))$.
Since $x_n\in R_-$ is a unit, we have an isomorphism $x_n^{a_n-1}:E^0((-a_n+1)\vec{x}_n))\simto E^0$.  Thus the image of $E^0$ by the equivalence $\HMF^{\bZ\oplus L_{n-1}}_{R_-}(W_Q)\xrightarrow{(i_{-}^*)^{-1}}\cW_{a_n-1}^-$ is isomorphic to the object $\upgamma F((-a_n+1)\vec{x}_n)$, and so we have 
\[
\Upphi(F)\cong (\Phi_+^{\proj}\circ i_+^*)(\upgamma F((-a_n+1)\vec{x}_n)).
\]
As we saw in the proof of Lemma \ref{psi coh des}, the functor 
\[
\Phi_+:\Dcoh_{\bG_m\times G_{n-1}}(Q_{+},1\times \chi_{f_{n-1}}, W_Q)\simto \Dcoh_{G_n}(\bA^n_x,\chi_{f_n},f_n)
\]
is given by $e_{\psi\times \pi}^*:\coh_{\bG_m\times G_{n-1}}Q_{+}\to \coh_{G_n}\bA^n_x$, where we use the same notation as in the proof of Lemma \ref{psi coh des}. Hence the corresponding functor 
\[
\Phi_+^{\proj}:\HMF^{\bZ\oplus L_{n-1}}_{R_+}(W_Q)\simto \cC_n
\]
is given by $r^*_{\sigma}:\proj^{\bZ\oplus L_{n-1}}R_+\to \proj^{L_n}S^n$, where $r:R_+\to S^{n}$ is the ring homomorphism given by substituting $1$ to $u$, and $\sigma:\bZ\oplus L_{n-1}\to L_n$ is the group homomorphism  defined by $\sigma(a,\sum_{i=1}^{n-1} b_i\vec{x}_i)=\sum_{i=1}^{n-1} b_i\vec{x}_i+a\vec{x}_n$. It is obvious that we have
\[
(r^*_{\sigma}\circ i_+^*)(\upgamma F((-a_n+1)\vec{x}_n))\cong \upphi F((-a_n+1)\vec{x}_n),
\]
and therefore we have an isomorphism $\Upphi(F)\cong \upphi F((-a_n+1)\vec{x}_n)$.
\end{proof}
\end{lem}

The following remark  will be necessary in  Section \ref{quiver computation}.

\begin{rem}\label{psi phi mor}
The assignment $E\mapsto \uppsi E$ for $E\in \cC_{n-1}$ defines  an exact  functor 
\[
\uppsi:\cC_{n-1}\to \cC_n,
\]
where, for a morphism $\alpha=(\alpha_1,\alpha_0):E\to F$  in $\cC_{n-1}$, the morphism
\[
\uppsi\alpha=\bigl((\uppsi \alpha)_1,(\uppsi\alpha)_0\bigr):\uppsi E\to \uppsi F
\]
is defined by
\begin{eqnarray*}
(\uppsi\alpha)_1&:&(\uppsi E)_1=\widetilde{E}_1\oplus \widetilde{E}_0(-\vec{y})
\xrightarrow{\left(\begin{array}{cc}\widetilde{\alpha}_1&0\\0&\widetilde{\alpha}_0\end{array}\right)}
\widetilde{F}_1\oplus \widetilde{F}_0(-\vec{y})=(\uppsi F)_1\\
(\uppsi \alpha)_0&:&(\uppsi E)_0=\widetilde{E}_0\oplus \widetilde{E}_1(\vec{f}-\vec{y})
\xrightarrow{\left(\begin{array}{cc}\widetilde{\alpha}_0&0\\0&\widetilde{\alpha}_1\end{array}\right)}
\widetilde{F}_0\oplus \widetilde{F}_1(\vec{f}-\vec{y})=(\uppsi F)_0.
\end{eqnarray*}
Note that  the following diagram commutes:
\[\xymatrix{
 \uppsi E\ar[rr]^{\sim}\ar[d]_{\uppsi \alpha}&&q_*(E^{L_n})\ar[d]^{q_*(\alpha^{L_n})}\ar@{=}[r]&\Uppsi^{\fmod}(E)\ar[d]^{\Uppsi^{\fmod}(\alpha)}\\
 \uppsi F\ar[rr]^{\sim}&&q_*(F^{L_n})\ar@{=}[r]&\Uppsi^{\fmod}(F),
 }\]
where the horizontal isomorphisms are given by the morphism $\beta$ in \eqref{mf resol psi}. This means that the functor $\Uppsi:\cC_{n-1}\to \cC_n$ is isomorphic to the functor  $\uppsi:\cC_{n-1}\to \cC_n$.

Similarly, the assignment $E\mapsto \upphi E$ defines a functor  
\[
\upphi:\cC_{n-2}\to \cC_{n}
\]
such that the composition $\Bigl(\bigl(-\bigr)\bigl((-a_n+1)\vec{x}_n\bigr)\Bigr)\circ\upphi$ is isomorphic to $\Upphi:\cC_{n-2}\to \cC_n$.
For a morphism $\alpha=(\alpha_1,\alpha_0):E\to F$ in $\cC_{n-2}$, the morphism 
\[
\upphi\alpha:\upphi E\to \upphi F
\]
is defined by 
\begin{eqnarray*}
(\upphi\alpha)_1&:&(\upphi E)_1=\widehat{E}_1\oplus \widehat{E}_0(-\vec{y})
\xrightarrow{\left(\begin{array}{cc}\widehat{\alpha}_1&0\\0&\widehat{\alpha}_0\end{array}\right)}
\widehat{F}_1\oplus \widehat{F}_0(-\vec{y})=(\upphi F)_1\\
(\upphi \alpha)_0&:&(\upphi E)_0=\widehat{E}_0\oplus \widehat{E}_1(\vec{f}-\vec{y})
\xrightarrow{\left(\begin{array}{cc}\widehat{\alpha}_0&0\\0&\widehat{\alpha}_1\end{array}\right)}
\widehat{F}_0\oplus \widehat{F}_1(\vec{f}-\vec{y})=(\upphi F)_0.
\end{eqnarray*}
\end{rem}

\subsection{Proof of Theorem \ref{strong thm}}
In this subsection, we prove Theorem \ref{strong thm} by induction on $n$. If $a_n=1$, by Theorem \ref{sod for chain} the category $\cC^n$ is orthogonally decomposed into  some copies of $\cC^{n-2}$, and so nothing to prove. Thus we may assume $a_n\geq 2$. First we compute the Milnor number $\mu_n:=\mu(\widetilde{f}_n)$ of $\widetilde{f}_n$. Set 
\[
b_i:=\prod_{j=0}^{n-i}a_{n-j}=a_na_{n-1}\cdots a_i.
\]
\begin{lem}\label{milnor number}
Assume $a_n\geq 2$. We have 
\[
\mu_n=b_1-b_2+\cdots+(-1)^{n-1}b_n+(-1)^{n}.
\]
\begin{proof}
Since $\widetilde{f}_n$ is a quasi-homogeneous polynomial with an isolated singularity only at the origin, the Milnor number $\mu_n$ is equal to $\prod_{i=1}^n(\frac{1}{q_i}-1)$, where $q_i$ is the rational weight of $x_i$ such that the degree of $\widetilde{f}_n$ with respect to $\{q_i\}$ is equal to $1$. Set $c_{n+1}:=1$, and for each $1\leq i\leq n$ set
\[
c_i:=b_i-b_{i+1}+\cdots+(-1)^{n-i}b_n+(-1)^{n-i+1}.
\]
Since $q_i=c_{i+1}/b_i$, we have
\[
\frac{1}{q_i}-1=\frac{c_i}{c_{i+1}}.
\]
Thus 
\[
\mu_n=\prod_{i=1}^n\left(\frac{1}{q_i}-1\right)=\prod_{i=1}^n\frac{c_i}{c_{i+1}}=c_1=b_1-b_2+\cdots+(-1)^{n-1}b_n+(-1)^{n}.
\]
\end{proof}
\end{lem}

\begin{prop}\label{fec}
The sequence $\cE^n$ is a full exceptional collection in $\cC_n$. If $a_n\geq 2$, the length $\nu_n$ of $\cE^n$ is equal to  the Milnor number $\mu_n$ of $\widetilde{f}_n$, where we set $\mu_0:=1$.
\begin{proof}
Since the equivalence $\cC_0=\HMF^{\bZ}_k(0)\cong \Db(\fmod k)$ by Corollary \ref{derived category}  maps the object $E^0\in \cC_0$  to the free module  $k\in \Db(\fmod k)$, $\cE^0$ is a full strong exceptional collection of length $1$. Since we have $\cC_1=\l \Image\Uppsi_0,\hdots,\Image\Uppsi_{a_1-2}\r$, $\cE_1$ is a full exceptional collection of length $\nu_1=a_1-1$.   
By induction on $n$, the former statement follows from Theorem \ref{sod for chain}, Lemma \ref{uppsi descri} and Lemma \ref{upphi descri}. 

By Theorem \ref{sod for chain} we have $\nu_n=(a_n-1)\nu_{n-1}+a_{n-1}\nu_{n-2}$. Since $\nu_0=1$ and $\nu_1=a_1-1$, this recursion shows the equality
\[
\nu_n=b_1-b_2+\cdots+(-1)^{n-1}b_n+(-1)^{n}.
\]
Thus the latter assertion follows from Lemma \ref{milnor number}.
\end{proof}
\end{prop}

Atsushi Takahashi pointed out that the semi-orthogonal decomposition in Theorem \ref{sod for chain}  may be related to the result of Gabrielov \cite{Gab} via mirror symmetry:

\begin{rem}
From the semi-orthogonal decompositions in Theorem \ref{sod for chain}, we have a recursion 
$\nu_n=(a_n-1)\nu_{n-1}+a_{n-1}\nu_{n-2}$. The equivalent recursion $\mu_n=(a_n-1)\mu_{n-1}+a_{n-1}\mu_{n-2}$ can be deduced from the geometry of the singularity $(\widetilde{f}_n,0)$ \cite[Theorem 1]{Gab}.
\end{rem}

 The following lemma implies that $\cE^1$ is a full strong exceptional collection.

\begin{lem}\label{psi part}
Let $E$ and $F$ be objects in $\cC_{n-1}$. For $0\leq i< j\leq a_n-2$, we have the following
\[
\Hom(\uppsi_iE,\uppsi_jF[l])\cong
\begin{cases}
\Hom(E,F[l])& j=i+1\\
0&\mbox{otherwise.}
\end{cases}
\]
In particular, if $\cE^{n-1}$ is strongly exceptional,  the sequence 
\[
\left(\uppsi_0\cE,\hdots, \uppsi_{b-2}\cE\right)
\]
 is also strongly exceptional.
\begin{proof} 
Applying the Serre functor and using the isomorphism $(\uppsi_{j} F)(\vec{x}_n)\cong (\uppsi_{j-1}F)[1]$, we obtain  isomorphisms of vector spaces
\begin{eqnarray*}
\Hom(\uppsi_iE, \uppsi_jF[l])&\cong &\Hom(\uppsi_j F[l],(\uppsi_iE)(-\vec{\bf x}^n)[n])\\
&\cong &\Hom(\uppsi_{j-1}F[l+1], \uppsi_i(E(-\vec{\bf x}^{n-1}))[n])
\end{eqnarray*}

First assume that $j>i+1$. The object $\uppsi_{j-1} F[l+1]$ lies in $\Image\Psi_{-j+1}$, and the object $\uppsi_i\left(E(-\vec{\bf x}^{n-1})\right)[n]$ lies  in $\Image \Psi_{-i}$. Since $-j+1<-i$, by the semi-orthogonal decomposition in Theorem \ref{sod for chain}, we have 
\[
\Hom(\uppsi_iE,\uppsi_jF[l])\cong\Hom(\uppsi_{j-1} F[l+1],\uppsi_i\left(E(-\vec{\bf x}^{n-1})\right)[n])=0.
\]
Next assume that $j=i+1$. Then we have isomorphisms
\begin{eqnarray*}
&&\Hom(\uppsi_iE,\uppsi_jF[l])\\
&\cong& \Hom(\uppsi_{i}  F[l+1],\uppsi_i\left(E(-\vec{\bf x}^{n-1})\right)[n])\\
&\cong& \Hom(F[l+1],E(-\vec{\bf x}^{n-1})[n])\\
&\cong& \Hom(E(-\vec{\bf x}^{n-1})[n],F(-\vec{\bf x}^{n-1})[l+n])\\
&\cong& \Hom(E,F[l]),
\end{eqnarray*}
where the second isomorphism follows from the fully faithfulness of the functor $\Psi_i$, and the third one follows from the Serre duality in $\cC_{n-1}$.
\end{proof}
\end{lem}

Next we compute the set of morphisms  $\Hom_{\cC_n}(\uppsi_i\cE^{n-1},\upphi_j\cE^{n-2}[l])$ for any $l\in \bZ$. For this, we need the following lemmas.

\begin{lem}\label{psi phi tri}
For each object $F\in\cC_n$, we have the following triangle in $\cC_{n+2}$.
\[
\upphi F(-\vec{z})\xrightarrow{z} \upphi F\to \uppsi^2F\to \upphi F(-\vec{z})[1].
\]
\begin{proof}
We use the notation \eqref{easy notation} and  set $g:=xy^{b-1}+z^c$. By construction, the object 
\[
\uppsi^2F=\left((\uppsi^2F)_1\xrightarrow{\varphi_1^{\uppsi^2F}}(\uppsi^2F)_0\xrightarrow{\varphi_0^{\uppsi^2F}}(\uppsi^2F)_1(\vec{f})\right)
\]
is given by
\begin{eqnarray*}
(\uppsi^2F)_1&:=&\widehat{F_1}\oplus\left(\widehat{F_0}(-\vec{y})\right)\oplus\left(\widehat{F}_0(-\vec{z})\right)\oplus\left(\widehat{F}_1(\vec{f}-\vec{y}-\vec{z})\right)\\
(\uppsi^2F)_0&:=&\widehat{F_0}\oplus\left(\widehat{F}_1(\vec{f}-\vec{y})\right)\oplus\left(\widehat{F}(\vec{f}-\vec{z})\right)\oplus\left(\widehat{F}_0(\vec{f}-\vec{y}-\vec{z})\right)
\end{eqnarray*}
\[
\varphi_1^{\uppsi^2F}:=\left(
\begin{array}{cccc}
\varphi_1&y&z&0\\
-xy^{b-1}&\varphi_0&0&z\\
-yz^{c-1}&0&\varphi_0&-y\\
0&-yz^{c-1}&xy^{b-1}&\varphi_1
\end{array}\right),
\hspace{5mm}
\varphi_0^{\uppsi^2F}:=\left(
\begin{array}{cccc}
\varphi_0&-y&-z&0\\
xy^{b-1}&\varphi_1&0&-z\\
yz^{c-1}&0&\varphi_1&y\\
0&yz^{c-1}&-xy^{b-1}&\varphi_0
\end{array}\right).
\]
By applying  elementary row
and column operations, we see that $\uppsi^2F$ is isomorphic to the following matrix factorization
\[(\uppsi^2F)_1\xrightarrow{
\left(
\begin{array}{cccc}
\varphi_1&y&z&0\\
-g&\varphi_0&0&z\\
0&0&-\varphi_0&y\\
0&0&-g&-\varphi_1
\end{array}\right)
}
(\uppsi^2F)_0\xrightarrow{
\left(
\begin{array}{cccc}
\varphi_0&-y&z&0\\
g&\varphi_1&0&z\\
0&0&-\varphi_1&-y\\
0&0&g&-\varphi_0
\end{array}\right)
}
(\uppsi^2F)_1(\vec{f}),
\]
which is nothing but the mapping cone of the multiplication $z:\upphi F(-\vec{z})\to \upphi F$ by $z$. Indeed, the following matrices 
\[
\alpha_1:=\left(
\begin{array}{cccc}
1&0&0&0\\
0&1&0&0\\
0&0&1&0\\
z^{c-1}&0&0&1
\end{array}\right)
\hspace{3mm}\mbox{ and }\hspace{3mm}
\alpha_0:=\left(
\begin{array}{cccc}
1&0&0&0\\
0&1&0&0\\
0&0&-1&0\\
-z^{c-1}&0&0&-1
\end{array}\right)
\]
defines an isomorphism $\alpha=(\alpha_1,\alpha_0):\uppsi^2F\to \Cone(z:\upphi F(-\vec{z})\to \upphi F)$.
\end{proof}
\end{lem}

\begin{lem}\label{5.D}
For $E\in \cC_{n-1}$ and $0\leq j\leq a_{n-1}-1$, set 
\begin{equation}\label{Ej notation}
\widetilde{E}_j:=(\uppsi E)(j\vec{x}_{n-1}+(-a_n+2)\vec{x}_n)[a_n-j-2]\in\cC_n.
\end{equation}
For an object $F\in \cC_{n-2}$ and any $l\in \bZ$, we have an isomorphism
\[
\Hom(\widetilde{E}_j, \upphi_0F[l])\cong \Hom(\widetilde{E}_j,(\uppsi^2F)((-a_n+1)\vec{x}_n)[a_n-1+l]).
\]

\begin{proof} By Lemma \ref{psi phi tri} and $\upphi_0F=(\upphi F)((-a_n+1)\vec{x}_n)[a_n-1]$, we have following  triangle 
\[
(\upphi _0F)(-\vec{x}_n)\to \upphi_0 F\to (\uppsi^2F)((-a_n+1)\vec{x}_n)[a_n-1]\to (\upphi_0 F)(-\vec{x}_n)[1].
\]
By the long exact sequence obtained by applying $\Hom(\widetilde{E}_j,-)$ to the above triangle, it suffices to show the  vanishing
\begin{equation}\label{vanish}
\Hom(\widetilde{E}_j,\upphi_0F(-\vec{x}_n)[l])=0
\end{equation}
for any $l\in \bZ$. By the Serre duality, we have isomorphisms
\begin{eqnarray*}
\Hom(\widetilde{E}_j,\upphi_0F(-\vec{x}_n)[l])
&\cong& \Hom(\upphi_0F(-\vec{x}_n)[l],\widetilde{E}_j(-\vec{\bf x}^n)[n])\\
&\cong& \Hom(\upphi_0F[l],\widetilde{E}_j(-\vec{\bf x}^{n-1})[n]).
\end{eqnarray*}
Since $\widetilde{E}_j(-\vec{\bf x}^{n-1})=\uppsi (E(-\vec{\bf x}^{n-2}+(j-1)\vec{x}_{n-1}))(-a_n+2)[a_n-j-2]\in \Image\Psi_{-a_n+2}$ and $\upphi_0F\in \Image\Phi_0$,  the semi-orthogonality in Theorem \ref{sod for chain} implies the vanishing \eqref{vanish}. 
\end{proof}
\end{lem}

Now we can compute the set $\Hom(\uppsi_iE,\upphi_j F[l])$.

\begin{lem}\label{last lemma}
Let $E\in \cC_{n-1}$ and $F\in \cC_{n-2}$ be objects. For $0\leq i\leq a_n-2 $ and $0\leq j\leq a_{n-1}-1$, we have the following
\[
\Hom(\uppsi_iE,\upphi_j F[l])\cong
\begin{cases}
\Hom(E,\uppsi_jF[l])& i=a_n-2\\
0 & \mbox{otherwise}.
\end{cases}
\]
\begin{proof}
First we show the vanishing $\Hom(\uppsi_iE,\upphi_j F[l])=0$ when $i<a_n-2$. Applying the Serre functor and using isomorphism $\uppsi_iE(-\vec{x}_n)[1]\cong \uppsi_{i+1}E$, we have 
\begin{eqnarray*}
\Hom(\uppsi_iE,\upphi_j F[l])&\cong& \Hom(\upphi_jF[l], (\uppsi_iE)(-\vec{\bf x}^{n})[n])\\
&\cong& \Hom(\upphi_jF[l], \uppsi_{i+1}(E(-\vec{\bf x}^{n-1}))[n-1])
\end{eqnarray*}
By the semi-orthogonal decomposition in Theorem \ref{sod for chain}, the object $\upphi_jF[l]$ is left orthogonal to the full subcategory $\Image\Psi_{i+1}$, and thus $\Hom(\upphi_jF[l], \uppsi_{i+1}(E(-\vec{\bf x}^{n-1}))[n-1])=0$.

Now it suffices to show an isomorphism 
\[
\Hom(\uppsi_{a_n-2}E,\upphi_j F[l])\cong\Hom(E,\uppsi_jF[l]).
\]
Using notation \eqref{Ej notation} and Lemma \ref{5.D}, we have isomorphisms
\begin{eqnarray}
\Hom(\uppsi_{a_n-2}E,\upphi_j F[l])
&=&\Hom(\uppsi E((-a_n+2)\vec{x}_n)[a_n-2],(\upphi_0F)(-j\vec{x}_{n-1})[l+j])\nonumber\\
&\cong&\Hom(\widetilde{E}_j,\upphi_0F[l])\nonumber\\
&\cong&\Hom(\widetilde{E}_j,(\uppsi^2F)((-a_n+1)\vec{x}_n)[a_n-1+l])\nonumber\\
&\cong&\Hom((\uppsi E)(j\vec{x}_{n-1}+\vec{x}_n),(\uppsi^2F)[l+j+1])\label{psi to pso2}.
\end{eqnarray}
Applying the Serre duality in $\cC_n$ and the fully faithfulness of $\Uppsi$, we  have isomorphisms
\begin{eqnarray}
\eqref{psi to pso2}&\cong&\Hom((\uppsi^2F)[l+j+1],(\uppsi E)(j\vec{x}_{n-1}+\vec{x}_n-\vec{\bf x}^n)[n])\nonumber\\
&\cong&\Hom((\uppsi^2F)[l+j+1],\uppsi\left(E(j\vec{x}_{n-1}-\vec{\bf x}^{n-1})\right)[n])\nonumber\\
&\cong&\Hom(\uppsi F[l+j],E(j\vec{x}_{n-1}-\vec{\bf x}^{n-1})[n-1])\label{psi to E}
\end{eqnarray}
Again by the Serre duality in $\cC_{n-1}$, we have the isomorphisms
\begin{eqnarray*}
\eqref{psi to E}&\cong&\Hom(E(j\vec{x}_{n-1}-\vec{\bf x}^{n-1})[n-1], (\uppsi F)(-\vec{\bf x}^{n-1})[l+j+n-1])\\
&\cong&\Hom(E,\uppsi F(-j\vec{x}_{n-1})[j+l])\\
&\cong&\Hom(E,\uppsi_jF[l]).
\end{eqnarray*}
This completes the proof. 
\end{proof}
\end{lem}

Now for Theorem \ref{strong thm}, it suffices to prove the following.

\begin{cor}
If $\cE^{n-1}$ and $\cE^{n-2}$ are full strong exceptional collections, then $\cE^n$ is also a full strong exceptional collection.
\begin{proof}
Write $\cE^n=\{E_1,\hdots, E_{\mu}\}$ ($\mu:=\mu(f_n)$). By Proposition \ref{fec}, it suffices to show that for any two objects $E_{s}, E_{s'}\in \cE^n$ with $s<s'$ and any non-zero integer $l\neq0$,  the vanishing
\[
\Hom(E_s,E_{s'}[l])=0
\] 
holds. Since the sequence $\cE^{n-1}$ is strongly exceptional, Lemma \ref{psi part} implies that the above vanishing holds if  $E_s\in \uppsi_{i}\cE^{n-1}$ and  $E_{s'}\in\uppsi_{i'}\cE^{n-1}$ for some $0\leq i,i'\leq a_n-2$.
Moreover, the above vanishing also holds when $E_s\in\upphi_{j}\cE^{n-2}$ and $E_{s'}\in \upphi_{j'}\cE^{n-2}$ for some $0\leq j, j',\leq a_{n-1}-1$, since if $j\neq j'$ then any elements $\upphi_{j}\cE^{n-2}$  is orthogonal to any elements in $\upphi_{j'}\cE^{n-2}$, and if $j=j'$  the vanishing follows from the assumption that the sequence $\cE^{n-2}$ is strongly exceptional.
 Therefore it is enough to show  the vanishing
\[
\Hom(\uppsi_iE,\upphi_j F[l])=0
\]
for any $E\in \cE^{n-1}$, $F\in \cE^{n-1}$, $0\leq i\leq a_n-2$, $0\leq j\leq a_{n-1}-1$ and $l\neq0$. By Lemma \ref{last lemma},   this vanishing holds for any $l\in \bZ$ and  $0\leq i<a_n-2$, and if $i=a_n-2$,  we have an isomorphism 
\begin{equation}\label{last step}
\Hom(\uppsi_{a_n-2}E,\upphi_j F[l])\cong \Hom(E,\uppsi_jF[l]).
\end{equation}
Since both of $E$ and $\uppsi_jF$ lie  in $\cE^{n-1}$ and the sequence $\cE^{n-1}$ is strongly exceptional, the right hand side of \eqref{last step} vanishes for any $l\neq0$.
\end{proof}
\end{cor}

\subsection{The quiver with relations associated to the tilting object in $\cC_n$}\label{quiver computation}

Using the results in the previous subsection, we compute the quiver with relation $(Q^n,I^n)$ which represents the endomorphism ring $\End_{\cC_n}(T_n)$ of the associated  tilting object 
\begin{equation}\label{tilting}T_n:=\bigoplus_{E_i\in\cE^n}E_i\in\cC_n.\end{equation}

\subsubsection{Low dimensional cases} 

Since $\cE^0$ consists of a unique exceptional object 
\[E^0=\left(0\to k\to 0\right)\in \cC_0,\] 
the quiver $Q^0$ has a unique vertex, no arrow and no relations. Pictorially, $Q^0$ is 
\[
\begin{tikzpicture}[scale=0.75]
\node (1) at (0.75,0) [DW] {};
\node (1b) at (0.75,0.55) [] {$E^0$};
%\node (1c) at (0.75,1.5) [DW] {};
%\draw [-] (1b) -- (1c);
\end{tikzpicture}
\]

Recall that $\cE^1=\left(\uppsi_0E^0,\hdots,\uppsi_{(a_1-2)}E^0\right)$. By Lemma \ref{psi part}, there is a non-zero  morphism from $\uppsi_iE$ to $\uppsi_jE$ if and only if  $j=i+1$, and in particular the non-zero morphisms are automatically irreducible. Since $\Hom(\uppsi_iE,\uppsi_{i+1}E)=1$,   the corresponding quiver $Q^1$ is 
\vspace{3mm}
\[
\begin{tikzpicture}[scale=0.75]
\node (1a) at (-6,0.55) [] {$\uppsi_0E^0$};
\node (1b) at (-6,0) [DW] {};
\node (2a) at (-3,0.55) [] {$\uppsi_1E^0$};
\node (2b) at (-3,0) [DW] {};
\node (3b) at (0,0) [] {\hspace{-2mm}\!};
\node (4a) at (3,0.55) [] {};
\node (4b) at (3,0) [] {};
\node (5a) at (6,0.55) [] {$\uppsi_{a_1-2}E^0$};
\node (5b) at (6,0) [DW] {};

%\node (1c) at (0.75,1.5) [DW] {};
%\draw [-] (1b) -- (1c);

\draw[->] (1b) --  node[above] {$\upalpha_1$} (2b);
\draw[->] (2b) -- node[above] {$\upalpha_2$}  (3b);
\draw[dotted] (3b) --  (4b);
\draw[->] (4b) --  node[above] {$\upalpha_{a_1-2}$} (5b);
\end{tikzpicture}\vspace{3mm}
\]
and the relations are
\[
I^1=\left\{ \alpha_{i+1}\circ\alpha_i=0\,\,\,\middle| \,\,\,1\leq i\leq a_1-3\right\}.
\]
\subsubsection{Lemmas for higher dimensional cases}\label{lem for}
Recall from Remark \ref{psi phi mor}, the assignments 
$E\mapsto \uppsi E$ and $E\mapsto \upphi E$ define functors $\uppsi:\cC_{n-1}\to \cC_n$ and $\upphi:\cC_{n-2}\to \cC_n$ respectively. Following  notation \eqref{psi phi obj}, for each $0\leq i\leq a_n-2$ and $0\leq j\leq a_{n-1}-1$, we define functors 
\begin{eqnarray*}
\uppsi_i&:&\cC_{n-1}\to\cC_n\\
\upphi_j&:&\cC_{n-2}\to\cC_n
\end{eqnarray*}
by $\uppsi_i:=\Bigl((-)(-i\vec{x}_n)[i]\Bigr)\circ\uppsi$ and $\upphi_j:=\Bigl((-)(-j\vec{x}_{n-1}+(-a_n+1)\vec{x}_n)[a_n+j-1]\Bigr)\circ\upphi$. 

 In order to compute the quiver with relations $(Q^n,I^n)$ in the case when $n\geq2$, we need the following lemmas.

\begin{lem}\label{non-zero lambda}
For each non-zero object $E\in\cC_{n-1}$ and any $0\leq i\leq a_n-3$, there is a non-zero map 
\[
\lambda_i^E:\uppsi_iE\to \uppsi_{i+1}E
\]
such that the following conditions hold:
\begin{itemize}
\item[$(1)$] For any morphism $\alpha=(\alpha_1,\alpha_0):E\to E'$ between non-zero objects  in $\cC_{n-1}$, the following diagram commutes:
\[\xymatrix{
 \uppsi_iE\ar[rr]^{\lambda_i^E}\ar[d]_{\uppsi_i\alpha}&&\uppsi_{i+1}E\ar[d]^{\uppsi_{i+1}\alpha}\\
 \uppsi_iE'\ar[rr]^{\lambda_i^{E'}}&&\uppsi_{i+1}E'.
 }\]
 In particular, the morphisms $\{\lambda_i^E\}$ define a functor morphism $\lambda_i\colon \uppsi_i\to\uppsi_{i+1}$.
\item[$(2)$] The composition $\lambda_i^{E'}\circ \uppsi_i\alpha$ is the zero map if and only if $\alpha$ is the zero map.
\end{itemize}
\begin{proof}
We only prove the result in the case when $i=0$, since $\lambda_i^E$ can be given by $\lambda_i^E:=\lambda_0^E(-i\vec{x}_n)[i]$.  Write $E=\left(E_1\xrightarrow{\varphi_1}E_0\xrightarrow{\varphi_0}E_1(\vec{f})\right)\in \cC_{n-1}$, and consider the following morphisms
\begin{eqnarray*}
(\lambda_0^E)_1&:&(\uppsi E)_1=\widetilde{E}_1\oplus \widetilde{E}_0(-\vec{y})
\xrightarrow{\left(\begin{array}{cc}0&1\\xy^{b-2}&0\end{array}\right)}
\widetilde{E}_0(-\vec{y})\oplus\widetilde{E}_1(\vec{f}-2\vec{y})=(\uppsi_1 E)_1\\
(\lambda_0^E)_0&:&(\uppsi E)_0=\widetilde{E}_0\oplus \widetilde{E}_1(\vec{f}-\vec{y})
\xrightarrow{\left(\begin{array}{cc}0&-1\\-xy^{b-2}&0\end{array}\right)}
\widetilde{E}_1(\vec{f}-\vec{y})\oplus\widetilde{E}_0(\vec{f}-2\vec{y})=(\uppsi_1 E)_0,
\end{eqnarray*}
where we use the similar notation as in \eqref{easy notation}, namely $x:=x_{n-1}$, $y:=x_n$ and $b:=a_n$.
These morphisms define a morphism 
\[\lambda_0^E:=\left((\lambda_0^E)_1, (\lambda_0^E)_0\right):\uppsi E\to \uppsi_{1}E.\] 
Since the object $E$ is a non-zero object in $\cC_{n-1}$, we may assume that the matrices $\varphi_1$ and $\varphi_0$ have no unit entry. On the other hand,  the matrices $(\lambda_0^E)_1$ and $(\lambda_0^E)_0$ have  unit entries. Thus the morphism $\lambda_0^E$ can not be homotopy equivalent to the zero map.

It suffices to  prove that the non-zero morphism $\lambda_0^E$ satisfies the conditions (1) and (2). The condition (1) can be checked by direct computations, and so we omit the details. We check the condition (2). It is obvious that $\alpha=0$ implies  $\lambda_0^{E'}\circ \uppsi\alpha=0$. Assume that $\lambda_0^{E'}\circ \uppsi_i\alpha=0$. We set $\beta:=\lambda_0^{E'}\circ \uppsi\alpha$ and write $E'=\left(E'_1\xrightarrow{\varphi'_1}E'_0\xrightarrow{\varphi'_0}E'_1(\vec{f})\right)$. Then $\beta$ is given by the following morphisms:
\begin{eqnarray*}
\beta_1&:&(\uppsi E)_1=\widetilde{E}_1\oplus \widetilde{E}_0(-\vec{y})
\xrightarrow{\left(\begin{array}{cc}0&\alpha_0\\xy^{b-2}\alpha_1&0\end{array}\right)}
\widetilde{E}'_0(-\vec{y})\oplus\widetilde{E}'_1(\vec{f}-2\vec{y})=(\uppsi_1 E')_1\\
\beta_0&:&(\uppsi E)_0=\widetilde{E}_0\oplus \widetilde{E}_1(\vec{f}-\vec{y})
\xrightarrow{\left(\begin{array}{cc}0&-\alpha_1\\-xy^{b-2}\alpha_0&0\end{array}\right)}
\widetilde{E}'_1(\vec{f}-\vec{y})\oplus\widetilde{E}'_0(\vec{f}-2\vec{y})=(\uppsi_1 E')_0,
\end{eqnarray*}
where we simply write $\alpha_i:\widetilde{E}_i\to\widetilde{E}'_i$ for $\widetilde{\alpha}_i$.
Since $\beta=0$, there are morphisms
\begin{eqnarray*}
g&:&(\uppsi E)_0=\widetilde{E}_0\oplus \widetilde{E}_1(\vec{f}-\vec{y})
\xrightarrow{\left(\begin{array}{cc}g_1&g_2\\g_3&g_4\end{array}\right)}
\widetilde{E}'_0(-\vec{y})\oplus\widetilde{E}'_1(\vec{f}-2\vec{y})=(\uppsi_1 E')_1\\
h&:&(\uppsi E)_1(\vec{f})=\widetilde{E}_1(\vec{f})\oplus \widetilde{E}_0(\vec{f}-\vec{y})
\xrightarrow{\left(\begin{array}{cc}h_1&h_2\\h_3&h_4\end{array}\right)}
\widetilde{E}'_1(\vec{f}-\vec{y})\oplus\widetilde{E}'_0(\vec{f}-2\vec{y})=(\uppsi_1 E')_0
\end{eqnarray*}
such that 
\begin{eqnarray}\label{beta1}
\beta_1&=&g\circ \varphi_1^{\uppsi E}+(\varphi_0^{\uppsi_1E'}(-\vec{f}))\circ(h(-\vec{f}))\nonumber\\
&=&\left(\begin{array}{cc}g_1&g_2\\g_3&g_4\end{array}\right)\left(\begin{array}{cc}\varphi_1&y\\-xy^{b-1}&\varphi_0\end{array}\right)+\left(\begin{array}{cc}-\varphi'_1&-y\\xy^{b-1}&-\varphi'_0\end{array}\right)\left(\begin{array}{cc}h_1&h_2\\h_3&h_4\end{array}\right)
\end{eqnarray}
and 
\begin{eqnarray}\label{beta0}
\beta_0&=&h\circ \varphi_0^{\uppsi E}+\varphi_1^{\uppsi_1E'}\circ g\nonumber\\
&=&\left(\begin{array}{cc}h_1&h_2\\h_3&h_4\end{array}\right)\left(\begin{array}{cc}\varphi_0&-y\\xy^{b-1}&\varphi_1\end{array}\right)+\left(\begin{array}{cc}-\varphi'_0&y\\-xy^{b-1}&-\varphi'_1\end{array}\right)\left(\begin{array}{cc}g_1&g_2\\g_3&g_4\end{array}\right)
\end{eqnarray}
By the equation \eqref{beta1}, we have 
\[
\alpha_0=g_1y+g_2\varphi_0-\varphi'_1h_2+yh_4.
\]
Since $(\alpha_i)|_{y=0}=\alpha_i$, $(\varphi_i)|_{y=0}=\varphi_i$ and $(\varphi'_i)|_{y=0}=\varphi'_i$, we have 
\[
\alpha_0=(g_2|_{y=0})\varphi_0+\varphi'_1(-(h_2)|_{y=0}).
\]
Similarly, by \eqref{beta0}, we have
\[
\alpha_1=(-(h_2)|_{y=0})\varphi_1+\varphi'_0(g_2|_{y=0}).
\]
Hence the morphism $\alpha$ is homotopy equivalent to the zero map.
\end{proof}
\end{lem}

\begin{lem}\label{non-zero mu}
For each non-zero object $F\in \cC_{n-2}$ and $0\leq j\leq a_{n-1}-2$, there is a non-zero morphism 
\[
\sigma_j^F:\uppsi_{a_n-2}\uppsi_j F\to \upphi_j F
\]
such that the following conditions hold:
\begin{itemize}
\item[$(1)$] For any morphism $\alpha=(\alpha_1,\alpha_0):F\to F'$ between non-zero objects  in $\cC_{n-2}$, the following diagram commutes:
\[\xymatrix{
 \uppsi_{a_n-2}\uppsi_j F\ar[rr]^{\sigma_j^F}\ar[d]_{\uppsi_{a_n-2}\uppsi_j\alpha}&&\upphi_{j}F\ar[d]^{\upphi_{j}\alpha}\\
\uppsi_{a_n-2}\uppsi_j F'\ar[rr]^{\sigma_j^{F'}}&&\upphi_{j}F'.
 }\]
 In particular, the morphisms $\{\sigma_j^F\}$ define a functor morphism $\sigma_j\colon \uppsi_{a_n-2}\uppsi_j \to \upphi_j$.
\item[$(2)$] The composition $\sigma_j^{F'}\circ  \uppsi_{a_n-2}\uppsi_j \alpha$ is the zero map if and only if $\alpha$ is the zero map.
\end{itemize}
\begin{proof}
We only give the non-zero morphism $\sigma_j^F$, since the other part can be checked by similar  computations and arguments as in the proof of Lemma \ref{non-zero lambda}.
We use the similar notation as in \eqref{easy notation}.
Then the morphism $\sigma_j^F:\uppsi_{c-2}\uppsi_j F\to \upphi_j F$ is defined by
\[
\sigma_j^F:=\sigma((-(c-2)\vec{z}-j\vec{y}))[c+j-2],
\]
where \[\sigma=(\sigma_1,\sigma_0):\uppsi^2F\to \upphi F(-\vec{z})[1]\]
is given by 
\[
\sigma_1:(\uppsi^2F)_1
\xrightarrow{\left(\begin{array}{cccc}0&0&1&0\\z^{c-1}&0&0&1\end{array}\right)}
(\upphi F(-\vec{z})[1])_1
\]
\[
\sigma_0:(\uppsi^2F)_0
\xrightarrow{\left(\begin{array}{cccc}0&0&1&0\\z^{c-1}&0&0&1\end{array}\right)}
(\upphi F(-\vec{z})[1])_0,
\]
where 
\begin{eqnarray*}
(\uppsi^2F)_1&=&\widehat{F}_1\oplus \widehat{F}_0(-\vec{y})\oplus\widehat{F}_0(-\vec{z})\oplus\widehat{F}_1(\vec{f}-\vec{y}-\vec{z})\\
(\uppsi^2F)_0&=&\widehat{F}_0\oplus \widehat{F}_1(\vec{f}-\vec{y})\oplus\widehat{F}_1(\vec{f}-\vec{z})\oplus\widehat{F}_0(\vec{f}-\vec{y}-\vec{z})\\
(\upphi F(-\vec{z})[1])_1&=&\widehat{F}_0(-\vec{z})\oplus\widehat{F}_1(\vec{f}-\vec{y}-\vec{z})\\
(\upphi F(-\vec{z})[1])_0&=&\widehat{F}_1(\vec{f}-\vec{z})\oplus\widehat{F}_0(\vec{f}-\vec{y}-\vec{z}).
\end{eqnarray*}
Since each matrix $\sigma_i$ contains  unit entries, $\sigma_j^F$ is not a zero-map.
\end{proof}
\end{lem}

\begin{lem}\label{non-zero nu}
For any non-zero object $F\in \cC_{n-2}$, there is a non-zero morphism
\[
\theta^F:\uppsi_{a_n-2}\uppsi_{a_{n-1}-2}F\to \upphi_{a_{n-1}-1}F
\]
such that the following conditions hold:
\begin{itemize}
\item[$(1)$] For any morphism $\alpha=(\alpha_1,\alpha_0):F\to F'$ between non-zero objects  in $\cC_{n-2}$, the following diagram commutes:
\[\xymatrix{
 \uppsi_{a_n-2}\uppsi_{a_{n-1}-2} F\ar[rr]^{\theta^F}\ar[d]_{\uppsi_{a_n-2}\uppsi_{a_{n-1}-2}\alpha}
 &&\upphi_{a_{n-1}-1}F\ar[d]^{\upphi_{a_{n-1}-1}\alpha}\\
\uppsi_{a_n-2}\uppsi_{a_{n-1}-1} F'\ar[rr]^{\theta^{F'}}
&&\upphi_{a_{n-1}-1}F'.
 }\]
 In particular, the morphisms $\{\theta^F\}$ define a functor morphism $\theta\colon \uppsi_{a_n-2}\uppsi_{a_{n-1}-2}\to \upphi_{a_{n-1}-1}$.
\item[$(2)$] The composition $\theta^{F'}\circ  \uppsi_{a_n-2}\uppsi_{a_{n-1}-2} \alpha$ is the zero map if and only if $\alpha$ is the zero map.
\end{itemize}
\begin{proof}
Again we only give the non-zero map $\theta^F$, and we use the similar notation as in \eqref{easy notation}. The morphism $\theta^F:\uppsi_{a_n-2}\uppsi_{a_{n-1}-2}F\to \upphi_{a_{n-1}-1}F$ is defined by 
\[
\theta^F:=\pi(-(b-2)\vec{y}-(c-2)\vec{z})[b+c-4],
\]
where 
\[\theta=(\theta_1,\theta_0):\uppsi^2\to \upphi F(-\vec{y}-\vec{z})[2]\]
is given by 
\[
\theta_1:(\uppsi^2F)_1
\xrightarrow{\left(\begin{array}{cccc}0&0&0&1\\0&-z^{c-1}&xy^{b-2}&0\end{array}\right)}
(\upphi F(-\vec{y}-\vec{z})[2])_1
\]
\[
\theta_0:(\uppsi^2F)_0
\xrightarrow{\left(\begin{array}{cccc}0&0&0&1\\0&-z^{c-1}&xy^{b-2}&0\end{array}\right)}
(\upphi F(-\vec{y}-\vec{z})[2])_0,
\]
where 
\begin{eqnarray*}
(\uppsi^2F)_1&=&\widehat{F}_1\oplus \widehat{F}_0(-\vec{y})\oplus\widehat{F}_0(-\vec{z})\oplus\widehat{F}_1(\vec{f}-\vec{y}-\vec{z})\\
(\uppsi^2F)_0&=&\widehat{F}_0\oplus \widehat{F}_1(\vec{f}-\vec{y})\oplus\widehat{F}_1(\vec{f}-\vec{z})\oplus\widehat{F}_0(\vec{f}-\vec{y}-\vec{z})\\
(\upphi F(-\vec{y}-\vec{z})[2])_1&=&\widehat{F}_1(\vec{f}-\vec{y}-\vec{z})\oplus\widehat{F}_0(\vec{f}-2\vec{y}-\vec{z})\\
(\upphi F(-\vec{y}-\vec{z})[2])_0&=&\widehat{F}_0(\vec{f}-\vec{y}-\vec{z})\oplus\widehat{F}_1(2\vec{f}-2\vec{y}-\vec{z}).
\end{eqnarray*}
Since each matrix $\theta_i$ contains  a unit entry, $\theta^F$ is not a zero-map.
\end{proof}
\end{lem}

\subsubsection{Irreducible morphisms between the exceptional objects}\label{lem irr}
In this section, using the above lemmas, we determine the irreducible morphisms between objects in $\cE^n$. For two objects $E,E'\in \cE^n$, we set 
\[
\irr(E,E'):=\dim_k\Irr(E,E'),
\]
where the vector space $\Irr(E,E')$ is defined as the following quotient space:
\[
\Irr(E,E'):=\Hom(E,E')/\Image\left(\bigoplus_{E''\in \cE^n\backslash\{E,E'\}}\Hom(E'',E')\times\Hom(E,E'')\xrightarrow{(-)\circ(-)}\Hom(E,E')\right).
\]
Then the number of arrows from a vertex $E$ to another vertex $E'$ in $Q^n$ is equal to the number $\irr(E,E')$.
By  Lemma \ref{psi part}, Lemma \ref{last lemma} and $\Hom(E^0,E^0)\cong k$, for any $E,E'\in \cE^n$ we have 
$
\hom(E,  E'):=\dim_{k}\Hom(E,E')\in \{0,1\},
$
and in particular we have 
\[
\irr(E,  E')\in \{0,1\}.
\]
The following lemmas determine the numbers $\irr(E,F)$ for all $E,E'\in \cE^{n}$.  Write 
\[\cE^{n-1}=\left(E_1,\hdots,E_{\mu_{n-1}}\right).\]

\begin{lem}\label{psi irr}
For $0\leq i\neq j\leq a_{n-1}-2$ and $E_r,E_s\in\cE^{n-1}$, we have
\[
\irr(\uppsi_iE_r, \uppsi_j E_s)=
\begin{cases}
1 & $r=s$ \mbox{ and } $j=i+1$\\
0 & \mbox{otherwise}
\end{cases}
\]
\begin{proof}
By Lemma \ref{psi part}  the vector space $\Hom(\uppsi_iE_r,\uppsi_jE_s)$  vanishes when $j\neq i+1$, so it is enough to show that 
\[
\irr(\uppsi_iE_r, \uppsi_{i+1} E_s)=
\begin{cases}
1 & r=s \\
0 & \mbox{otherwise}
\end{cases}
\]

First we prove $\irr(\uppsi_iE_r, \uppsi_{i+1} E_r)=1$. Since $\Hom(\uppsi_iE_r, \uppsi_{i+1} E_r)\cong \Hom(E_r,E_r)\cong k$, we have $\hom(\uppsi_iE_r, \uppsi_{i+1} E_r)=1$. It is enough to show that a non-zero map $\alpha:\uppsi_i E_r\to \uppsi_{i+1}E_r$ is an irreducible morphism. Assume that there are an exceptional object $E'\in \cE^n$ and non-isomorphisms
$\beta:\uppsi_i E_r\to E'$ and $\gamma: E'\to \uppsi_{i+1}E_r$ such that $\alpha=\gamma\circ \beta$. Since $\beta$ and $\gamma$ are non-zero morphisms, either of the following cases holds:
\begin{itemize}
\item[(1)] $E'\cong \uppsi_iE_t$ for some $E_t\in \cE^{n-1}$ with $r\leq t$.
\item[(2)] $E'\cong \uppsi_{i+1}E_u$ for some $E_u\in \cE^{n-1}$ with $u\leq r$.
\end{itemize}
Assume the former case $(1)$ holds. Since any non-zero endomorphism of an object in $\cE^n$ has to be an isomorphism, we see that $r<t$. Then there $\gamma$ has to be the zero-map since $\Hom(\uppsi_{i}E_t,\uppsi_{i+1}E_r)\cong \Hom(E_t,E_r)=0$. Thus the latter case (2) must hold. But since $u<r$, the morphism $\beta$ has to be the zero map, since $\Hom(\uppsi_iE_r,\uppsi_{i+1}E_u)\cong \Hom(E_r,E_u)=0$. This is a contradiction, and thus  $\alpha$ is irreducible.

We only need to prove the vanishing $\irr(\uppsi_iE_r,\uppsi_{i+1}E_s)=0$ if $r\neq s$. Since $\hom(\uppsi_iE_r,\uppsi_{i+1}E_s)=\hom(E_r,E_s)$, we may assume  $\hom(E_r,E_s)=1$. Then there is a non-zero morphism $\alpha:E_r\to E_s$. By Lemma \ref{non-zero lambda}, the composition
\[
\uppsi_iE_r\xrightarrow{\lambda_i^{E_r}}\uppsi_{i+1}E_r\xrightarrow{\uppsi_{i+1}\alpha}\uppsi_{i+1}E_s
\]
is a non-zero map. Since $\hom(\uppsi_iE_r,\uppsi_{i+1}E_s)=1$, this composition spans  $\Hom(\uppsi_iE_r,\uppsi_{i+1}E_s)$, and thus we have $\irr(\uppsi_iE_r,\uppsi_{i+1}E_s)=0$.
\end{proof}
\end{lem}

\begin{lem} For $0\leq i\leq a_{n}-2$, $0\leq j\leq a_{n-1}-1$ and $E\in\cE^{n-1}$ and $F\in \cE^{n-2}$, we have
\[
\irr(\uppsi_iE, \upphi_j F)=
\begin{cases}
 & i=a_n-2 \mbox{ and } E=\uppsi_jF\\
1&\mbox{\hspace{0mm}or}\\
& i=a_n-2,\,\, j=a_{n-1}-1\mbox{ and } E=\uppsi_{a_{n-1}-2}F\\
&\\
0 & \mbox{otherwise}
\end{cases}
\]
\begin{proof}
This follows from a similar argument as in the proof of Lemma \ref{psi irr} by using Lemma \ref{last lemma}, Lemma \ref{non-zero mu} and Lemma \ref{non-zero nu}.
\end{proof}
\end{lem}

\subsubsection{Construction of $(Q^n,I^n)$}
Now we give an inductive construction of $(Q^n,I^n)$.  We denote by $Q^n_0$ and $Q^n_1$ the set of vertices of $Q^n$ and the set of arrows in $Q^n$ respectively. Assume that the quiver with relations $(Q^k,I^k)$ is given for all $k<n$. For each $0\leq i\leq a_n-2$, we set
\begin{eqnarray*}
\uppsi_iQ_0^{n-1}&:=&\{\uppsi_iv\mid v\in Q_0^{n-1}\}\\
\uppsi_iQ_1^{n-1}&:=&\{\uppsi_i\alpha\mid \alpha\in Q^{n-1}_1\}.
\end{eqnarray*}
Then these defines a quiver 
\[
\uppsi_iQ^{n-1}:=\left(\uppsi_i{Q^{n-1}_0},\uppsi_i{Q^{n-1}_1}\right)
\]
that is isomorphic to the quiver $Q^{n-1}$. If we denote by $\uppsi_iI^{n-1}$ the relations of arrows in $\uppsi_iQ^{n-1}$ corresponding to $I^{n-1}$,  we have 
 the quiver with relations
\[
\uppsi_i(Q^{n-1},I^{n-1}):=(\uppsi_iQ^{n-1},\uppsi_iI^{n-1})
\]
that is isomorphic to $(Q^{n-1},I^{n-1})$. Similarly, for each $0\leq j\leq a_{n-1}-1$, we consider the quiver with relations
\[
\upphi_j(Q^{n-2},I^{n-2}):=(\upphi_jQ^{n-2},\upphi_jI^{n-2})
\]
that is isomorphic to $(Q^{n-2},I^{n-2})$. More precisely, the quiver $\upphi_jQ^{n-2}=(\upphi_jQ^{n-2}_0,\upphi_jQ^{n-2}_1)$  is given by 
\begin{eqnarray*}
\upphi_jQ_0^{n-2}&:=&\{\upphi_jv\mid v\in Q_0^{n-2}\}\\
\upphi_jQ_1^{n-2}&:=&\{\upphi_j\alpha\mid \alpha\in Q^{n-2}_1\},
\end{eqnarray*}
and its relations $\upphi_jI^{n-2}$  corresponds to $I^{n-2}$. 

Then the set of vertices $Q^n_0$ is defined by the following disjoint union
\[
Q^n_0:=\left(\bigsqcup_{i=0}^{a_n-2}\uppsi_iQ_0^{n-1}\right)\bigsqcup\left(\bigsqcup_{j=0}^{a_{n-1}-1}\upphi_jQ_0^{n-2}\right),
\]
and the set of arrows $Q^n_1$ is defined by 
\[
Q^n_1:=\left(\bigsqcup_{i=0}^{a_n-2}\uppsi_iQ_1^{n-1}\right)\bigsqcup\left(\bigsqcup_{j=0}^{a_{n-1}-1}\upphi_jQ_1^{n-2}\right)\bigsqcup\left(\bigsqcup_{i=0}^{a_n-3}\Lambda_i\right)\bigsqcup\left(\bigsqcup_{j=0}^{a_{n-1}-2}\Sigma_j\right)\bigsqcup\Theta,\]
where the sets $\Lambda_i$, $\Sigma_j$ and $\Theta$ are defined by
\begin{eqnarray*}
\Lambda_i&:=&\left\{ \lambda_i^v:\uppsi_1v\to \uppsi_{i+1}v\middle| v\in Q_0^{n-1}\right\}\\
\Sigma_j&:=&\left\{ \sigma_j^v:\uppsi_{a_n-2}\uppsi_jv\to \upphi_{j}v\middle| v\in Q_0^{n-2}\right\}\\
\Theta&:=&\left\{\theta^v:\uppsi_{a_n-2}\uppsi_{a_{n-1}-2}v\to \upphi_{a_{n-1}-1}v\middle| v\in Q_0^{n-2}\right\}.
\end{eqnarray*}
By setting $Q^{n-1}_{(i)}:=\uppsi_iQ^{n-1}$ and $Q^{n-2}_{[j]}:=\upphi_jQ^{n-2}$, we can draw a rough picture of the quiver $Q^n$:

\[{\small
\begin{tikzpicture}[scale=0.75]

\node (1) at (-7.5,0)  {$Q^{n-1}_{(0)}$};
\node (2) at (-4.5,0)  {$Q^{n-1}_{(1)}$};
\node (3) at (-1.5,0) [] {};
\node (4) at (0.5,0) [] {};
\node (5) at (3.5,0)  {$Q^{n-1}_{(a_n-2)}$};
\node (6) at (8,0)  {$\displaystyle \bigsqcup_{j=0}^{a_{n-1}-1}Q^{n-2}_{[j]}$};

\draw[->] (1) --  node[above] {$\lambda_0$} (2);
\draw[->] (2) -- node[above] {$\lambda_1$}  (3);
\draw[dotted] (3) --  (4);
\draw[->]  (4) -- node[above] {$\lambda_{a_n-3}$}  (5);
\draw[->] (5) --  node[above] {$\left(\sqcup_{j}\sigma_j\right)\sqcup\theta$} (6);
\end{tikzpicture}\vspace{3mm}}
\]
We can draw more precise picture as follows: To ease notation, we set $Q^{n-2}_{(i,i')}:=\uppsi_{i}\uppsi_{i'}Q^{n-2}$, $Q^{n-3}_{(i,j]}:=\uppsi_i\upphi_jQ^{n-3}$, $b_n:=a_n-2$, $c_n:=b_{n-1}=a_{n-1}-2$ and $d_n:=c_{n-1}+1=a_{n-2}-1$. Then the quiver $Q^n$ can be described by the following:

{\small
\[
\begin{tikzpicture}[scale=0.75]

\node (1) at (-7.5,0)  {$Q^{n-2}_{(0,0)}$};
\node (2) at (-4.5,0)  {$Q^{n-2}_{(1,0)}$};
\node (3) at (-1.5,0) [] {};
\node (4) at (0.5,0) [] {};
\node (5) at (3.5,0)  {$Q^{n-2}_{(b_n,0)}$};
\node (6) at (6.5,1)  {$Q^{n-2}_{[0]}$};

\draw[->] (1) --  node[above] {$\lambda_0$} (2);
\draw[->] (2) -- node[above] {$\lambda_1$}  (3);
\draw[dotted] (3) --  (4);
\draw[->]  (4) -- node[above] {$\lambda_{b_n-1}$}  (5);
\draw[->] (5) --  node[above] {$\sigma_0$} (6);

\node (1a) at (-7.5,-2)  {$Q^{n-2}_{(0,1)}$};
\node (2a) at (-4.5,-2)  {$Q^{n-2}_{(1,1)}$};
\node (3a) at (-1.5,-2) [] {};
\node (4a) at (0.5,-2) [] {};
\node (5a) at (3.5,-2)  {$Q^{n-2}_{(b_n,1)}$};
\node (6a) at (6.5,-1)  {$Q^{n-2}_{[1]}$};

\draw[->] (1a) --  node[above] {$\lambda_0$} (2a);
\draw[->] (2a) -- node[above] {$\lambda_1$}  (3a);
\draw[dotted] (3a) --  (4a);
\draw[->]  (4a) -- node[above] {$\lambda_{b_n-1}$}  (5a);
\draw[->] (5a) --  node[above] {$\sigma_1$} (6a);

\node (1b) at (-7.5,-4)  {};
\node (2b) at (-4.5,-4)  {};
\node (3b) at (-1.5,-4) [] {};
\node (4b) at (0.5,-4) [] {};
\node (5b) at (3.5,-4)  {};
\node (6b) at (6.5,-3)  {};

\node (1c) at (-7.5,-6)  {};
\node (2c) at (-4.5,-6)  {};
\node (3c) at (-1.5,-6) [] {};
\node (4c) at (0.5,-6) [] {};
\node (5c) at (3.5,-6)  {};
\node (6c) at (6.5,-5)  {};

\node (1d) at (-7.5,-8)  {$Q^{n-2}_{(0,c_{n})}$};
\node (2d) at (-4.5,-8)  {$Q^{n-2}_{(1,c_{n})}$};
\node (3d) at (-1.5,-8) [] {};
\node (4d) at (0.5,-8) [] {};
\node (5d) at (3.5,-8)  {$Q^{n-2}_{(b_n,c_{n})}$};
\node (6d) at (6.5,-7)  {$Q^{n-2}_{[c_{n}]}$};

\draw[->] (1d) --  node[above] {$\lambda_0$} (2d);
\draw[->] (2d) -- node[above] {$\lambda_1$}  (3d);
\draw[dotted] (3d) --  (4d);
\draw[->]  (4d) -- node[above] {$\lambda_{b_n-1}$}  (5d);
\draw[->] (5d) --  node[above] {$\sigma_{c_n}$} (6d);

\node (1e) at (-7.5,-10)  {$\displaystyle \bigsqcup_{j=0}^{d_n}Q^{n-3}_{(0,j]}$};
\node (2e) at (-4.5,-10)  {$\displaystyle \bigsqcup_{j=0}^{d_n}Q^{n-3}_{(1,j]}$};
\node (3e) at (-1.5,-10) [] {};
\node (4e) at (0.5,-10) [] {};
\node (5e) at (3.5,-10)  {$\displaystyle \bigsqcup_{j=0}^{d_n}Q^{n-3}_{(b_n,j]}$};
\node (6e) at (6.7,-9)  {$Q^{n-2}_{[c_{n}+1]}$};

\draw[->] (1e) --  node[above] {$\lambda_0$} (2e);
\draw[->] (2e) -- node[above] {$\lambda_1$}  (3e);
\draw[dotted] (3e) --  (4e);
\draw[->]  (4e) -- node[above] {$\lambda_{b_n-1}$}  (5e);
\draw[->] (5d) --  node[above] {$\theta$} (6e);

\draw[->] (1) -- (1a);
\draw[->] (2) -- (2a);
\draw[->] (5) -- (5a);

\draw[->] (1a) -- (1b);
\draw[->] (2a) -- (2b);
\draw[->] (5a) -- (5b);

\draw[dotted] (1b) -- (1c);
\draw[dotted] (2b) -- (2c);
\draw[dotted] (5b) -- (5c);
\draw[dotted] (6b) -- (6c);

\draw[->] (1c) -- (1d);
\draw[->] (2c) -- (2d);
\draw[->] (5c) -- (5d);

\draw[->] (1d) -- (1e);
\draw[->] (2d) -- (2e);
\draw[->] (5d) -- (5e);
\end{tikzpicture}\vspace{3mm}
\]
}

Next we construct the relations $I^n$. For this, consider the following set of relations
\begin{eqnarray*}
J^{\lambda}_{\rm null}&:=&\bigsqcup_{i=0}^{a_n-4}\left\{ \lambda_{i+1}^v\circ\lambda^v_{i}=0\,\,\middle|\,\,v\in Q_0^{n-1}\right\}\\
J^{\lambda}_{\rm comm}&:=&\bigsqcup_{i=0}^{a_n-3}\left\{ (\uppsi_{i+1}\alpha)\circ\lambda_{i}^v=\lambda_i^{v'}\circ(\uppsi_{i}\alpha)\,\,\middle|\,\,\alpha:v\to v'\in Q_1^{n-1}\right\}\\
J^{\sigma}_{\rm null}&:=&\bigsqcup_{j=0}^{a_{n-1}-1}\left\{ \sigma_{j}^{v}\circ\lambda^{\uppsi_jv}_{a_n-3}=0\,\,\middle|\,\,v\in Q_0^{n-2}\right\}\\
J^{\sigma}_{\rm comm}&:=&\bigsqcup_{j=0}^{a_{n-1}-1}\left\{ (\upphi_{j}\alpha)\circ\sigma_{j}^v=\sigma_j^{v'}\circ(\uppsi_{a_n-2}\uppsi_j\alpha)\,\,\middle|\,\,\alpha:v\to v'\in Q_1^{n-2}\right\}\\
J^{\theta}_{\rm null}&:=&\left\{ \theta^v\circ\lambda^{\uppsi_{(a_{n-1}-2)}v}_{a_n-3}=0\,\,\middle|\,\,v\in Q_0^{n-2}\right\}\\
J^{\theta}_{\rm comm}&:=&\left\{ (\upphi_{a_{n-1}-1}\alpha)\circ\theta^v=\theta^{v'}\circ(\upphi_{a_n-2}\uppsi_{a_{n-1}-2}\alpha)\,\,\middle|\,\,\alpha:v\to v'\in Q_1^{n-2}\right\}.
\end{eqnarray*}
Then the set of relations $I^n$ is generated by the following relations
\[
\left(\bigsqcup_{i=0}^{a_n-2}\uppsi_iI^{n-1}\right)\bigsqcup\left(\bigsqcup_{j=0}^{a_{n-1}-1}\upphi_jI^{n-2}\right)\bigsqcup\left(\bigsqcup_{\ast\in\{\lambda,\sigma,\theta\}}J^{\ast}_{\rm null}\right)\bigsqcup\left(\bigsqcup_{\ast\in\{\lambda,\sigma,\theta\}}J^{\ast}_{\rm comm}\right).\]

\vspace{3mm}
\begin{exa}
The quiver $Q^2$ is 
{\small
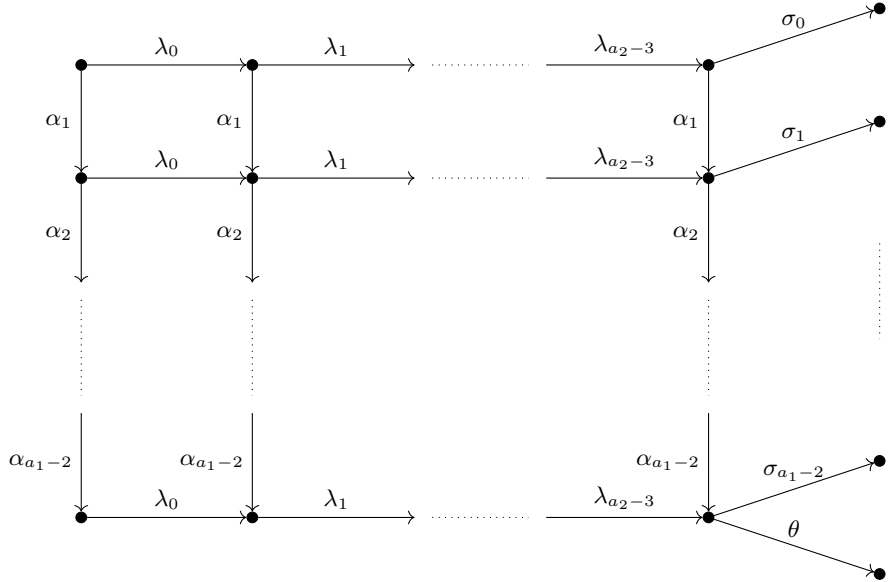
\[
\begin{tikzpicture}[scale=0.75]

\node (1) at (-7.5,0)  [DW]{};
\node (2) at (-4.5,0)  [DW]{};
\node (3) at (-1.5,0) [] {};
\node (4) at (0.5,0) [] {};
\node (5) at (3.5,0)  [DW]{};
\node (6) at (6.5,1)  [DW]{};

\draw[->] (1) --  node[above] {$\lambda_0$} (2);
\draw[->] (2) -- node[above] {$\lambda_1$}  (3);
\draw[dotted] (3) --  (4);
\draw[->]  (4) -- node[above] {$\lambda_{a_2-3}$}  (5);
\draw[->] (5) --  node[above] {$\sigma_0$} (6);

\node (1a) at (-7.5,-2) [DW]{};
\node (2a) at (-4.5,-2)  [DW]{};
\node (3a) at (-1.5,-2) [] {};
\node (4a) at (0.5,-2) [] {};
\node (5a) at (3.5,-2) [DW]{};
\node (6a) at (6.5,-1)  [DW]{};

\draw[->] (1a) --  node[above] {$\lambda_0$} (2a);
\draw[->] (2a) -- node[above] {$\lambda_1$}  (3a);
\draw[dotted] (3a) --  (4a);
\draw[->]  (4a) -- node[above] {$\lambda_{a_2-3}$}  (5a);
\draw[->] (5a) --  node[above] {$\sigma_1$} (6a);

\node (1b) at (-7.5,-4)  {};
\node (2b) at (-4.5,-4)  {};
\node (3b) at (-1.5,-4) [] {};
\node (4b) at (0.5,-4) [] {};
\node (5b) at (3.5,-4)  {};
\node (6b) at (6.5,-3)  {};

\node (1c) at (-7.5,-6)  {};
\node (2c) at (-4.5,-6)  {};
\node (3c) at (-1.5,-6) [] {};
\node (4c) at (0.5,-6) [] {};
\node (5c) at (3.5,-6)  {};
\node (6c) at (6.5,-5)  {};

\node (1d) at (-7.5,-8)  [DW]{};
\node (2d) at (-4.5,-8)  [DW]{};
\node (3d) at (-1.5,-8) [] {};
\node (4d) at (0.5,-8) [] {};
\node (5d) at (3.5,-8)  [DW]{};
\node (6d) at (6.5,-7)  [DW]{};

\draw[->] (1d) --  node[above] {$\lambda_0$} (2d);
\draw[->] (2d) -- node[above] {$\lambda_1$}  (3d);
\draw[dotted] (3d) --  (4d);
\draw[->]  (4d) -- node[above] {$\lambda_{a_2-3}$}  (5d);
\draw[->] (5d) --  node[above] {$\sigma_{a_1-2}$} (6d);

\node (6e) at (6.5,-9)  [DW]{};

\draw[->] (5d) --  node[above] {$\theta$} (6e);

\draw[->] (1) -- node[left] {$\alpha_1$} (1a);
\draw[->] (2) --  node[left] {$\alpha_1$} (2a);
\draw[->] (5) -- node[left] {$\alpha_1$} (5a);

\draw[->] (1a) -- node[left] {$\alpha_2$} (1b);
\draw[->] (2a) -- node[left] {$\alpha_2$} (2b);
\draw[->] (5a) -- node[left] {$\alpha_2$} (5b);

\draw[dotted] (1b) -- (1c);
\draw[dotted] (2b) -- (2c);
\draw[dotted] (5b) -- (5c);
\draw[dotted] (6b) -- (6c);

\draw[->] (1c) -- node[left] {$\alpha_{a_1-2}$} (1d);
\draw[->] (2c) -- node[left] {$\alpha_{a_1-2}$}  (2d);
\draw[->] (5c) -- node[left] {$\alpha_{a_1-2}$}  (5d);

\end{tikzpicture}\vspace{3mm}
\]
}
and the relations $I^2$ is given  by 
\[
I^2=\left\langle
\begin{array}{ccc}\alpha_{k+1}\alpha_k=0&\alpha_k\lambda_l=\lambda_l\alpha_k&\\
\lambda_{i+1}\lambda_i=0&\sigma_j\lambda_{a_2-3}=0&\theta\lambda_{a_2-3}=0\end{array}
\middle| 
\begin{array}{cc}
1\leq k\leq a_1-3& 0\leq l\leq a_2-3\\
0\leq i\leq a_2-4&0\leq j\leq a_1-2
\end{array}
\right\rangle.
\]
\end{exa}

\vspace{3mm}
Consider the path algebra $kQ^n$ of $Q^n$, and denote by the same notation $I^n\subset kQ^n$  the two-sided ideal of $kQ^n$ associated to the relations $I^n$ in $Q^n$. 
\begin{thm}
Let $T_n\in \cC_n$ be the tilting object in \eqref{tilting}. Then the endomorphism ring
$\End_{\cC_n}(T_n)$ is isomorphic to the algebra $kQ^n/I^n$. In particular, we have an equivalence 
\[
\cC_n\simto\Db(\fmod kQ^n/I^n).
\]
\begin{proof}
The former statement follows from the above results in  Section \ref{lem for} and Section \ref{lem irr}. We prove the latter statement.  Since $\cC_n$ is idempotent complete and algebraic, by  \cite{keller}  (see also \cite[Proposition 2.3]{ai}) or  \cite{bk}, we have an exact equivalence  
\[
\cC_n\cong \Kb(\proj kQ^n/I^n).
\]
Since the tilting object $T_n$ is a strong generator of $\cC_n$ in the sense of \cite{rou}, $\Kb(\proj kQ^n/I^n)\cong \Db(\fmod kQ^n/I^n)$ by \cite[Proposition 7.25]{rou}.
\end{proof}
\end{thm}

\section{Further applications of Theorem \ref{main result}}

\subsection{Kuznetsov--Perry's semi-orthogonal decomposition revisited}

As a special case of Corollary \ref{1d reduction}, we obtain Kuznetsov--Perry type semi-orthogonal decompositions (cf. \cite{kp}) for ramified cyclic covers of weighted projective spaces. In this subsection, we prove the following theorem.

\vspace{2mm}
 Let $n>0$, $N>1$ and $c>0$ be  positive integers. Let $W\in S:=k[x_1,\hdots,x_n]$ be a quasi-homogeneous polynomial of degree $Nc$ with respect to $\deg(x_i)=a_i\in\bZ_{>0}$, and denote by ${\pmb \upmu_N}\subset \bG_m$ the subgroup generated by a primitive $N$-th root  of unity $\zeta\in\bG_m$. Consider a ${\pmb \upmu_N}$-action and a $\bG_m$-action on $\Spec S[t]\cong\bA^n_{x}\times \bA^1_t$ defined by  $\zeta\cdot({x},t):=({x},\zeta t)$ and $g\cdot (x,t):=(g^{a_1}x_1,\hdots,g^{a_n}x_n,g^ct)$, where $x=(x_1,\hdots,x_n)$ is a coordinate of $\bA^n_x$.  Then the ${\pmb \upmu_N}$-action on $\bA^n_{x}\times \bA^1_t$ induces a ${\pmb \upmu_N}$-action on $\HMF^{\bZ}_{S[t]}(W+t^N)$. Then we have the following.

\begin{thm}\label{KP revisited}
 There are fully faithful functors 
\[\Phi_i:\HMF^{\bZ}_S(W)\hookto \HMF^{\bZ}_{S[t]}(W+t^N)^{{\pmb \upmu_N}}\] for $-N+2\leq i\leq 0$ and a semi-orthogonal decomposition
 \[  \HMF^{\bZ}_{S[t]}(W+t^N)^{{\pmb \upmu_N}}= \l \,{\sf Im}(\Phi_0),\hdots,{\sf Im}(\Phi_{-N+2})\,\r ,\]
 where $\HMF^{\bZ}_{S[t]}(W+t^N)^{{\pmb \upmu_N}}$ denotes the ${\pmb \upmu_N}$-equivariant category of $\HMF^{\bZ}_{S[t]}(W+t^N)$.
 \end{thm}

 In order to prove Theorem \ref{KP revisited}, we need to prepare some results. Set $X:=\mathbb{A}^n_x$ and $G:=\mathbb{G}_m \times {\pmb \upmu_N}$. We define characters  $\phi:\bG_m\to\bG_m$ and $\psi:\bG_m\times {\pmb \upmu_N}\to \bG_m$  by $\phi(g)=1$ and $\psi(g,\zeta):=g^c\zeta$. Fix $(a_1, \hdots, a_n) \in \mathbb{Z}_{>0}^n$. Note that  $W+t^N$ is a $\chi_{Nc}$-semi-invariant regular function on $X \times \mathbb{A}^1_t$.
By Corollary \ref{1d reduction}.(1), we have the following semi-orthogonal decomposition.

\begin{cor} \label{Kuznetsov-Perry pre} 
Assume that $N>1$. There is a fully faithful functor
\[
\Phi:\Dcoh_{\mathbb{G}_m}(\mathbb{A}^n_x,\chi_{Nc},W)\hookrightarrow\Dcoh_{\mathbb{G}_m \times{\pmb \upmu_N}}(\mathbb{A}^n_x\times\bA^1_t,\chi_{Nc} \times 1,W+t^N),\]
and we have the following semi-orthogonal decomposition:
\[
\Dcoh_{\mathbb{G}_m \times{\pmb \upmu_N}}(\mathbb{A}^n_x \times\bA^1_t,\chi_{Nc} \times 1,W+t^N)=\langle\,{\sf Im}(\Phi_0),\hdots,{\sf Im}(\Phi_{-N+2})  \,\rangle,
\]
Here ${\sf Im}(\Phi_i)$ denotes the essential image of the composition $((-)\otimes\cO(\psi^i))\circ\Phi$.
\end{cor}

We refine Corollary \ref{Kuznetsov-Perry pre} in terms of equivariant categories. For equivariant categories, see \cite{elagin2} or Appendix A. 
We prove the following proposition.
\begin{prop}\label{equivariant thm} Let $X$ be a smooth variety, and suppose that  a reductive affine algebraic $G$ and a finite group $H$ act on $X$. Let $W:X\to \bA^1$ be a non-constant semi-invariant regular function with respect to some character $\chi:G\times H\to \bG_m$. Assume that  $\Dcoh_{G\times H}(X,\chi,W)$ is idempotent complete. Then we have an equivalence
\[
\Dcoh_{G\times H}(X,\chi,W)\simto \Dcoh_G(X,\chi|_{G\times\{1\}},W)^H.
\]
\end{prop}

Before giving  the proof of Proposition \ref{equivariant thm}, we prove Lemma \ref{equivariant verdier} and Proposition \ref{equiv sing}.
In Lemma \ref{equivariant verdier}, we see the compatibility between equivariant categories and Verdier quotients.
\begin{lem}\label{equivariant verdier}
Let $\cT$ be a triangulated category with a dg-enhancement, and  let $\cS$ a thick subcategory of $\cT$. Let $G$ be a finite group acting on $\cT$. Assume that $\cS$ is stable under $G$-action, i.e. for every $g\in G$ the autoequivalences $\sigma_g:\cT\simto\cT$ defining the $G$-action preserves the subcategory $\cS$. Then we have a  fully faithful functor
\[
\Sigma:\cT^G/\cS^G\hookrightarrow (\cT/\cS)^G.
\]
Moreover, if $\cT^G/\cS^G$ is idempotent complete, the above functor is an equivalence.
\begin{proof}
Since the subcategory $\cS$ is  $G$-stable, the functor $p^*:\cT^G\to\cT$  maps the subcategory $\cS^G$ to the subcategory $\cS$, and the functor $p_*:\cT\to \cT^G$ maps $\cS$ to $\cS^G$, where the functors $p^*$ and $p_*$ are defined by Definition \ref{forget}. Hence these functors induces the functors
\[
\cT^G/\cS^G\xrightarrow{\overline{p^*}} \cT/\cS\xrightarrow{\overline{p_*}}\cT^G/\cS^G,
\]
and the adjunction $(p^*\dashv p_*)$ induces the adjunction $(\overline{p^*}\dashv\overline{p_*})$. By the argument in the proof of \cite[Proposition 3.11.(1)]{elagin2}, the adjunction morphism $\eta:{\rm id}\to p_*p^*$ is a split mono, i.e. there exists a functor morphism $\zeta:p_*p^*\to{\rm id}$ such that $\zeta\circ\eta={\rm id}$. These functor morphisms naturally induces functor morphisms $\overline{\eta}:{\rm id}\to \overline{p_*}\overline{p^*}$ and $\overline{\zeta}:\overline{p_*}\overline{p^*}\to{\rm id}$ such that $\overline{\zeta}\circ\overline{\eta}={\rm id}$. Since $\overline{\eta}$ is nothing but the adjunction morphism of the adjoint pair $(\overline{p^*}\dashv\overline{p_*})$, by Proposition \ref{comparison theorem}, the comparison functor 
\[
\Gamma:\cT^G/\cS^G\to (\cT/\cS)_{\bT(\overline{p^*},\overline{p_*})}
\]
is fully faithful. Moreover, if $\cT^G/\cS^G$ is idempotent complete, the functor $\Gamma$ is an equivalence.

On the other hand, the $G$-stability of $\cS$ also induces a $G$-action on the Verdier quotient $\cT/\cS$, and this $G$-action defines an adjoint pair $(q^*\dashv q_*)$ of functors 
\[
(\cT/\cS)^G\xrightarrow{q^*} \cT/\cS\xrightarrow{q_*}(\cT/\cS)^G.
\]
Since two comonads $\bT(\overline{p^*},\overline{p_*})$ and  $\bT(q^*,q_*)$ on $\cT/\cS$ are naturally isomorphic, we have a natural equivalence
\[
(\cT/\cS)_{\bT(\overline{p^*},\overline{p_*})}\cong(\cT/\cS)_{\bT(q^*,q_*)},
\] 
and the latter category $(\cT/\cS)_{\bT(q^*,q_*)}$ is equivalent to $(\cT/\cS)^G$ by Proposition \ref{equivariant comparison}. This finishes the proof.
\end{proof}
\end{lem}

%Here we recall singularity categories. Let $X$ be a quasi-projective scheme with an action from an affine algebraic group $G$. We denote by $\Perf_GX\subseteq \Db(\coh_GX)$ the subcategory consisting of  $G$-equivariant perfect complexes.  Then we define the {\it $G$-equivariant singularity category} $\D^{\rm sg}_G(X)$ to be the Verdier quotient $\Db(\coh_GX)/\Perf_GX$. 

\begin{prop}\label{equiv sing}
Let $Y$ be a quasi-projective scheme. Suppose that a reductive affine algebraic $G$ and a finite group $H$ act on $Y$. Assume that the $G\times H$-equivariant singularity category $\D^{\rm sg}_{G\times H}(Y)$ of $Y$ is idempotent complete. Then we have an equivalence
\[
\D^{\rm sg}_{G\times H}(Y)\cong \D^{\rm sg}_G(Y)^H.
\]
\begin{proof}
First, note that there is a natural equivalence $\coh_{G\times H}Y\cong (\coh_GY)^H$, where $(\coh_GY)^H$ denotes the $H$-equivariant category of the abelian category $\coh_GY$. By \cite[Theorem 7.1]{elagin2}, we have a natural equivalence $\Db(\coh_{G\times H}Y)\cong \Db(\coh_GY)^H$,  and it is easy to see that this equivalence restricts to an equivalence $\Perf_{G\times H}Y\cong \Perf_GY^H$.  Then we obtain the result by the following sequence of equivalences;
\begin{eqnarray}
\D^{\rm sg}_{G\times H}(Y)&=&\Db(\coh_{G\times H}Y)/\Perf_{G\times H}Y\nonumber\\
&\cong &\Db(\coh_GY)^H/ \Perf_GY^H\nonumber\\
&\cong &\bigl(\Db(\coh_GY)/\Perf_GY\bigr)^H\nonumber\\
&=&\D^{\rm sg}_G(Y)^H,\nonumber
\end{eqnarray}
where the third line follows from Lemma \ref{equivariant verdier}.
\end{proof}
\end{prop}

{\it Proof of }Proposition \ref{equivariant thm}.
Denote by $X_0$ the zero scheme of $W$. By Theorem \ref{fact sing}, we have  equivalences $\Dcoh_{G\times H}(X,\chi,W)\cong \D^{\rm sg}_{G\times H}(X_0)$ and $\Dcoh_G(X,\chi|_{G\times\{1\}},W)\cong \D^{\rm sg}_G(X_0)$. Hence  the result follows from Proposition \ref{equiv sing}.\\

By Corollary \ref{Kuznetsov-Perry pre}, Proposition \ref{equivariant thm} and Proposition \ref{coh mod}, we have Theorem \ref{KP revisited}.

\subsection{Case of Thom--Sebastiani sum}
In this subsection, we consider special cases of Theorem \ref{main result}, where  the sum $W+F$    is of Thom--Sebastiani type.
We use the notation in Section \ref{general sod section}.

Set $X:=\mathbb{A}^n_x$, $G:=\mathbb{G}_m \times {\pmb \upmu_N}$, and define characters $\phi:\bG_m\to\bG_m$ and $\psi:G\to \bG_m$ by $\phi(g):=1$ and $\psi(g,\zeta):=g\zeta$ for $g\in \bG_m$. Then $\overline{\chi}=\overline{\psi^N}=\chi_N$.  Fix $(a_1,\hdots, a_n) \in \mathbb{Z}_{>0}^n$. We define  actions from $\mathbb{G}_m$ on $\bA^n_x$ and on $\bA^m_t$ such that $\deg(x_i)=a_i$ and $\deg(t_i)=d_i$, we also consider an action from ${\pmb \upmu_N}$ on $\mathbb{A}^n_x\times \mathbb{A}^m_t$ defined by 
\[\zeta\cdot(x_1,\hdots,x_n,t_1,\hdots,t_m):=(x_1,\hdots,x_n,\zeta t_1,\hdots,\zeta t_m).\]
Assume that $W \in k[\mathbf{x}]:=k[x_1, \cdot \cdot \cdot, x_n]$ and $F \in k[\mathbf{t}]:=k[t_1,\hdots, t_m]$ are quasi-homogeneous polynomials of degree $N$ with respect to the weights $\deg(x_i)=a_i$ and $\deg(t_i)=d_i$. Then $W+F$ is a $\chi$-semi-invariant regular function on $\bA^n_x \times \mathbb{A}^m_t$, and $Z_F$ is the product $X \times V_F$, where $V_F \subset \bP (\bm d)$ is the hypersurface defined by $F$. 
 For simplicity, we write  $k[\mathbf{x},\mathbf{t}]:=k[x_1,\hdots,x_n,t_1,\hdots,t_m]$. 
By Proposition \ref{coh mod}, Theorem \ref{main result} and Proposition \ref{equivariant thm}, we have the following theorem.

\begin{thm}\label{TS sum} Set $\mu:=\sum_{i=1}^md_i$.
\begin{itemize}
\item[$(1)$] If $N<\mu$, there are fully faithful functors
\[
\Phi:\HMF^{\bZ}_{k[\mathbf{x},\mathbf{t}]}(W+F)^{\pmb \upmu_N} \hookrightarrow \Dcoh_{\mathbb{G}_m}(\bA^n_x \times V_F ,\chi_N,W),
\]
\[
\Psi_i:\HMF^{\bZ}_{k[\mathbf{x}]}(W)\hookrightarrow \Dcoh_{\mathbb{G}_m}(\bA^N_x \times V_F ,\chi_N,W)\]
for $N-\mu+1 \leq i \leq 0$
and there is a semi-orthogonal decomposition 
\[
\Dcoh_{\mathbb{G}_m}(\bA^n_x \times V_F ,\chi_n,W)=\l\, {\sf Im}(\Psi_{N-\mu+1}),\hdots,{\sf Im}(\Psi_{0}), {\sf Im}(\Phi) \, \r.
\]

\vspace{2mm}
\item[$(2)$] If $N=\mu$, we have an equivalence
\[
\Phi: \HMF^{\bZ}_{k[\mathbf{x},\mathbf{t}]}(W+F)^{\pmb \upmu_N} \simto \Dcoh_{\mathbb{G}_m}(\bA^N_x \times V_F ,\chi_N,W).
\]

\vspace{2mm}
\item[$(3)$] If $N>\mu$, 
there are fully faithful functors
\[
\Phi: \Dcoh_{\mathbb{G}_m}(\bA^n_x \times V_F ,\chi_N,W)  \hookrightarrow \HMF^{\bZ}_{k[\mathbf{x},\mathbf{t}]}(W+F)^{\pmb \upmu_N},
\]
\[
\Psi_i: \HMF^{\bZ}_{k[\mathbf{x}]}(W) \hookrightarrow \HMF^{\bZ}_{k[\mathbf{x},\mathbf{t}]}(W+F)^{\pmb \upmu_N}\]
for $\mu-N+1 \leq i \leq 0$ and there is a semi-orthogonal decomposition 
\[
\HMF^{\bZ}_{k[\mathbf{x},\mathbf{t}]}(W+F)^{\pmb \upmu_N}=\l\, {\sf Im}(\Psi_0),\hdots,{\sf Im}(\Psi_{\mu-N+1}), {\sf Im}(\Phi) \, \r.
\]

\end{itemize}
\end{thm}

%Using Morita product, we have the decomposition of the category $\Dcoh_{\bG_m}(\bA^n_x \times V_F ,\chi_N,W)$ in Theorem \ref{TS sum}.

\begin{rem}\label{Morita decomposition}
Note  that we have $(\bA^n_x \times V_F ,\chi_N,W)^{\bG_m}\cong(\bA^n_x,\chi_N,W)^{\bG_m}\boxplus(V_F,\chi_1,0)^{\bG_m}$ via the algebraic group isomorphism $\bG_m\simto \bG_m\times_{\bG_m}\times \bG_m; a\mapsto (a,\chi_N(a))$. By \cite{orlov3} and \cite[Lemma 4.8]{degree d},  the categories $\HMF^{\bZ}_{k[\mathbf{x}]}(W)$ and $\HMF^{\bZ}_{k[\mathbf{x},\mathbf{t}]}(W+F)$ in Theorem \ref{TS sum} are idempotent complete. Hence by Lemma \ref{equiv vs idem}, the category $\HMF^{\bZ}_{k[\mathbf{x},\mathbf{t}]}(W+F)^{\pmb \upmu_N}$ in Theorem \ref{TS sum} is also idempotent complete. Consequently,  by Theorem \ref{TS sum} and  \cite[Lemma 4.8]{degree d} again, $\Dcoh_{\bZ}(\bA^n_x \times V_F ,\chi_N,W)$ is also idempotent complete. Due to Proposition \ref{dg enh} and Theorem \ref{TS sum}, we have 
\begin{eqnarray}
\Dcoh_{\bG_m}(\bA^n_x \times V_F ,\chi_N,W)
&\cong& \bigl[\inj_{\bG_m}(\bA^n_x \times V_F ,\chi_N,W)\bigl]\nonumber\\
&\cong& \bigl[\inj_{\bG_m}(\bA^n_x,\chi_N,W)\circledast \inj_{\bG_m}(V_F,\chi_1,0)\bigl].\nonumber
 \end{eqnarray}
\end{rem}

The decomposition in Remark \ref{Morita decomposition} gives an application to derived categories of products of Calabi-Yau hypersurfaces. 
\begin{rem}
 Thus, if $N=n$, we have 
\begin{eqnarray}
\Dcoh_{\bG_m}(\bA^n_x \times V_F ,\chi_n,W)
&\cong& \bigl[\inj_{\bG_m}(\bA^n_x,\chi_n,W)\circledast \inj_{\bG_m}(V_F,\chi_1,0)\bigl]\nonumber\\
 &\cong& \bigl[\inj_{\bG_m}(V_{W},\chi_1,0)\circledast\inj_{\bG_m}(V_F,\chi_1,0)\bigl]\nonumber\\
&\cong& \bigl[\inj_{\bG_m}(V_{W}\times V_{F},\chi_1,0)\bigl]\nonumber\\
&\cong& \Db\bigl(\coh V_{F}\times V_W\bigr),\nonumber
\end{eqnarray}
where the first equivalence is by Remark \ref{Morita decomposition} and the second equivalence follows from Orlov's LG/CY correspondence (see \cite[Theorem 6.13]{bfk} for dg-version of it). In particular, if $n=\mu=N$, i.e. both of $V_F$ and $V_W$ are Calabi-Yau, we have an equivalence
\[
\Db\bigl(\coh V_{F}\times V_W\bigr)\cong \HMF^{\bZ}_{k[\mathbf{x},\mathbf{t}]}(W+F)^{\pmb \upmu_n}.
\]
\end{rem}

The situation in this subsection appeared in previous works \cite{bfk}, \cite{lim}.

\begin{rem}
Assume that $a_1=\hdots=a_n=d_1=\hdots=d_m=1$ and $N \geq \mathrm{max}\{n,m\}$. Then we have  $\mu=m$ and $N \geq \mu$.
In \cite{lim}, Lim studied a semi-orthogonal decomposition of the derived category $\Db([V_{W+F}/\pmb \upmu_N])$ of the quotient stack $[V_{W+F}/\pmb \upmu_N]$, where $V_{W+F} \subset \bP^{n+m}$ is the hypersurface defined by $W+F$. 
By \cite[Example 3.10]{bfk}, if $N \geq \mathrm{max}\{n, m\}$, there is a fully faithful functor  
$\HMF^{\bZ}_{k[\mathbf{x},\mathbf{t}]}(W+F)^{\pmb \upmu_N} \hookrightarrow  \Db([V_{W+F}/\pmb \upmu_N])$.
\end{rem}

\appendix

\section{Comodules over comonads and equivariant category}
In this appendix,  following \cite{elagin1,elagin2}, we recall definitions and basic properties of comodules over comonads and equivairant categories.  The only new result in this appendix is Lemma \ref{equiv vs idem}.

\subsection{Comodules over comonads}
Let $\mathcal{C}$ be a category. We begin by recalling the definitions of comonads  on $\mathcal{C}$ and comodules over a comonad.

% 2.1

\begin{dfn}
A {\it comonad} $\mathbb{T}=(T,\varepsilon,\delta)$ on  $\mathcal{C}$ consists of a functor $T:\mathcal{C}\rightarrow\mathcal{C}$ and functor morphisms $\varepsilon:T\rightarrow {\rm id}_{\mathcal{C}}$ and $\delta:T\rightarrow T^2$ such that the following diagrams are commutative$:$

\[\xymatrix{
&T \ar[r]^{\delta} \ar[d]_{\delta} \ar@{=}[dr]^{{\rm id}_T} & T^2 \ar[d]^{T\varepsilon}&&&T \ar[r]^{\delta} \ar[d]_{\delta} & T^2 \ar[d]^{T\delta}\\
&T^2 \ar[r]^{\varepsilon T} & T && &T^2 \ar[r]^{\delta T} & T^3
}\]

\end{dfn}

\begin{exa}\label{adj}
Let $P=(P^*\dashv P_*)$ be  adjoint  functors $P^*:\mathcal{C}\rightarrow\mathcal{D}$ and $P_*:\mathcal{D}\rightarrow\mathcal{C}$, and let $\eta_{P} :{\rm id}_{\mathcal{C}}\rightarrow P_*P^*$ and $\varepsilon_{P}:P^*P_*\rightarrow {\rm id}_{\mathcal{D}}$ be the adjunction morphisms. Set $T_P:=P^*P_*$ and $\delta_P:=P^*\eta_P P_*$. Then  the triple $\mathbb{T}(P):=(T_P,\varepsilon_{P},\delta_P)$ is a comonad on $\mathcal{D}$.
\end{exa}

\begin{dfn} \label{comodule}
Let $\mathbb{T}=(T,\varepsilon,\delta)$ be a comonad on $\mathcal{C}$. A {\it comodule} over $\mathbb{T}$ is a pair $(C,\theta_C)$ of an object $C\in\mathcal{C}$ and a morphism $\theta_{C}:C\rightarrow T(C)$ such that 

$(1)$ $\varepsilon(C)\circ\theta_C={\rm id}_C$, and

$(2)$  the following diagram is commutative:

$$\begin{CD}
C@>{\theta_C}>>T(C)\\
@V{\theta_C}VV @VVT({\theta_C})V\\
T(C)@>{\delta(C)}>>T^2(C).
\end{CD}$$

\end{dfn}

Given a comonad $\mathbb{T}$ on $\mathcal{C}$, we define  the category $\mathcal{C}_{\mathbb{T}}$ of comodules over the comonad $\mathbb{T}$:

\begin{dfn}
Let $\mathbb{T}=(T,\varepsilon,\delta)$ be a comonad on $\mathcal{C}$.  The category $\mathcal{C}_{\mathbb{T}}$ of comodules over $\mathbb{T}$ is the category whose  objects are comodules over $\mathbb{T}$ and  whose sets of morphisms  are defined as follows;
$${\rm Hom}_{\cC_{\bT}}((C_1,\theta_{C_1}),(C_2,\theta_{C_2})):=\{f\in\Hom_{\cC}(C_1,C_2) \mid  T(f)\circ\theta_{C_1}=\theta_{C_2}\circ f \}.$$
For a full subcategory $\mathcal{B}\subseteq\mathcal{C}$, we define the full subcategory  $\mathcal{C}_{\mathbb{T}}^{\mathcal{B}}\subseteq\mathcal{C}_{\mathbb{T}}$ as
$$\textrm{Ob}(\mathcal{C}_{\mathbb{T}}^{\mathcal{B}}):=\{(C,\theta_C)\in\textrm{Ob}(\mathcal{C}_{\mathbb{T}})\mid C\cong B  \textrm{ for some } B\in\mathcal{B}\}.$$

\end{dfn}

\begin{rem}\label{rem2.5}
Let $(C,\theta_C)\in\mathcal{C}_{\mathbb{T}}^{\mathcal{B}}$. By definition, there exist an object $B\in\mathcal{B}$ and an isomorphism $\varphi : C\xrightarrow{\sim} B$. If we set $\theta_B:=T(\varphi)\theta_C\varphi^{-1}$, then the pair $(B,\theta_B)$ is an object of $\mathcal{C}_{\mathbb{T}}^{\mathcal{B}}$ and $\varphi$ gives an isomorphism from $(C,\theta_C)$ to $(B,\theta_B)$  in $\mathcal{C}_{\mathbb{T}}^{\mathcal{B}}$.

\end{rem}

For a comonad which is given by an adjoint pair $(P^*\dashv P_*)$, we have a canonical functor, called comparison functor, from  the domain of $P^*$ to the category of comodules over the comonad.

\begin{dfn}\label{comparison functor}
Notation is same as in Example \ref{adj}. For an adjoint pair $P=(P^*\dashv P_*)$, we define a  functor $$\Gamma_{P}:\mathcal{C}\rightarrow\mathcal{D}_{\mathbb{T}(P)}$$ as follows: For any $C\in\mathcal{C}$ we define $\Gamma_{P}(C):=(P^*(C),P^*(\eta_{P}(C)))$, and for any  morphism $f$ in $\mathcal{C}$ we define $ \Gamma_{P}(f):=P^*(f).$
This functor is called the {\it  comparison functor} of $P$. 
\end{dfn}

The following proposition gives  sufficient conditions for a comparison functor to be fully faithful or an equivalence.

\begin{prop}[{\cite[Theorem 3.9, Corollary 3.11]{elagin1}}]\label{comparison theorem}
Notation is same as above.
\begin{itemize}

\item[$(1)$] If for any $C\in\mathcal{C}$, the morphism $\eta_P(C):C\rightarrow P_*P^*(C)$ is a split mono, i.e. there is a  morphism $\zeta_C:P_*P^*(C)\rightarrow C$ such that $\zeta\circ\eta_P(C)={\rm id}_{{C}}$, then the comparison functor $\Gamma_{P}:\mathcal{C}\rightarrow\mathcal{D}_{\mathbb{T}(P)}$ is fully faithful. 
\item[$(2)$] If  $\mathcal{C}$ is idempotent complete and the functor morphism $\eta_P :{\rm id}_{\mathcal{C}}\rightarrow P_*P^*$ is a split mono, i.e. there exists a functor morphism $\zeta:P_*P^*\rightarrow {\rm id}_{\mathcal{C}}$ such that $\zeta\circ\eta={\rm id}$, then $\Gamma_{P}:\mathcal{C}\rightarrow\mathcal{D}_{\mathbb{T}(P)}$ is an equivalence.

\end{itemize}
\end{prop}

Next we recall linearizable functors which induce natural functors between categories of comodules following \cite{H1}.
Let $\mathcal{A}$ (resp. $\mathcal{B}$) be a category and  $\mathbb{T}_{\mathcal{A}}=(T_{\mathcal{A}},\varepsilon_{\mathcal{A}},\delta_{\mathcal{A}})$ (resp. $\mathbb{T}_\mathcal{B}=(T_{\mathcal{B}},\varepsilon_{\mathcal{B}},\delta_{\mathcal{B}})$) a command on $\mathcal{A}$ (resp.   $\mathcal{B}$).

\begin{dfn}\label{equivariant}A functor $F:\mathcal{A}\rightarrow\mathcal{B}$ is said to be {\it linearizable} with respect to $\mathbb{T}_{\mathcal{A}}$ and $\mathbb{T}_{\mathcal{B}}$, or just {\it linearizable}, if there exists an isomorphism of functors $$\Omega:F  T_{\mathcal{A}}\xrightarrow{\sim}T_{\mathcal{B}}  F$$ such that the following two diagrams of functor morphisms are commutative :
\[\xymatrix{
(1)& FT_{\mathcal{A}} \ar[rr]^{\Omega} \ar[rd]_{F\varepsilon_{\mathcal{A}}} & &  T_{\mathcal{B}}F \ar[ld]^{\varepsilon_{\mathcal{B}}F}&&(2)&FT_{\mathcal{A}} \ar[rrr]^{\Omega} \ar[d]_{F\delta_{\mathcal{A}}} &&&   T_{\mathcal{B}}F \ar[d]^{\delta_{\mathcal{B}}F} \\ && F &&&&FT^2_{\mathcal{A}} \ar[rrr]^{T_{\mathcal{B}}\Omega\circ\Omega T_{\mathcal{A}}} & && T^2_{\mathcal{B}}F  \\
}\]
\vspace*{2mm}We call  the pair $(F,\Omega)$  a {\it linearized functor} with respect to $\mathbb{T}_{\mathcal{A}}$ and $\mathbb{T}_{\mathcal{B}}$, and the isomorphism $\Omega:F  T_{\mathcal{A}}\xrightarrow{\sim}T_{\mathcal{B}}  F$ of functors  is called a {\it linearization} of $F$ with respect to $\mathbb{T}_{\mathcal{A}}$ and $\mathbb{T}_{\mathcal{B}}$.
\end{dfn}

If $F:\mathcal{A}\rightarrow\mathcal{B}$ is a linearizable functor with a linearization $\Omega:F T_{\mathcal{A}}\xrightarrow{\sim}T_{\mathcal{B}} F$, we have an induced functor 
$${F_{\Omega}}:\mathcal{A}_{\mathbb{T}_{\mathcal{A}}}\rightarrow\mathcal{B}_{\mathbb{T}_{\mathcal{B}}}$$
defined by $$F_{\Omega}(A,\theta_A):=(F(A),\Omega(A)\circ F(\theta_A))\hspace{3mm}\textrm{and}\hspace{3mm}F_{\Omega}(f):=F(f).$$

\vspace{2mm}
The following proposition is a special case of \cite[Proposition 2.10]{H1}.

\begin{prop}[cf. {\cite[Proposition 2.10]{H1}}]\label{main prop 2}
Let $F:\mathcal{A}\rightarrow\mathcal{B}$ be a linearizable functor with a linearization $\Omega:F T_{\mathcal{A}}\xrightarrow{\sim}T_{\mathcal{B}} F$. If $F:\cA\to \cB$ is  fully faithful (resp. an equivalence), then the induced functor $F_{\Omega}:\cA_{\bT_{\cA}}\to\cB_{\bT_{\cB}}$ is also fully faithful (resp. an equivalence).
\end{prop}

\vspace{3mm}
\subsection{Equivariant category}
In this section all categories are $k$-linear. Let $\cC$ be an additive category and  $G$
 a finite group.

\begin{dfn}
A {\it (right) action} of $G$ on $\cC$ is given by the following data (i) and (ii):
\begin{itemize}
\item[(i)] For every $g\in G$, an autoequivalence $\sigma_g:\cC\simto\cC$.
\item[(ii)] For every $g,h\in G$, functor isomorphisms $c_{g,h}:\sigma_g\circ\sigma_h\simto\sigma_{hg}$ such that the diagram
\[\begin{tikzcd}
\sigma_f\sigma_g\sigma_h\arrow[rr, "{c_{g,h}}"]\arrow[d,"c_{f,g}"']&&\sigma_f\sigma_{hg}\arrow[d, "c_{f,gh}"]\\
\sigma_{gf}\sigma_h\arrow[rr, "c_{gf,h}"]&&\sigma_{hgf}\,
\end{tikzcd}\]
commutes for all $f,g,h\in G$.
\end{itemize}
\end{dfn} 

For a $G$-action on $\cC$, the {\it equivariant category} $\cC^G$ of $\cC$ is defined as follows: An object of $\cC^G$ is a  pair $(F,(\uptheta_g)_{g\in G})$ of an object $F\in \cC$ and $(\uptheta_g)_{g\in G}$ is a family of isomorphisms
\[
\uptheta_g:F\simto\sigma_g(F)
\]
such that the diagram
\[\begin{tikzcd}
F\arrow[rr, "{\uptheta_g}"]\arrow[d,"\uptheta_{hg}"']&&\sigma_g(F)\arrow[d, "\sigma_g(\uptheta_h)"]\\
\sigma_{hg}(F)&&\sigma_{g}(\sigma_h(F))\arrow[ll, "c_{g,h}"]\,
\end{tikzcd}\]
commutes for all $g,h\in G$. A morphism $f:(F_1,(\uptheta^1_g))\to(F_2,(\uptheta^2_g))$ in $\cC^G$ is defined as a morphism $f:F_1\to F_2$ in $\cC$ such that $f$ is compatible with $\uptheta^i$.

\begin{rem}
If a category $\cC$ has a (right) $G$-action, then for any two objects $(F_1,(\uptheta_g^1))$  and $(F_2,(\uptheta^2_g))$, the the equivariant structures $(\uptheta^1_g)$ and $(\uptheta_g^2)$ defines  a natural (right) $G$-action on the set $\Hom_{\cC}(F_1,F_2)$ as follows: For any $\varphi\in \Hom_{\cC}(F_1,F_2)$ and $g\in G$, the morphism $\varphi\cdot g\in \Hom_{\cC}(F_1,F_2)$ is defined as the following composition of morphisms;
\[
F_1\xrightarrow{\uptheta^1_g}\sigma_g(F_1)\xrightarrow{\sigma_g(\varphi)}\sigma_g(F_2)\xrightarrow{(\uptheta_g^2)^{-1}}F_2.
\]
By definition, the set of morphisms in the equivariant category $\cC^G$ is equal to the $G$-invariant subspace of the set of morphisms in $\cC$;
\[\Hom_{\cC^G}((F_1,(\uptheta_g^1)),(F_2,(\uptheta_g^2)))=\Hom_{\cC}(F_1,F_2)^G.\]
\end{rem}

\begin{dfn}\label{forget} Define a functor 
\[
p^*:\cC^G\to \cC
\]
to be the forgetful functor, and define a functor \[p_*:\cC\to \cC^G\] by 
$p_*(F):=\Bigl(\bigoplus_{h\in G}\sigma_h(F), (\uptheta_g)\Bigl)$.
\end{dfn}

Then by \cite[Lemma 3.8]{elagin2} we have the  adjunctions
\[
p_*\dashv p^*\dashv p_*.
\]
Recall that $\bT(p^*,p_*)$ denotes the comonad induced by the adjunction $p^*\dashv p_*$.

\begin{prop}[{\cite[Proposition 3.11]{elagin2}}]\label{equivariant comparison}
The comparison functor
\[\cC^G\to\cC_{\bT(p^*,p_*)}\]
is an equivalence.

\end{prop}

Note that any $G$-action on $\cC$ naturally extends to the $G$-action on its idempotent completion $\widehat{\cC}$ (see \cite[Proposition 3.13]{elagin2}).

\begin{lem}\label{equiv vs idem}
Let $\cC$ be an additive category with an action of a finite group $G$. Then there is an additive equivalence
\[
\widehat{\cC^G}\cong(\widehat{\cC}\,)^G.
\]
In particular, if $\cC$ is idempotent complete,  so is $\cC^G$. Furthermore, if $G$ is abelian, $\cC^G$ is idempotent complete if and only if $\cC$ is idempotent complete.

\end{lem}

Before we prove this, we recall basic properties of idempotent completion.
An additive functor $F:\cA\to \cB$   between additive categories induces an additive functor $\widehat{F}:\widehat{\cA}\to \widehat{\cB}$ defined by $\widehat{F}(A,e):=(F(A),F(e))$. If $F_1,F_2:\cA\to \cB$ are additive functors and $\varphi:F_1\to F_2$ is a functor morphism, then $\varphi$ induces a functor morphism $\widehat{\varphi}:\widehat{F_1}\to \widehat{F_2}$ defined by
\[
\widehat{\varphi}(A,e):=F_2(e)\circ\varphi(A)\circ F_1(e).
\]
Since  $\varphi(A)\circ F_1(e)=F_2(e)\circ\varphi(A)$, we have  $\widehat{\varphi}(A,e)=F_2(e)\circ\varphi(A)=\varphi(A)\circ F_1(e)$. This implies that the equality $\widehat{\varphi_2\circ\varphi_1}=\widehat{\varphi_2}\,\circ\,\widehat{\varphi_1}$ of functor morphisms, and for any additive functor $F:\cA\to \cB$, we have $\widehat{{\rm id}_{F}}={\rm id}_{\widehat{F}}$ since the identity morphism  of an object $(A,e)\in \widehat{\cA}$ is the idempotent $e:A\to A$.
Therefore, we see that if $P^*:\cA\to \cB$ and $P_*:\cB\to \cA$ are an adjoint pair, the induced functors $\widehat{P^*}:\widehat{\cA}\to \widehat{\cB}$ and $\widehat{P_*}:\widehat{\cB}\to \widehat{\cA}$ are also an adjoint pair.

\vspace{2mm}

{\it Proof of Lemma \ref{equiv vs idem}.} The latter statement follows from the former one by \cite[Theorem 4.2]{elagin2}. Let $p^*:\cC^G\to \cC, p_*:\cC\to \cC^G$ and $q^*:(\widehat{\cC})^G\to\widehat{\cC}, q_*:\widehat{\cC}\to (\widehat{\cC})^G$ be the adjoint pairs defined in Definition \ref{forget}. Then the adjoint pair $(p^*\dashv p_*)$ induces the adjoint pair $(\widehat{p^*}\dashv \widehat{p_*})$. Then we obtain two comonads   $\bT(q^*,q_*)$ and  $\bT(\widehat{p^*},\widehat{p_*})$ on $\widehat{C}$, and these comonads are tautologically isomorphic to each other. Hence by Proposition \ref{equivariant comparison}  there is a sequence of  equivalences
\[
(\widehat{\cC})^G\cong \widehat{\cC}_{\bT(q^*,q_*)}\cong \widehat{\cC}_{\bT(\widehat{p^*},\widehat{p_*})}.
\]
Thus we only need to show the comparison functor \[\widehat{\cC^G}\to \widehat{\cC}_{\bT(\widehat{p^*},\widehat{p_*})}\] is an equivalence. For this, we use Proposition \ref{comparison theorem}.(2). Since $\widehat{\cC^G}$ is  idempotent complete, it suffices to verify that the functor morphism $\widehat{\eta}:{\rm id}\to \widehat{p_*}\widehat{p^*}$ is a split mono. By the proof of  \cite[Proposition 3.11.(1)]{elagin2}, the functor morphism $\eta:{\rm id }\to p_*p^*$ splits, and so there exists a functor morphism $\gamma:p_*p^*\to {\rm id }$ such that $\gamma\circ\eta={\rm id}$. It is easy to see that the induced functor morphism $\widehat{\gamma}:\widehat{p_*}\widehat{p^*}\to {\rm id}$ is a retraction of $\widehat{\eta}$. \qed

\vspace{2mm}
Let $\cT$ be a triangulated category with an action of  a finite group $G$, where the autoequivalences $\sigma_g:\cT\simto \cT$ of the action are supposed to be exact equivalences.
Then the equivariant category $\cT^G$ has  natural shift functors and a class of distinguished triangles induced by the triangulated structure of $\cT$. The following result guarantees that, if $\cT$ has a dg-enhancement,  $\cT^G$ is a triangulated category with respect to the natural shift functors and distinguished triangles.

\begin{prop}[{\cite[Corollary 6.10]{elagin2}}]
Notation is same as above. If  $\cT$ has a dg-enhancement, then  $\cT^G$ is a triangulated category  with respect to the natural shift functors and distinguished triangles.
\end{prop}

\vspace{10mm}

\end{document}